\magnification=\magstep1
\input amstex
\documentstyle{amsppt}
\pagewidth{6.5truein}
\pageheight{9truein}
\loadbold

\def\R{{\bold R}}
\def\Q{{\bold Q}}

\def\N{{\bold N}}
\def\Z{{\bold Z}}

\def\Hom{\mathop{\roman{Hom}}\nolimits}
\def\st{\mathop{\roman{st}}\nolimits}
\def\lk{\mathop{\roman{lk}}\nolimits}
\def\Iso{\mathop{\roman{Iso}}\nolimits}
\def\Aiso{\mathop{\roman{Aiso}}\nolimits}
\def\id{\mathop{\roman{id}}\nolimits}
\def\dim{\mathop{\roman{dim}}\nolimits}
\def\Int{\mathop{\roman{Int}}\nolimits}
\def\Im{\mathop{\roman{Im}}\nolimits}
\def\Sing{\mathop{\roman{Sing}}\nolimits}
\def\graph{\mathop{\roman{graph}}\nolimits}
\def\codim{\mathop{\roman{codim}}\nolimits}
\def\Dom{\mathop{\roman{Dom}}\nolimits}

\topmatter
\title PL and differential topology in o-minimal structure\endtitle 
\author Masahiro SHIOTA \endauthor 
\address 
Graduate School of Mathematics, Nagoya University, Chikusa, Nagoya, 
464-8602, Japan\endaddress
\abstract
Arguments on PL\,(=piecewise linear) topology work over any ordered field in the same way as over $\R$, and those on 
differential topology do over a real closed field $R$ in an o-minimal structure that expands $(R,<,0,1,+,\cdot)$. 
One of the most fundamental properties of definable sets is that a compact definable set in $R^n$ is definably homeomorphic to a polyhedron (see [v]). 
We show uniqueness of the polyhedron up to PL homeomorphisms (o-minimal Hauptvermutung). 
Hence a compact definable topological manifold admits uniquely a PL manifold structure and is, so to say, tame. 
We also see that many problems on PL and differential topology over $R$ can be translated to those over $\R$. 
\endabstract 
\subjclass 
03C64, 57Q05, 57Q25
\endsubjclass
\keywords
PL topology, differential topology, o-minimal structure 
\endkeywords
\thanks
This research was partially supported by Fields Institute 
\endthanks
\endtopmatter

\document
The main theorem of the present paper is the theorem 2.1 of o-minimal Hauptvermutung in \S 2.
\head \S 1. PL topology over an ordered field \endhead
  Let $R$ denote an ordered field. 
For simplicity of notation we assume $R\supset\R$, but the arguments in the other case will be clear by the context. 
In this section we show that the known results on PL topology hold over $R$ (remark 1.3). 
A key to proof of o-minimal Hauptvermutung is the idea of a regular neighborhood and the Alexander trick in PL topology (see [R-S]). 
Hence we need to clarify PL topology over $R$. 
This is the other purpose of this section. 
A {\it simplex} in $R^n$, $n\in\N$, is defined as in $\R^n$ (we set $\N=\{0,1,...\}$). 
A {\it polyhedron} in $R^n$ is a subset of $R^n$ which is locally a finite union of simplexes. 
A {\it compact} polyhedron is a finite union of simplexes. 
(Note that a simplex is, in general, not compact in the usual sense when we give a topology to $R$ by open intervals.) 
All other terminology of PL topology (e.g., in [R-S]) is defined in the same way as in the real case. 
Here we need to be careful not to use topological but non-PL terminology (see the definition of a PL homotopy below). 
For a simplicial complex $K$, $K^r$ denotes the $r$-skeleton of $K$ (p\. 15 in [R-S]), and $|K|$ does the underlying 
polyhedron to $K$ (p\. 14 in [R-S]). 
A {\it full} subcomplex $L$ of $K$ is a subcomplex such that each $\sigma\in K$ with vertices all in $L$ is a simplex in $L$ 
(see p.~ 30 in [R-S] for properties of a full subcomplex). 
From now the fundamentals of PL topology in [R-S] are assumed to be known. \par
¡¡¡¡A simplex and, hence, a polyhedron $X$ in $\R^n$ are naturally extended to a simplex and a polyhedron in $R^n$, denoted by 
$X_R$ and called the {\it$R$-extension} of $X$. 
Considering graph we extend an $\R$-PL map $f:X\to Y$ between polyhedra in $\R^n$ to an $R$-PL map $f_R:X_R\to Y_R$, which 
we call the $R$-{\it extension} of $f$. 
For a simplicial complex $K$ in $\R^n$, let $K_R$ denote the extended simplicial complex in $R^n$ $\{\sigma_R:\sigma\in 
K\}$. 
We also define the extended simplicial map $\phi_R:K_R\to L_R$ of a simplicial map $\phi:K\to L$ between simplicial 
complexes in $\R^n$. 
Conversely, for simplicial maps $\phi_i:K_i\to L_i$ between {\bf finite} or {\bf countable} simplicial complexes in $R^n$ 
there exist $\R$-simplicial maps $\psi_i:M_i\to N_i$ between simplicial complexes in $\R^n$ and $R$-simplicial 
isomorphisms $\pi_i:K_i\to M_{iR}$ and $\tau_i:L_i\to N_{iR}$ such that $\tau_i\circ\phi_i=\psi_{iR}\circ\pi_i$ and if 
some of $K_i$'s and $L_i$'s coincide then the corresponding ones $(M_i,\pi_i)$ and $(N_i,\tau_i)$ coincide because the 
category of finite or countable simplicial complexes and simplicial maps is combinatorial and does not depend on the 
choice of $R$. 
We reduce compact polyhedra in $R^n$ and PL maps between them to those in $\R^n$ as follows. 
\proclaim{Lemma 1.1}
Let $f_i:X_i\to Y_i$ be a finite number of PL maps between compact polyhedra in $R^n$. 
Here there may be distinct $i$ and $j$ with $X_i=Y_j$ or $Y_i=Y_j$. 
Assume that $X_i\not=X_j$ for $i\not=j$ and if there is a sequence $X_{i_1}\overset f_{i_1}\to\longrightarrow Y_{i_1}=
X_{i_2}\overset f_{i_2}\to\longrightarrow\cdots\overset f_{i_k}\to\longrightarrow Y_{i_k}$ then $X_{i_1}\not=Y_{i_k}$. 
Then there exist $\R$-PL maps $g_i:U_i\to V_i$ between polyhedra in $\R^m$ for some $m\in\N$ and $R$-PL homeomorphisms 
$\pi_i:X_i\to U_{iR}$ and $\tau_i:Y_i\to V_{iR}$ such that $\tau_i\circ f_i=g_{iR}\circ\pi_i$, if $X_i=Y_j$ then $U_i=
V_j$ and $\pi_i=\tau_j$, and if $Y_i=Y_j$ then $V_i=V_j$ and $\tau_i=\tau_j$. \par
  The relative case $f_i:(X_i,X_{i,j})_j\to(Y_i,Y_{i,j})_j$ holds. 
To be precise, let $X_{i,j}$ and $Y_{i,j}$ be a finite number of unions of a finite number of open simplexes included in 
$X_i$ and $Y_i$ such that $f_i(X_{i,j})\subset Y_{i,j}$. 
Then there exist unions of a finite number of open simplexes $U_{i,j}$ of $U_i$ and $V_{i,j}$ of $V_i$ such that $\pi_i(
X_{i,j})=U_{i,jR}$ and $\tau_i(Y_{i,j})=V_{i,jR}$. 
\endproclaim
\demo{Proof}
Theorem 2.15 in [R-S] states that in the real number case there exist simplicial decompositions $K_i$ of $X_i$ and $L_i$ of 
$Y_i$ such that $f_i:K_i\to L_i$ are simplicial, $K_i=L_j$ if $X_i=Y_j$, and $L_i=L_j$ if $Y_i=Y_j$. 
The $R$-case is proved in the same way.
In the relative case we can choose $K_i$ and $L_i$ so that each $X_{i,j}$, or $Y_{i,j}$, is the unions of some open simplexes 
in $K_i$, or $L_i$ respectively. 
Regard $f_i:K_i\to L_i$ combinatorially and abstractly, and realize them in some $\R^m$. 
Then the lemma follows. \qed
\enddemo
The assumption in lemma 1.1 cannot be dropped. 
Let $X$ be a compact polyhedron in $R^n$ of positive dimension, consider all $R$-PL maps $f:X\to X$, and define an 
equivalence relation on them by setting $f_1\thicksim f_2$ if there exists a PL homeomorphism $\pi$ of $X$ with $f_1
\circ\pi=\pi\circ f_2$ (i.e., if $f_1$ and $f_2$ are PL conjugate). 
Then the cardinal number of the equivalence classes is $\# R$. 
Hence for $X$ in $\R^n$, the equivalence classes of $\R$-PL maps between $X$ do not coincide with the equivalence 
classes of $R$-PL maps between $X_R$ if $\# R>\#\R$. \par
  This is, however, the case for homotopy classes as follows. 
Since the simplicial approximation theorem (see [H]) does not hold for general $R$ in the above-mentioned topology, we 
clarify the definition of a homotopy. 
Let $X$ and $Y$ be compact polyhedra in $R^n$. 
$R$-PL maps $f,\,g:X\to Y$ are {\it $R$-PL homotopic} ({\it$R$-PL isotopic}) if there exists an $R$-PL map $F:X\times
[0,\,1]\to Y$ such that $F(\cdot,0)=f$ and $F(\cdot,1)=g$ (and $F(\cdot,t)$ is a PL imbedding for each $t\in[0,\,1]$). 
(We write $F(\cdot,t)$ as $f_t(\cdot)$, $0\le t\le 1,$ and call it a {\it homotopy} of, or from, $f_0$ to $f_1$.) 
Let $\Hom_R(X,Y)\ (\Iso_R(X,Y))$ denote the $R$-PL homotopy (isotopy) classes. 
Let $X_1\subset X$ be a compact subpolyhedron and $f:X_1\to Y$ a PL map. 
Let $\Hom_R(X,Y)_f\ (\Iso_R(X,Y)_f)$ denote the equivalence classes by $R$-PL homotopies (isotopies) $F:X\times[0,\,1]
\to Y$ such that $F(\cdot,t)=f$ on $X_1$ for each $t\in[0,\,1]$. 
An {\it$R$-PL isotopy of} $X$ is a PL map $F:X\times[0,\,1]\to X$ such that $F(\cdot,0)=\id$ and $F(\cdot,t)$ is a 
homeomorphism of $X$ for each $t\in[0,\,1]$. 
(Note an $R$-PL isotopy $F:X\times[0,\,1]\to X$ of the identity map does not mean an $R$-PL isotopy of $X$. 
In the former case $F(\cdot,t)$ is not necessarily a homeomorphism of $X$. 
In order to distinguish them we call an $R$-PL map $F:X\times[0,\,1]\to Y$ an {\it isotopy through homeomorphisms} if each 
$F(\cdot,t)$ is a homeomorphism from $X$ to $Y$ and $R$-PL maps $f,g:X\to Y$ {\it isotopic through homeomorphisms} if there 
exists an isotopy of $f$ to $g$ through homeomorphisms.) 
We define the {\it ambient $R$-PL isotopy classes} $\Aiso_R(X,Y)$ ($\Aiso_R(X,Y)_f$) by $R$-PL isotopies of $Y$ 
(fixing $f(X_1)$, respectively) as at p\. 37 in [R-S]. 
Then 
\proclaim{Lemma 1.2}
Let $X_1$ be a compact subpolyhedron of $X$ and $f:X_1\to Y$ a PL map. 
Then the following former four natural maps are bijective, and the latter two ones are injective. 
$$
\gather
\Hom_\R(X,Y)\to\Hom_R(X_R,Y_R),\ \Hom_\R(X,Y)_f\to\Hom_R(X_R,Y_R)_{f_R},\\
\Iso_\R(X,Y)\to\Iso_R(X_R,Y_R),\ \Iso_\R(X,Y)_f\to\Iso_R(X_R,Y_R)_{f_R},\\
\Aiso_\R(X,Y)\to\Aiso_R(X_R,Y_R),\ \Aiso_\R(X,Y)_f\to\Aiso_R(X_R,Y_R)_{f_R}. 
\endgather
$$
The relative case also holds as in lemma 1.1. 
\endproclaim
Here the maps $\Aiso_\R(X,Y)\to\Aiso_R(X_R,Y_R)$ and $\Aiso_\R(X,Y)_f\to\Aiso_R\allowmathbreak(X_R,Y_R)_{f_R}$ are not 
necessarily surjective. 
Let $X=[0,\,1]$ and $Y$ a polyhedron in $\R^2$ of the form of the alphabet Y. 
Let $f$ be a PL imbedding of $X_R$ into $Y_R$ carrying a number of $(0,\,1)_R-X$ to the singular point of $Y$. 
Then the ambient $R$-PL isotopy class of $f$ is not the $R$-extension of any ambient $\R$-PL isotopy class. 
\demo{Proof}
{\it Proof of surjectivity of the map $\Hom_\R(X,Y)\to\Hom_R(X_R,Y_R)$.}
Let $K$ and $L$ be $\R$-simplicial decompositions of $X$ and $Y$, respectively, and $g:X_R\to Y_R$ an $R$-PL map. 
We show that $g$ is $R$-PL homotopic to the $R$-extension of some $\R$-PL map from $X$ to $Y$. 
By Theorem 2.15 in [R-S] there exist $R$-simplicial subdivisions $M$ of $K_R$ and $N$ of $L_R$ such that $g:M\to N$ is 
simplicial. 
Hence it suffices to prove the following statement. \par
{\it Statement.} 
There exists an $R$-PL isotopy $\Phi:X_R\times[0,\,1]\to X_R$ of $X_R$ preserving $K_R$ (i.e., $\Phi(\sigma,t)\subset\sigma$ 
for each $\sigma\in K_R$ and $t\in[0,\,1]$) and whose finishing homeomorphism $\Phi(\cdot,1)$ is a simplicial 
isomorphism from $M$ to the $R$-extension of some $\R$-simplicial subdivision $K'$ of $K$. \par
  {\it Proof of the statement.}
Assume there exist non-real vertices in $M$ (i.e., vertices outside of $|K|$). 
It suffices to decrease the number of non-real vertices. 
For each vertex $v$ in $M$, let $\tau(v)\in K_R$ such that $v\in\Int\tau(v)$. 
There is at least one non-real vertex $v$ in $M$ such that $\Int|\st(v,M)|\cap\tau(v)$ includes a real point, say $v'$, 
where $\st(v,M)$ and $\Int|\st(v,M)|$ denote the simplicial complex generated by $\sigma\in M$ with $v\in\sigma$ and $\cup\{
\Int\sigma:\sigma\in\st(v,M),v\in\sigma\}$ respectively. 
Consider a simplicial subdivision $C$ of the cell complex $\{\sigma\times\{0\},\sigma\times\{1\},\sigma\times[0,\,1]:
\sigma\in\st(v,M)\}$ without introducing new vertices (the $R$-case of Proposition 2.9 in [R-S], which is proved in the 
same way). 
(A cell means that in the terminology of PL topology.) 
Define a simplicial complex $C'$ in $R^n$ and a simplicial isomorphism $\Psi_v:C\to C'$ so that $\Psi_v(\cdot,0)=\id,
\ \Psi_v(v,1)=(v',1)$ and $\Psi_v(\cdot,1)=\id$ on $\lk(v,M)$, where $\lk(v,M)=\{\sigma\in\st(v,M):v\not\in\sigma\}$. 
Then $|C|=|C'|$, $C'$ is a new simplicial decomposition of $|\st(v,M)|\times[0,\,1]$ with vertex $(v',1)$ and $\Psi_v=
\id$ on $|\lk(v,M)|\times[0,\,1]$. 
Extend identically the $R$-PL homeomorphism $\Psi_v:|C|\to|C|$ to $X_R\times[0,\,1]\to X_R\times[0,\,1]$ and compose 
the extension and the projection $X_R\times[0,\,1]\to X_R$. 
Then we have an $R$-PL isotopy $\Phi_v:X_R\times[0,\,1]\to X_R$ of $X_R$ preserving $K_R$, fixing $X_R-\Int|\st(v,M)|$ 
and whose finishing homeomorphism carries $v$ to $v'$ and is a simplicial isomorphism from $M$ to some simplicial 
subdivision of $K_R$. 
Thus $v$ becomes a real vertex and the number of non-real vertices decreases. \par
  We strengthen the statement as follows. 
Let $M_1$ be a subcomplex of $M$ which is the $R$-extension of some $\R$-simplicial complex. 
Then we can choose $\Phi$ so that $\Phi(\cdot,t)=\id$ on $|M_1|$ for each $t\in[0,\,1]$. 
This is clear by the above proof. \par
  {\it Proof of injectivity of $\Hom_\R(X,Y)\to\Hom_R(X_R,Y_R)$.}
Let $g,\,h:X\to Y$ be $\R$-PL maps such that $g_R$ and $h_R$ are $R$-PL homotopic. 
Let $F:X_R\times[0,\,1]\to Y_R\times[0,\,1]$ be an $R$-PL map such that $F(\cdot,0)=(g_R,0)$, $F(\cdot,1)=(h_R,1)$ and 
$F(X_R\times\{t\})\subset Y_R\times\{t\}$ for each $t\in[0,\,1]$. 
Let $K$ and $L$ be $\R$-simplicial decompositions of $X\times\{0,1\}$ and $Y\times\{0,1\}$, respectively, such that $F|
_{X\times\{0,1\}}:K\to L$ is simplicial. 
Let $M$ and $N$ be $R$-simplicial decompositions of $X_R\times[0,\,1]$ and $Y_R\times[0,\,1]$, respectively, such that 
$F:M\to N$ is simplicial and $M|_{X_R\times\{0,1\}}$ and $N|_{Y_R\times\{0,1\}}$ are subdivisions of $K_R$ and $L_R$. 
Then we can assume that $M=K_R$ on $X_R\times\{0,1\}$ and $N=L_R$ on $Y_R\times\{0,1\}$ for the following reason. \par
  Extend $F$ to $\tilde F:X_R\times[-1,\,2]\to Y_R\times[-1,\,2]$ trivially, i.e., $\tilde F(\cdot,t)=(g_R,t)$ for $t\in
[-1,\,0]$ and $\tilde F(\cdot,t)=(h_R,t)$ for $t\in[1,\,2]$, and extend trivially $F|_{X\times\{0,1\}}:K\to L$ to a 
cellular map $\tilde F|_{X\times([-1,0]\cup[1,2])}:\tilde K\to\tilde L$, where $\tilde K$ is the union of 
$$
\multline
K,\quad\{\sigma\times\{-1\},\sigma\times\{0\},\sigma\times[-1,\,0]:\sigma\times\{0\}\in K\}\\
\text{and}\ \{\sigma\times\{1\},\sigma\times\{2\},\sigma\times[1,\,2]:\sigma\times\{1\}\in K\}
\endmultline
$$
and $\tilde L$ is defined by $L$ in the same way. 
Extend $N$ to a family of cells $N\cup\tilde L_R$ and then subdivide the family to an $R$-simplicial decomposition $\tilde N$ of $Y_R
\times[-1,\,2]$ so that $\tilde N^0=N^0\cup\tilde L_R^0$ (Proposition 2.9 in [R-S]). 
(Though $N\cup\tilde L_R$ is not necessarily a cell complex, the proof of Proposition 2.9 in [R-S] shows we can subdivide 
$N\cup\tilde L_R$ to an $R$-simplicial complex without introducing new vertices.) 
Then $\tilde N=\tilde L_R$ on $Y_R\times\{-1,2\}$. 
Set 
$$
\multline
\hat M=\\
\{\sigma_1\cap\tilde F^{-1}(\sigma_2):\sigma_1\in M\ \text{or}\ \sigma_1\in\tilde K_R\ \text{with}\ \sigma_1\cap X_R
\times([-1,\,0)\cup(1,\,2])\not=\emptyset,\ \sigma_2\in\tilde N\}. 
\endmultline
$$
In the same way as above we subdivide $\hat M$ to a simplicial complex without new vertices $\tilde M$. 
Then $\tilde M=\tilde K_R$ on $X_R\times\{-1,2\}$ and $\tilde F:\tilde M\to\tilde N$ is simplicial. 
Hence by the linear homeomorphism from $[-1,\,2]$ to $[0,\,1]$ carrying $-1,2$ to $0,1$, respectively, we can assume from the 
beginning that $M=K_R$ on $X_R\times\{0,1\}$ and $N=L_R$ on $Y_R\times\{0,1\}$. \par
  By the strengthened statement we can modify $F:M\to N$ without changing $M|_{X_R\times\{0,1\}}$ and $N|_{Y_R\times\{0,
1\}}$ so that $M$ and $N$ are the $R$-extensions of some $\R$-simplicial complexes. 
Hence $g$ and $h$ are $\R$-homotopic. \par
   The other claims in the lemma are proved in the same way.
\qed
\enddemo

By the above lemmas we have 
\proclaim{Remark 1.3}
Let a property on a finite number of compact polyhedra $(X_i,X_{i,j})_j$ in $R^n$ and a hypothesis on them be stated by 
only a finite combination of terms of PL topology in [R-S] or PL microbundle in [M]. 
Assume that the property and the hypothesis are not concerned with special numbers in $R-\Q$ and preserved under $R$-PL 
homeomorphisms $(X_i,X_{i,j})_j\to(X_i',X'_{i,j})_j$. 
Let $(Y_i,Y_{i,j})_j$ be polyhedra in $\R^n$ such that $(X_i,X_{i,j})_j$ are $R$-PL homeomorphic to $(Y_{iR},Y_{i,jR})_j$. 
Then the property holds under the hypothesis for $(X_i,X_{i,j})_j$ if and only if it does for $(Y_i,Y_{i,j})_j$. 
\endproclaim

The term {\bf cardinal number}, for example, is not such a term. 
We cite some properties which satisfy the conditions. 
(i) $X_1$ and $X_2$ are PL homeomorphic; 
(ii) $X_1$ is a PL manifold; 
(iii) a PL submanifold $X_2$ of a PL manifold $X_2$ is locally flat or admits a normal PL microbundle; 
(iv) $(X_1,X_{1,i})_i$ and $(X_2,X_{2,i})_i$ have the same (simple) homotopy type; 
(v) a PL manifold $X_1$ is obtained by a finite sequence of PL surgeries from a PL manifold $X_2$; 
(vi) $X_1$ is PL contractible; 
(vii) PL manifolds $X_1$ and $X_2$ are PL cobordant. 
By the remark, (co)homology groups, homotopy groups, linking numbers, Whitehead groups, codordims groups, etc., do not depend 
on the choice of $R$, the PL Poincar\' e conjecture is true if and only if it is so for $\R$, and, moreover, I do not know 
well-known theorems on compact polyhedra which hold for $\R$ but not for general $R$. 

\head \S 2. o-minimal Hauptvermutung \endhead
A manifold means that without boundary. 
From now on $R$ denotes a real closed field, and we assume $R\supset\R$ for simplicity of notation. 
We fix an o-minimal structure on $R$ which expands $(R,<,0,1,+,\cdot)$. 
See [v] for the definition and fundamental properties of an o-minimal structure. 
The conjecture that two homeomorphic polyhedra in $\R^n$ are PL homeomorphic is known as Hauptvermutung. 
Milnor showed that Hauptvermutung is not true for polyhedra and Siebenmann did for PL manifolds. 
Existence of such polyhedra is due to {\bf wildness} of some homeomorphisms. 
On the other hand the category of definable sets and definable maps is topologically {\bf tame} as claimed in [S$_3$] and [v]. 
Indeed o-minimal Hauptvermutung is true. 
\proclaim{Theorem 2.1}
If compact polyhedra $(X_1,X_{1,i})_{i=1,...,k}$ and $(X_2,X_{2,i})_{i=1,...,k}$ in $R^n$ are definable homeomorphic then they 
are PL homeomorphic. 
\endproclaim
The theorem was proved in [S-Y] for $(\R,<,0,1,+,\cdot)$ (i.e., in the real semialgebraic structure), in [C$_1$] for $(R,<,0,1,+,
\cdot)$ and then in [S$_3$] for $\R$ in any o-minimal structure. 
[S$_3$] explains an o-minimal structure on $\R$, and some of arguments there work on general $R$. 
The proof of o-minimal Hauptvermutung is, however, false for $R$. 
The proof in [S-Y] and [S$_3$] used essentially an approximation theorem in the proof of uniqueness of $C^\infty$ 
triangulation by Cairnes-Whitehead and the approximation theorem does not hold if $R$ is larger than $\R$. 
[C$_1$] proved by complexities of semialgebraic sets, the Tarski-Seidenberg principle and the result of [S-Y]. \par
  A {\it compact} set in $R^n$ stands for a closed and bounded subset of $R^n$. 
(A closed, bounded and polyhedral subset of $R^n$ is not necessarily a compact polyhedron in the sense in \S 1, i.e., a finite union of simplexes. 
This does not cause confusions. 
Indeed we are interested in only definable sets and a closed, bounded, definable polyhedral subset of $R^n$ is a compact polyhedron.) 
Let $\{E_i\}_{i=1,...,k}$ be a finite stratification (i.e., partition) of a compact definable set $E$ in $R^n$ into definable 
$C^1$ manifolds. 
(The strata of a stratification are always, except once, of class $C^1$.)  
$\{E_i\}_i$ is {\it compatible with} a finite family of definable sets $\{X_j\}_j$ in $R^n$ if each $E\cap X_j$ is the union of 
some $E_i$. 
$\{E_i\}_i$ satisfies the {\it frontier condition} if $E_j\cap(\overline E_i-E_i)\not=\emptyset$ implies $E_j\subset\overline E_i-E_i$ for 
each $i$ and $j$. 
Assume the closure $\overline E_1=E$. 
Let $e_1$ be a point of $E_1$. 
We say that a substratification $\{e_1\}\cup\{E'_i,E_i\}_{i=2,...,k}$ of $\{E_i\}_{i=1,...,k}$ is {\it obtained by starring at} 
$e_1$ {\it from} $\{E_i
\}_{i=1,...,k}$ if there is a definable homeomorphism $\rho$ from the cone $a*(E-E_1)$ with vertex $a$ and base $E-E_1$ to $E$ 
such that $\rho|_{E_i}=\id$ and $\rho(a*E_i)=\{ e_1\}\cup E'_i\cup E_i$ for $i=2,...,k$. 
(We regard $a*\emptyset$ as $\{a\}$.) 
Let $g:X\to Y$ be a definable $C^0$ map between definable sets. 
A {\it definable} $C^1$ {\it stratification} of $g$ is a pair of finite stratifications $\{X_i\}_i$ of $X$ and $\{Y_j\}_j$ 
of $Y$ into definable $C^1$ manifolds such that for each $i$, the restriction $g|_{X_i}$ is a surjective $C^1$ submersion 
to some $Y_j$. 
We write $g:\{X_i\}_i\to\{Y_j\}_j$. 
Let $A_k\subset X$ and $B_l\subset Y$ be a finite number of definable sets. 
Then (II.1.17) in [S$_3$], whose $R$-case is proved in the same way, states that $g$ admits a definable $C^1$ stratification 
$\{X_i\}_i\to\{Y_j\}_j$ such that $\{X_i\}_i$ and $\{Y_j\}_j$ are compatible with $\{A_k\}_k$ and $\{B_l\}_l$, respectively, 
and satisfy the frontier condition. \par
  We know that any compact definable set is definably homeomorphic to some compact polyhedron. 
Moreover, given a finite simplicial complex $K$ in $R^n$ and a finte number of compact definable sets $\{X_i\}_i$ in $|K|$, then there exists a definable isotopy $\tau_t,\ 0\le t\le 1$, of $|K|$ preserving $K$ such that $\{\tau^{-1}(X_i)\}_i$ are polyhedra (see Theorem II.2,1, Remark II.2.4 and Lemma II.2.7 in [S$_3$], whose $R$-case is proved in the same way). 
We call this fact the {\it triangulation theorem of definable sets} and $\tau_1:\{\tau^{-1}_1(X_i)\}_i\to\{X_i\}_i$ a {\it definable triangulation} of $\{X_i\}_i$. 
(From now we omit to mention [S$_3$].) \par
  We prove theorem 2.1 in process of proving the following theorem of triangulations of definable continuous functions. 
\proclaim{Theorem 2.2}
(1) Let $f:X\to R$ be a definable continuous function on a compact definable set in $R^n$. 
Then there exist a compact polyhedron $Y$ in  $R^n$ and a definable homeomorphism $\pi:Y\to X$ such that $f\circ\pi$ is PL. 
Given finitely many compact definable subsets $X_i$ of $X$, then we can choose $Y$ and $\pi$ so that $\pi^{-1}(X_i)$ are polyhedra. 
\par
  (2) Moreover if $X$ is the underlying polyhedron to a finite simplicial complex $P$ then we can choose $Y$ and $\pi$ in (1) so that $Y=X$ 
and there exists a definable isotopy $\pi_t,\ 0\le t\le 1$, of $X$ from $\id$ to $\pi$ preserving $P$. \par
  (3) If $f$ is PL on polyhedral $X$ from the beginning in (3), $X_i\cap f^{-1}(R-(s_1,\,s_2))$ are polyhedra for some $s_1<s_2\in R$ and 
$s'_1$ and $s'_2$ are numbers in $R$ with $s'_1<s_1<s_2<s'_2$ then we can choose $\pi$ and $\pi_t$ so that $f\circ\pi_t=f$ for 
$t\in[0,\,1
]$ and $\pi_t=\id$ on $f^{-1}(R-(s'_1,\,s'_2))$. 
\endproclaim 
\proclaim{Complement}
Moreover $\pi:Y\to X$ is unique in the following sense. 
Let $\pi':Y'\to X$ be another definable homeomorphism with the same properties as $\pi:Y\to X$ in (1). 
Then there exists a definable isotopy $\omega_t:Y'\to Y,\ 0\le t\le 1$, of $\pi^{-1}\circ\pi':Y'\to Y$ through homeomorphisms 
such that $\omega_1$ is PL, $f\circ\pi\circ\omega_t=f\circ\pi'$ for $t\in[0,\,1]$ and $\omega_t(\pi^{\prime-1}(
X_i))=\pi^{-1}(X_i)$ for each $i$ and $t$. 
\endproclaim 
We call $f\circ\pi:Y\to R$, or $\pi:Y\to X$, in theorem 2.2,(1) a {\it definable triangulation} of $f$. 
[S$_2$] showed a semialgebraic triangulation of a semialgebraic $C^0$ function on a compact semialgebraic set in $R^n$. 
The idea of the following proof comes from its proof. 
The theorem is partially proved in [C$_2$]. 
\demo{Proof of theorem 2.2,(2) in a special case}
By the triangulation theorem of definable sets we assume $X$ is the underlying polyhedron to a finite simplicial complex $P$ and $\{X_i\}_i=P$. 
We will triangulate $\graph f$ so that the triangulation is the graph of some triangulation of $f$. 
Assume $n>1$ and proceed by induction on $n$ because the case of $n=0$ is trivial and that of $n=1$ is clear by the following fact. 
A definable $C^0$ function on an interval in $R$ is monotone on each stratum of some finite stratification of the interval. \par
  Let $p:R^n\times R\to R$ denote the projection. 
Set $A=\graph f$ and $A_t=\{x\in R^n:(x,t)\in A\}$ for each $t\in R$. 
Forget the property $\dim A\le n$ for the sake of induction process. 
Let $\{A_j\}_j$ be a finite stratification of $A$ into definable connected $C^1$ manifolds which is compatible with $\{X_i
\times R\}_i$ and satisfies the frontier condition and such that $p|_A:\{A_j\}_j\to\{p(A_j)\}_j$ is a definable $C^1$ 
stratification of $p|_A$ (II.1.17). 
Let $A'$ denote the union of $A_j$ with $\dim A_{j t}<n$ for each $t\in R$. 
Set $S^{n-1}=\{\lambda\in R^n:|\lambda|=1\}$. 
Let $T_t\subset S^{n-1}$ denote the set of singular directions $\lambda$ for $A'_t$, i.e., $\lambda\in S^{n-1}$ such that 
$A'_t\cap(a+R\lambda)$ has interior points in the line $a+R\lambda$ in $R^n$ for some $a\in A'_t$. 
Set $T=\{(\lambda,t)\in S^{n-1}\times R:\lambda\in T_t\}$. 
Then Lemma II.2.2$'$, which is proved for general $R$ in the same way, states that $\overline T$ is definable and $\overline T
\cap(S^{n-1}\times\{t\})$ is of dimension\,$<n-1$ for each $t$. 
Hence $\dim \overline T<n$ and the restriction to $\overline T$ of the projection $S^{n-1}\times R\to S^{n-1}$ is a finite-to-one 
map except onto some definable subset $S$ of $S^{n-1}$ of dimension\,$<n-1$. 
Choose $s\in S^{n-1}-S$. 
Then $\{s\}\times R\cap \overline T$ is finite. 
Therefore, by changing linearly the coordinate system of $R^n$ we can assume that $(1,0,...,0)\in R^n$ is not a singular 
direction for $A'_t$ for each $t\in R$ except a finite number of points, say $B$. 
It follows that the restriction to $A'-R^n\times B$ of the projection $p_1:R^n\times R\to R^{n-1}\times R$ forgetting the 
first factor is a finite-to-one map. 
Note each $A_j$ outside of $A'$ is open in $R^n\times R$ or included and open in $R^n\times\{t\}$ for some $t$. \par
  After then keeping the same notation we substratify $\{A_j\}_j$ so that $p_1|_A:\{A_j\}_j\allowmathbreak\to\{p_1(A_j)\}_j$ is 
a definable $C^1$ stratification of $p_1|_A$ (II.1.17). 
Note first given a substratification $\{D_k\}_k$ of $\{p_1(A_j)\}_j$ into definable $C^1$ manifolds with the frontier condition, then 
$p_1|_A:\{A_j\cap p_1^{-1}(D_k)\}_{j,k}\to\{D_k\}_k$ is a definable $C^1$ stratification of $p_1|_A$ and $\{A_j\cap p_1^{-1}(D_k
)\}_{j,k}$ satisfies the frontier condition. 
Secondary, $p_1|_{A_j-R^n\times B}$ is a diffeomorphism onto $p_1(A_j-R^n\times B)$ for each $A_j$ in $A'$. 
Thirdly, for each $A_j$ outside of $A'$ there are two $A_{j_1}$ and $A_{j_2}$ such that $p_1(A_{j_1})=p_1(A_{j_2})=p_1(A_j)$ 
and $A_{j_1}\cup A_{j_2}\subset A'\cap(\overline A_j-A_j)$. 
By the last property we will see that we can ignore such $A_j$'s. \par 
  Repeat the same arguments for definable sets $p_1(A'-R^n\times B)$ and $\{p_1(A_j-R^n\times B):A_j\subset A'\}_j$ and so on, 
and let $R^{n-1}\times R
\overset p_2\to\longrightarrow R^{n-2}\times R\longrightarrow\cdots\overset p_n\to\longrightarrow R$ be the projections 
forgetting the respective first factors. 
Then by changing linearly the coordinate systems of $R^{n-1},...,R^2$ in sequence we modify $A$ as follows. 
There exist a finite set $C$ in $R$ and finite stratifications with the frontier condition $\{A_{j,k}:k=1,...\}$ of $p_{j-1}\circ
\cdots\circ p_1(A)-R^{n+1-j}\times C$ into definable connected $C^1$ manifolds for each $j=1,...,n+1$ such that $\{A_{1,k}\}_k$ 
is a substratification of $\{A_j-R^n\times C\}_j$, $p_j|_{p_{j-1}\circ\cdots\circ p_1(A)}:\{A_{j,k}\}_k\to\{A_{j+1,k}\}_k$ is 
a definable $C^1$ stratification for each $j\le n$ and, moreover, for each $(j,k_0)$ with $j\le n$ 
$$
\multline
\overline A_{j,k_0}-R^{n+1-j}\times C=
\{(x_j,...,x_n,t)\in R\times(\overline A_{j+1,k_1}-R^{n-j}\times C):\\
\phi_{j,k_0}(x_{j+1},...,x_n,t)\le x_j\le\psi_{j,k_0}(x_
{j+1},...,x_n,t)\}
\endmultline
$$
for some definable $C^0$ functions $\phi_{j,k_0}$ and $\psi_{j,k_0}$ on some $\overline A_{j+1,k_1}-R^{n-j}\times C$ with $\phi_
{j,k_0}=\psi_{j,k_0}$ on $\overline A_{j+1,k_1}-R^{n-j}\times C$ or $\phi_{j,k_0}<\psi_{j,k_0}$ on $A_{j+1,k_1}$. 
It follows from the frontier condition that the following sets are elements of $\{A_{j,k}\}_k$ for each stratum $A_{j+1,k_2}$ in $\overline A_{j+1,k_1}$. 
$$
\gather
\{(x_j,...,x_n,t)\in R\times A_{j+1,k_2}:\phi_{j,k_0}(x_{j+1},...,x_n,t)<x_j<\psi_{j,k_0}(x_{j+1},...,x_n,t)\},\\
\{(x_j,...,x_n,t)\in R\times A_{j+1,k_2}:x_j=\phi_{j,k_0}(x_{j+1},...,x_n,t)\},\\
\{(x_j,...,x_n,t)\in R\times A_{j+1,k_2}:x_j=\psi_{j,k_0}(x_{j+1},...,x_n,t)\}.
\endgather
$$
\par
  {\it Note 1.} For a finite subset $C_0$ of $R$ we can choose the above linear transformations of $R^n,...,R^2$ so that $C_0
\cap C=\emptyset$. \par
  Assume for a while $C=\emptyset$. 
We wish to substratify $\{A_{j,k}\}_k$ in the same way as the barycentric subdivision of a cell complex. 
Choose one point $a_{j,k}$ in each $A_{j,k}$ so that $\{p_j(a_{j,k})\}_k=\{a_{j+1,k}\}_k$ for each $j\le n$. 
Substratify $\{A_{n+1,k}\}_k$ to $\{A'_{n+1,k'}\}_{k'}$ which is defined to be the family of $a_{n+1,k}$ for all $k$ and the connected components of $p_n\circ\cdots\circ p_1(A)-C-\cup\{a_{n+1,k}\}_k$. 
Assume by downward induction on $l$ that $\{A_{j,k}\}_k$ is substratified to $\{A'_{j,k'}\}_{k'}$ for each $j\ge l+1$. 
Then define $\{A'_{l,k'}\}_{k'}$ to be the family of $a_{l,k}$ and the connected components of $p_l^{-1}(A'_{l+1,k'})\cap A
_{l,k}-\{a_{l,k}\}$ for all $k$ and $k'$. 
Note if $p_l^{-1}(A'_{l+1,k'})\cap A_{l,k}$ does not include $a_{l,k}$ then it is connected. 
In the other case it is of dimension\,$\le 1$ and its image under $p_l$ is the point $p_l(a_{l,k})$. 
Note also $\{\overline{A'}_{j,k'}\}_{k'}$ continue to be described by graphs as $\{\overline A_{j,k}\}_k$ are so. 
Let $\phi'_{j,k'}$ and $\psi'_{j,k'}$ denote the functions. \par
  Next keeping the graph property we substratify $\{A'_{j,k'}\}_{k'}$ by triple induction without introducing new strata of 
dimension 0. 
Set $\{A''_{n+1,k''}\}_{k''}=\{A'_{n+1,k'}\}_{k'}$. 
Let $l$ and $c$ be positive and non-negative integers respectively. 
Assume by induction that $\{A'_{j,k'}\}_{k'}$ and $\{A'_{l,k'}:A'_{l,k'}\subset A_{l,k}\}_{k'}$ are substratified to $\{A''_{j,k''}\}_{k''}$ and $\{A''_{l,k''}:A''_{l,k''}\subset A_{l,
k}\}_{k''}$, respectively, for each $j\ge l+1$ and  each 
$A_{l,k}$ of dimension $<c$ so that $\{A''_{j,k''}:A''_
{j,k''}\subset\overline A_{j,k}\}_{k''}$ is obtained by starring at $a_{j,k}$ from $\{A_{j,k}\}\cup\{A''_{j,k''}:A''_{j,k''}
\subset\overline A_{j,k}-A_{j,k}\}_{k''}$ for each $A_{j,k}$ with $j\ge l+1$, or with $j=l$ and $\dim A_{j,k}<c$. 
Consider one $A_{l,k}$ of dimension $c$ and $\{A'_{l,k'}:A'_{l,k'}\subset A_{l,k}\}_{k'}$. 
If $c=0$ then we set clearly $\{A''_{l,k''}:A''_{l,k''}\subset A_{l,k}\}_{k''}=\{a_{l,k}\}$. 
Hence assume $c>0$. 
Then there are two possible cases to consider\,: $A_{l,k}\cap p_l^{-1}(p_l(a_{l,k}))=\{a_{l,k}\}$, or $A_{l,k}\cap p_l^{-1}(p_
l(a_{l,k}))$ is of dimension 1. 
In the former case we set $\{A''_{l,k''}:A''_{l,k''}\subset A_{l,k}\}_{k''}=\{A_{l,k}\cap p_l^{-1}(A''_{l+1,k''})\}_{k''}$. 
Then the starring condition is satisfied, to be precise, $\{A''_{l,k''}:A''_{l,k''}\subset\overline A_{l,k}\}_{k''}$ is 
obtained by starring at $a_{l,k}$ from $\{A_{l,k}\}\cup\{A''_{l,k''}:A''_{l,k''}\subset\overline A_{l,k}-A_{l,k}\}_{k''}$. 
In the latter case we substratify by starring at $a_{l,k}$ as follows. 
By induction hypothesis $\{A''_{l+1,k''}:A''_{l+1,k''}\subset p_l(\overline A_{l,k})\}_{k''}$ is obtained from $\{p_l(A_{l,
k})\}\cup\{A''_{l+1,k''}:A''_{l+1,k''}\subset p_l(\overline A_{l,k})-p_l(A_{l,k})\}_{k''}$ by starring at $p_l(a_{l,k})$. 
Let $A''_{l,k''}$ be such that $A''_{l,k''}\subset\overline A_{l,k}-A_{l,k}$ and $p_l(A''_{l,k''})\subset p_l(\overline A_
{l,k})-p_l(A_{l,k})$. 
Then by hypothesis of starring there exists uniquely $A''_{l+1,k''_0}$ such that $A''_{l+1,k''_0}\subset p_l(A_{l,k})$ 
and $\overline{A''}_{l+1,k''_0}\cap(p_l(\overline A_{l,k})-p_l(A_{l,k}))=p_l(\overline{A''}_{l,k''})$. 
Let $\phi''_{l,k''}$ and $\psi''_{l,k''}$ be the functions which describe $\overline{A''}_{l,k''}$ on $p_l(\overline{A''}
_{l,k''})$. 
Extend $\phi''_{l,k''}$ and $\psi''_{l,k''}$ to definable $C^0$ functions $\phi''_{l,k,k''}$ and $\psi''_{l,k,k''}$ on 
$\overline{A''}_{l+1,k''_0}$ so that they are of class $C^1$ on $A''_{l+1,k''_0}$, 
$$
\gather
\big(\phi''_{l,k,k''}(p_l(a_{l,k})),p_l(a_{l,k})\big)=\big(\psi''_{l,k,k''}(p_l(a_{l,k})),p_l(a_{l,k})\big)=a_{l,k},\\
\phi''_{l,k,k''}=\psi''_{l,k,k''}\ \text{if}\ \phi''_{l,k''}=\psi''_{l,k''},\quad \phi''_{l,k,k''}<\psi''_{l,k,k''}\ 
\text{on}\ A''_{l+1,k''_0}\ \text{otherwise},\\
\phi_{l,k}<\phi''_{l,k,k''}\le\psi''_{l,k,k''}<\psi_{l,k}\ \text{on}\ \overline{A''}_{l+1,k''_0}-p_l(\overline{A''}_{l,k''}). 
\endgather
$$
(Here the extensions are easily constructed if they are not required to be of class $C^1$ on $A''_{l+1,k''_0}$. 
Let $\phi''_{l,k,k''}$ and $\psi''_{l,k,k''}$ be such extensions. 
Then Theorem II.5.2 says that $\phi''_{l.k.k''}|_{A''_{l+1,k''_0}}$ and $\psi''_{l,k,k''}|_{A''_{l+1,k''_0}}$ are approximated 
by definable $C^1$ functions with the same properties so that the approximations are extended to $\overline{A''_{l+1,k''_0}}$ 
and the extensions coincide, respectively, with $\phi''_{l,k,k''}$ and $\psi''_{l,k,k''}$ at $\overline{A''_{l+1,k''_0}}-A''_
{l+1,k''_0}$.) 
Moreover we choose $\{\phi''_{l,k,k''},\psi''_{l,k,k''}:A''_{l,k''}\subset\overline A_{l,k}-A_{l,k},\,p_l(A''_{l,k''})
\subset p_l(\overline A_{l,k})-p_l(A_{l,k})\}_{k''}$ by induction on dimension of $A''_{l,k''}$ so that if $A''_{l,k''_1}\subset\overline{A''}_{l,k''_2}\subset\overline A_{l,k}-A_{l,k}$ with $p_l(A''_{l,k''_2})\subset p_l(\overline A_{l,k})-p_l(A_{l,
k})$ and if $A''_{l+1,k''_{10}}$ and $A''_{l+1,k''_{20}}$ are given from $A''_{l,k''_1}$ and $A''_{l,k''_2}$, respectively, 
as $A''_{l+1,k''_0}$ is given from $A''_{l,k''}$ then 
$$
\phi''_{l,k,k''_1}=\phi''_{l,k,k''_2}\quad\psi''_{l,k,k''_1}=\psi''_{l,k,k''_2}\quad\text{on}\ \overline{A''}_{l+1,k''_{10}}, 
$$
in a word $\{\phi''_{l,k,k''},\psi''_{l,k,k''}\}_{k''}$ is compatible. 
Now let $A_{l,k}$ be divided to the family of the sets $\{a_{l,k}\}$, two connected components of $A_{l,k}\cap p_l^{-1}(p_l(a_{l,k}))$, the graphs of $\phi''
_{l,k,k''}|_{A''_{l+1,k''_0}}$ and $\psi''_{l,k,k''}\allowmathbreak|_{A''_{l+1,k''_0}}$ and 
$$
\gather
\{(x_l,...,x_n,t)\in R\times A''_{l+1,k''_0}:\phi_{l,k,k''}(x_{l+1},...,x_n,t)<x_l<\phi''_{l,k,k''}(x_{l+1},...,x_n,t)\},\\
\{(x_l,...,x_n,t)\in R\times A''_{l+1,k''_0}:\phi''_{l,k,k''}(x_{l+1},...,x_n,t)<x_l<\psi''_{l,k,k''}(x_{l+1},...,x_n,t)\},\\
\{(x_l,...,x_n,t)\in R\times A''_{l+1,k''_0}:\psi''_{l,k,k''}(x_{l+1},...,x_n,t)<x_l<\psi_{l,k}(x_{l+1},...,x_n,t)\}
\endgather
$$
for all possible $k''$ and $k''_0$ given from $k$ and $k''$ as above. 
Let $\{A''_{l,k''}:A''_{l,k''}\subset A_{l,k}\}_{k''}$ denote the division. 
Then clearly the starring condition is satisfied. 
Thus we obtain the required substratification $\{A''_{j,k''}\}_{k''}$. \par
  We will construct simplicial complexes $K_1,...,K_{n+1}$ in $R^n\times R,...,R$, respectively, and definable homeomorphisms 
$\tau_1:|K_1|\to A,\,\tau_2:|K_2|\to p_1(A),...,\tau_{n+1}:|K_{n+1}|\to p_n\circ\cdots\circ p_1(A)$ such that 
$$
\gather
K_j^0=\{A''_{j,k''}\!:\dim A''_{j,k''}=0\}_{k''},\quad\{\tau_j(\Int\sigma)\!:\sigma\in K_j\}=\{A''_{j,k''}\}_{k''}\quad\text{for 
each}\ j,\\
\tau_j=\id\ \text{on}\ K_j^0,\quad p_j\circ\tau_j=\tau_{j+1}\circ p_j\quad\text{for}\ j\le n,\quad\tau_{n+1}=\id,
\endgather
$$
$\tau_1$ is of the form 
$$
\gather
\tau_1(x_1,...,x_n,t)=\big(\tau_{1,1}(x_1,...,x_n,t),\tau_{1,2}(x_2,...,x_n,t),...,\tau_{1,n}(x_n,t),t\big),\\
\tau_j(x_j,...,x_n,t)=\big(\tau_{1,j}(x_j,...,x_,t),...,\tau_{1,n}(x_n,t),t\big),\quad j=2,...,n,
\endgather
$$
and $\tau_{j+1}^{-1}\circ p_j\circ\tau_j=p_j|_{|K_j|}:K_j\to K_{j+1}$ are simplicial maps for $j\le n$. 
Note for each $\sigma\in K_j$ the vertices spanning $\sigma$ are just the {\it vertices} of $\tau_j(\sigma)$, to be precise, 
$\{A''_{j,k''}:A''_{j,k''}\subset\tau_j(\sigma),\ \dim A''_{j,k''}=0\}_{k''}$. 
Hence define $K_j$ for each $j$ to be the family of cells spanned by the vertices of $\overline{A''}_{j,k''}$ for all $A''_{
j,k''}$. 
Let $K^c_j$ be the subfamily of cells for $A''_{j,k''}$ of dimension $\le c$ for each $c$. 
Then we need to see $K_j$ and $K^c_j$ are simplicial complexes. \par
  As above let $l$ and $c$ be positive and non-negative integers. 
Assume by induction we have shown $K_j,\ j\ge l+1$, and $K^{c-1}_l$ are simplicial complexes and constructed $\tau_j,\ j\ge 
l+1$, and a definable homeomorphism $\tau^{c-1}_l:|K^{c-1}_l|\to\cup\{A''_{l,k''}:\dim A''_{l,k''}<c\}_{k''}$ with the required 
conditions. 
Note the vertices of $\overline{A''}_{j,k''}$ then span a simplex of the same dimension for $j\ge l+1$ and for $j=l$ and 
$k''$ with $\dim A''_{j,k''}<c$ since $K_j,\ j\ge l+1$, and $K^{c-1}_l$ are simplicial complexes. 
We need to see $K^c_l$ is a simplicial complex and define $\tau^c_l$. 
If $c=0$, $K^c_l$ is a finite set and, hence, a simplicial complex, and $\tau^c_l$ should be id and satisfies the conditions. 
Assume $c>0$ and choose one $A''_{l,k''}$ of dimension $c$. 
There are two possible cases to consider: $p_l|_{A''_{l,k''}}$ is injective or not. 
If $p_l|_{A''_{l,k''}}$ is injective, then $\overline{A''}_{l,k''}$ is the graph of some definable $C^0$ function on some 
$\overline{A''}_{l+1,k''_0}\,(=p(\overline{A''}_{l,k''}$)), the vertices of $\overline{A''}_{l+1,k''_0}$ span a simplex of 
dimension $c$ by induction hypothesis and, hence, the vertices of $\overline{A''}_{l,k''}$ span a simplex of dimension $c$. 
If $p_l|_{A''_{l,k''}}$ is not injective, then $p_l(a_0)=p_l(a_1)$ for two vertices $a_0$ and $a_1$ of $\overline{A''}_{l,k''}$, 
$p_l(a_1),...,p_l(a_c)$ are distinct one another for the other vertices $a_2,...,a_c$ of $\overline{A''}_{l,k''}$, $p_l(a_1),...,p_l(a_c)$ 
span a simplex of dimension $c-1$ by induction hypothesis and, hence, $a_0,...,a_c$ span a simplex of dimension $c$. 
In the same way we see that given two $A''_{l,k''_1}$ and $A''_{l,k''_2}$ of dimension\,$\le c$ then the intersection of the two simplexes 
spanned by the vertices of $\overline{A''}_{l,k''_1}$ and of $\overline{A''}_{l,k''_2}$ is a common face of them. 
Therefore $K^c_l$ is a simplicial complex. \par
  It remains to define $\tau^c_l$ on each simplex $\sigma$ spanned by the vertices of $\overline{A''}_{l,k''}$ of dimension 
$c$. 
If $p_l|_{A''_{l,k''}}$ is injective, let $\overline{A''}_{l,k''}$ be the graph of a definable $C^0$ function $\phi$ on 
$\overline{A''}_{l+1,k''_0}$. 
Then set 
$$
\tau^c_l(x_l,...,x_n,t)=\big(\phi\circ\tau_{l+1}(x_{l+1},...,x_n,t),\tau_{l+1}(x_{l+1},...,x_n,t)\big)\ \text{for}\ (x_l,...,x_n,t)
\in\sigma.
$$
Clearly $\tau^c_l|_\sigma$ satisfies the conditions. 
Assume $p_l|_{A''_{l,k''}}$ is not injective. 
Let $A''_{l,k''}$ lie between the graphs of two definable $C^0$ functions $\phi\le\psi$ on $\overline{A''}_{l+1,k''_0}$. 
Let $\sigma_\phi,\sigma_\psi\in K^{c-1}_l$ such that $p_l(\sigma_\phi)=p_l(\sigma_\psi)=p_l(\sigma)$, $\tau^{c-1}_l(\sigma_
\phi)=\graph\phi$ and $\tau^{c-1}_l(\sigma_\psi)=\graph\psi$. 
It is natural to define $\tau^c_l$ on $\sigma$ to be an extension of $\tau^{c-1}_l|_{\sigma_\phi\cup\sigma_\psi}$ and 
linear on each segment $\sigma\cap p_l^{-1}(a),\ a\in|K_{l-1}|$, to be precise, 
$$
\gather
\tau^c_l\big(sx_l+(1-s)x'_l,x_{l+1},...,x_n,t\big)=\qquad\qquad\qquad\qquad\qquad\qquad\qquad\qquad\\
\qquad\qquad\qquad\qquad s\tau^{c-1}_l(x_l,x_{l+1},...,x_n,t)+(1-s)\tau^{c-1}_l(x'_l,x_{l+1},...,x_n,t)\\
\text{for}\ (x_l,x_{l+1},...,x_n,t)\in\sigma_\phi,\ (x'_l,x_{l+1},...,x_n,t)\in\sigma_\psi\quad \text{and}\ s\in[0,\,1].
\endgather
$$\par
  Thus we construct $K_1,...,K_{n+1},\tau_1,...,\tau_{n+1}$. 
Let $q:R^n\times R\to R^n$ denote the projection. 
Then $q\circ\tau_1$ is a definable homeomorphism from $|K_1|$ to $X$ since $\tau(|K_1|)$ is the graph of $f$, 
$$
f\circ q\circ\tau_1(x_1,...,x_n,t)=t\quad\text{for}\ (x_1,...,x_n,t)\in|K_1|
$$
by the form of $\tau_1$, and, hence, $f\circ q\circ\tau_1:|K_1|\to R$ is PL. 
Hence $q\circ\tau_1:|K_1|\to X$ is a definable triangulation of $f$ in the case of $C=\emptyset$. \par
  {\it Note 2.} Later we will perturb $p_1,...,p_{n-1}$ so that $C$ vanishes. 
Then the above arguments work. 
To be precise, let $R^n\times R\overset\tilde p_1\to\longrightarrow R^{n-1}\times R\longrightarrow\cdots\overset\tilde p_{n-1}
\to\longrightarrow R^2$ be definable $C^1$ maps of the form $\tilde p_j(x_j,...,x_n,t)=(\tilde p'_j(x_j,...,x_n,t),t)$ which 
are close to the projections $p_1,...,p_n$ in the compact-open $C^1$ topology and keep all the properties except being linear 
projections. 
(We adopt the compact-open $C^1$ topology for approximations in this proof.) 
Set 
$$
\hat p_j(x_1,...,x_n,t)\!=\!(x_1,...,x_j,\tilde p_j(x_j,...,x_n,t))\ \text{for}\ (x_1,...,x_n,t)\!\in\! R^n\times R,\ 
j\!=\!1,...,n-1,
$$
and $\eta=\hat p_{n-1}\circ\cdots\hat p_1$, whose restriction to $A$ is a definable $C^1$ imbedding into $R^n\times R$. 
Here a definable $C^1$ imbedding is, by definition, a definable $C^0$ map extensible to a definable $C^1$ imbedding of an open definable neighborhood of $A$ in $R^n\times R$ into $R^n\times R$. 
If $A$ is a compact polyhedron then a definable $C^1$ imbedding of $A$ into $R^n\times R$ means a definable $C^1$ imbedding of 
some simplicial decomposition of $A$. 
See p\. 72 and 73 of [S$_3$] for the definition of a definable $C^1$ imbedding of a simplicial complex and its fundamental 
properties. 
Especially Lemma I.3.14 says that a definable $C^1$ approximation of a definable $C^1$ imbedding of a finite simplicial complex 
into a Euclidean space is an imbedding, whose $R$-case is proved in the same way. 
Now come back to general $A$ and translate $A$ by $\eta$, proceed the same arguments for $\eta(A)$, assume $C=\emptyset$ for $\eta(A)$, and let $\tau_1:|K_1|\to
\eta(A)$ be the triangulation. 
Then $q\circ\eta^{-1}\circ\tau_1:|K_1|\to X$ is a definable triangulation of $f$. \par
  Continue to assume $C=\emptyset$, and set $\pi(x)=q\circ\tau_1^{-1}(x,f(x))$ for $x\in X$. 
Then $\pi$ is a map onto $X$ and injective if and only if $q|_{|K_1|}$ is a map onto $X$ and injective respectively. 
We wish to see $q|_{|K_1|}$ is a homeomorphism onto $X$. 
Indeed, if so, $\pi^{-1}:X\to X$ is the required definable triangulation of $f$. 
For a while we assume $q|_A$ is injective. \par
  {\it Surjectivity of $q|_{|K_1|}$.} First $\{q(A''_{1,k''})\}_{k''}$ is a stratification of $X$ into definable $C^1$ manifolds 
compatible with $P$ and with the frontier condition. 
Next for each $\sigma\in K_1$, $q(\sigma)$ is the cell spanned by the images under $q$ of the vertices of $\overline{A''}_
{1,k''}\,(=\tau_1(\sigma))$. 
Hence $\pi(\sigma_0)\subset\sigma_0$ for each $\sigma_0\in P$. 
Therefore $\Im\pi\subset X$. 
We will prove $\Im\pi=X$. 
By induction on dimension, let $\sigma_0\in P$ and assume $\pi(\sigma_1)=\sigma_1$ for each proper face $\sigma_1$ of 
$\sigma_0$ and $\pi|_{\partial\sigma_0}:\partial\sigma_0\to\partial\sigma_0$ is definably homotopic to id preserving the 
faces. 
Then it suffices to see $\pi(\sigma_0)=\sigma_0$ and $\pi|_{\sigma_0}:\sigma_0\to\sigma_0$ is definably homotopic to id so 
that the homotopy is an extension of the given homotopy of $\pi|_{\partial\sigma_0}$. \par
  We prove $\pi(\sigma_0)=\sigma_0$ by reduction to absurdity. 
Assume $\pi(\sigma_0)\not=\sigma_0,\ 0\in\Int\sigma_0$ and $0\not\in\pi(\sigma_0)$. 
Let $\partial L$ denote the simplicial complex consisting of proper faces of $\sigma_0$, and set 
$$
L=\{0,\,0*\sigma,\,\sigma:\sigma\in\partial L\},\quad L/2=\{\sigma/2:\sigma\in L\}.
$$
Let $\pi_0:\sigma_0\to\partial\sigma_0$ be a definable $C^0$ map such that $\pi_0|_{\partial\sigma_0}:\partial\sigma_0\to
\partial\sigma_0$ is definably homotopic to $\id$ preserving $\partial L$, e.g., the composite of $\pi|_{\sigma_0}$ with the 
retraction $\pi(\sigma_0)\ni tx\to x\in\partial\sigma_0$ for $(x,t)\in\partial\sigma_0\times(0,\,1]$ with $tx\in\pi(\sigma_0)$. 
By this homotopy property we can modify $\pi_0$ so that 
$$
\gather
\pi_0(tx)=x\quad\text{for}\ (x,t)\in\partial\sigma_0\times[1/2,\,1]. \\
\quad\pi_0=\id\quad\text{on}\ \partial\sigma_0,\quad 0*\sigma-\Int\sigma_0/2=\pi^{-1}_0(\sigma)-\Int\sigma_0/2\quad\text{for}\ 
\sigma\in\partial L. \tag"Then"
\endgather
$$
Moreover we can assume each $\pi^{-1}_0(\sigma)$ for $\sigma\in\partial L$ is a polyhedron for the following reason. \par
  Apply the triangulation theorem of definable sets to $\{\sigma_0/2,\sigma_0/2\cap\pi^{-1}(\sigma_1),\sigma_2:\sigma_1\in\partial L,\,\sigma_2
\in L/2\}$. 
Then there exists a definable homeomorphism $\xi_1$ of $\sigma_0/2$ preserving $L/2$ such that $\xi_1^{-1}(\sigma_0/2\cap\pi^{-1}
_0(\sigma_1))$ are polyhedra for $\sigma_1\in\partial L$. 
Note $\xi_1(\sigma)=\sigma$ for $\sigma\in L/2$. 
Set 
$$
\xi_t(x)=(1-t)x+t\xi_1(x)\quad\text{for}\ (x,t)\in\partial\sigma_0\times[1/2,\,1]. 
$$
Then $\xi_t:\sigma_0/2\to\sigma_0/2,\ 0\le t\le 1$, is a definable homotopy of id to $\xi_1$ preserving $L/2$. 
Extend $\xi_1$ to a definable map $\xi$ between $\sigma_0$ by 
$$
\xi(tx)=2t\xi_{2-2t}(x/2)\quad\text{for}\ (x,t)\in\partial\sigma_0\times[1/2,\,1]. 
$$
Then $\xi=\id$ on $\partial\sigma_0$, $\xi$ is preserving $L$ and $\xi^{-1}(\pi_0^{-1}(\sigma_1))$ is a polyhedron for each $\sigma_1\in\partial L$. 
Consider $\pi_0\circ\xi:\sigma_0\to\partial\sigma_0$ in place of $\pi_0$. 
Then we can assume from the beginning $\pi_0^{-1}(\sigma_1)$ is a polyhedron for each $\sigma_1\in\partial L$. \par
  Let $L'$ be a simplicial subdivision of $L$ such that each $\pi^{-1}_0(\sigma_1)$ is the union of some simplexes in $L'$. 
Then for each $\sigma\in L'$ there exists uniquely $\delta_\sigma\in\partial L$ such that $\pi_0(\sigma)\subset\delta_\sigma$ 
and $\pi_0(\Int\sigma)\subset\Int\delta_\sigma$. 
Define a map $\pi'_0:L^{\prime 0}\to\partial\sigma_0$ by $\pi'_0(v)=v$ for $v\in L^{\prime 0}\cap\partial\sigma_0$ and so that 
$\pi'_0(v)$ is the barycenter of $\delta_\sigma$ for $v\in L^{\prime 0}-\partial\sigma_0$. 
Then $\pi'_0$ carries the vertices of each $\sigma\in L'$ into $\delta_\sigma$, and, hence, $\pi'_0$ is extended to a PL map $\pi'
_0:\sigma_0\to\partial\sigma_0$ so that $\pi'_0$ is linear on each simplex in $L'$. 
Note $\pi'_0=\id$ on $\partial\sigma_0$. 
Regard $\sigma$ as an $\R$-simplex and apply lemma 1.2 to $\pi'_0:\sigma_R\to\partial\sigma_R$. 
Then $\pi'_0$ is PL homotopic to the $R$-extension of some $\R$-PL map $\pi''_0:\sigma\to\partial\sigma$. 
Here the $R$-extension of $\pi''_0|_{\partial\sigma}:\partial\sigma\to\partial\sigma$ is PL homotopic to the identity map of 
$\partial\sigma_R$. 
Once more, apply lemma 1.2 to these two maps. 
Then $\pi''_0|_{\partial\sigma}$ is PL homotopic to id. 
That contradicts the fact $H_{\dim\sigma_0}(\sigma_0,\partial\sigma_0)=\Z$ in the real case. 
Hence $\pi(\sigma_0)=\sigma_0$. \par
  We extend the given homotopy of $\pi|_{\partial\sigma_0}$ by the method of cone extension (II.1.16) called the Alexander trick. 
Let $\Int\sigma_0$ contain 0. 
Let $g_\partial:\partial\sigma_0\times[-1,\,1]\to\partial\sigma_0$ be a definable $C^0$ map such that $g_\partial(\cdot,-1)=\pi|
_{\partial\sigma_0}$ and $g_\partial(\cdot,1)=\id$. 
Extend $g_\partial(\cdot,-1)$ to $\pi|_{\sigma_0}$ on $\sigma_0$ and $g_\partial(\cdot,1)$ to id on $\sigma_0$, and set 
$$
g(s(x,t))=sg_\partial(x,t)\quad\text{for}\ (x,t)\in\partial(\sigma_0\times[-1,\,1]),\ s\in[0,\,1].
$$
Then $g$ is a well-defined definable $C^0$ map from $\sigma_0\times[-1,\,1]$ to $\sigma_0$ and fulfills the requirements. 
Thus we have shown surjectivity and also that $\pi$ is preserving $P$. \par
  {\it Injectivity of $q|_{|K_1|}$.} 
We will modify the above linear changing of the coordinate systems of $R^n,...,R^2$ so that $q|_{|K_1|}$ is injective. 
For each $j=1,...,n+1$, let $E_j$ denote the union of $A_{j,k}$ with $\dim A_{j,k t}<n+1-
j$ for each $t\in R$, where $A_{j,kt}=\{x\in R^{n+1-j}:(x,t)\in A_{j,k}\}$ as above. 
Let $q_j:R^{n+1-j}\times R\to R^{n+1-j}$ denote the projection. 
Remember $E_1$ and $q_1$ equal $A'$ and $q$, respectively, and $q_1|_{E_1}$ is injective. 
We, however, do not know whether the other $q_j|_{E_j}$ are finite-to-one maps. 
Assume they are so and, moreover, using (II.1.17) choose $\{A_{j,k}\}_k$ and $\{A''_{j,k''}\}_{k''}$ so that $q_j|_{E_j}:\{A_{j,k}
:A_{j,k}\subset E_j\}_k\to\{q_j(A_{j,k}):A_{j,k}\subset E_j\}_k$ and $q_j|_{E_j}:\{A''_{j,k''}:A''_{j,k''}\subset E_j\}_{k''}\to
\{q_j(A''_{j,k''}):A''_{j,k''}\subset E_j\}_{k''}$ are definable $C^1$ stratifications of $q_j|_{E_j}$. 
Then by induction on $n-j$ and by the method of construction of $\tau_1$ we see 
$$
q_j\circ\tau_j^{-1}(\overline{A''}_{j,k''_1})=q_j\circ\tau_j^{-1}(\overline{A''}_{j,k''_2})\quad\text{if}\quad q_j(A''_{j,k''_1})=q_j(A''_{j,k''_2}), 
$$
the vertices of each $q_j(\overline{A''}_{j,k''})$ span a simplex of the same dimension as $A''_{j,k''}$, say $\sigma_{j,k''
}$, which equals $q_j\circ\tau_j^{-1}(\overline{A''}_{j,k''})$, $\sigma_{j,k''_1}$ is a proper face of $\sigma_{j,k''_2}$ if 
$A''_{j,k''_1}\subset \overline{A''}_{j,k''_2}-A''_{j,k''_2}$, and then $\{\sigma_{j,k''}\}_{k''}$ is a simplicial complex. 
In the case of $j=1$ there are no $k''_1\not=k''_2$ such that $q_1(A''_{1,k''_1})=q_1(A''_{1,k''_2})$ since $q_1|_A$ is 
injective. 
Hence $q_1|_{|K_1|}:K_1\to\{\sigma_{1,k''}\}_{k''}$ is a simplicial isomorphism. 
Therefore $q|_{|K_1|}$ is injective. \par
  We changed linearly the coordinate systems of $R^n,...,R^2$ in sequence. 
We need to choose the linear transformations so that the assumption---$q_j|_{E_j}$ are finite-to-one maps---is satisfied. 
It suffices to consider only $A\ (\subset R^n\times R)$ assuming $q|_A$ is a finite-to-one map in place of injectivity and to 
find a direction $\lambda\in S^{n-1}$ not singular for $A_t$ for each $t\in R$ except a finite number of points and such that 
$q_2|_{E_2}$ is a finite-to-one map. 
Go back to the stratification $\{A_j\}_j$ of $A$ such that $p|_A:\{A_j\}_j\to\{p(A_j)\}_j$ is a definable $C^1$ 
stratification of $p|_A$. 
From the beginning we assume $q|_{A_j}$ is a $C^1$ imbedding for each $j$. 
Since we can ignore $A_j$ with $\dim A_{jt}=n$ for some $t\in R$ we assume $\dim A_{jt}<n$ for any $j$ and $t$. 
Note then $p|_{A_j}$ is regular for each $A_j$ of dimension $n$. 
By definition $p_1:R^n\times R\to R^{n-1}\times R$ is the projection to the direction $\lambda\in S^{n-1}$. 
Then we have the {\bf canonical} method of substratification of $\{A_j\}_j$ to $\{A_{1,j}\}_j$ so that $p_1|_A:\{A_{1,j}\}_j\to\{
A_{2,j}\}_j$ is a definable $C^1$ stratification of $p_1|_A$ (see II.1.17). 
It follows from the definition of the canonical substratification that  
$$
\gather
\cup\{A_{2,j}:\dim A_{2,j}<n\}_j=\qquad\qquad\qquad\qquad\qquad\qquad\qquad\\
p_1\big(\cup\{A_j:\dim A_j<n\}_j\cup\cup\{\Sing p_1|_{A_j}:\dim A_j=n\}_j\big),
\endgather
$$
where $\Sing$ denotes the singular point set of a $C^1$ map. 
Clearly $E_2$ is included in $\cup\{A_{2,j}:\dim A_{2,j}<n\}_j$ and we can ignore $A_{2,j}$ outside of $E_2$. 
Hence it suffices to choose $\lambda$ not singular for $A_t$ for each $t\in R$ except a finite number of points such that 
the restriction of $q_2$ to the right side of the above equality is a finite-to-one map. 
Therefore what we prove is that for each $A_j$, the set of $\lambda\in S^{n-1}$ such that the restriction of $q_2$ to 
$p_1(A_j)$ if $\dim A_j<n$, or to $p_1(\Sing p_1|_{A_j})$ otherwise, is not a finite-to-one map is included in a definable 
set of dimension $<n-1$. 
(By the above method of construction of $\{A_{j,k}\}_k,\ j=1,...,n+1$, we substratify the canonical stratification $\{A_{2,j}
\}_k$ and ignore $A_{2,j}$ of dimension $n$ in the canonical one. 
Hence $q_2|_{E_2}$ is a finite-to-one map for any substratification if $q_2|_{E_2}$ is a finite-to-one map for the canonical 
one.) \par
  Assume $\dim A_j<n$ and consider $q_2|_{p_1(A_j)}$. 
Let $p'_1:R^n\to R^{n-1}$ denote the projection to the direction $\lambda$. 
Then $\dim q_1(A_j)<n$ and, hence, $p'_1\circ q_1\ (=q_2\circ p_1)$ on $A_j$ is a finite-to-one map for $\lambda\in S^{n-1}$ 
except a definable subset of $S^{n-1}$ of dimension $<n-1$. 
We choose $\lambda$ outside of this definable subset. 
If $q_2\circ p_1|_{A_j}$ is a finite-to-one map, so is $q_2|_{p_1(A_j)}$. 
Thus the problem is solved. \par
  {\it Note 3.} We will approximate $p_1$ by a definable $C^1$ map $\tilde p_1:R^n\times R\to R^{n-1}\times R$ and define 
$\hat p_1$ as in note 2 so that $p'_1$ on $q_1\circ\hat p_1(A_j)$ is a finite-to-one map. 
Then $p|_{\hat p_1(A)}:\{\hat p_1(A_k)\}_k\to\{p(A_k)\}_k$ is a definable $C^1$ stratification of $p|_{\hat p_1(A)}$, and 
we can replace $A$ with $\hat p_1(A)$ keeping the property that $q_2|_{p_1(A_j)}$ is a finite-to-one map. \par
  Assume $\dim A_j=n$. 
Let $T_{(x,t)}A_j$ denote the tangent space of $A_j$ at a point $(x,t)\in A_j$. 
By definition $p_1|_{A_j}$ is singular at $(x,t)$ if and only if $(\lambda,0)\in T_{(x,t)}A_j$. 
Hence $q_2|_{p_1(\Sing p_1|_{A_j})}$ is not a finite-to-one map if and only if there exist an open interval $I$ in $R$ and a 
definable $C^1$ curve $\xi:I\to A_j$ of the form $\xi(t)=(\xi'(t),t)$ such that $(\lambda,0)\in T_{\xi(t)}A_j$ for $t\in I$ 
and $p'_1\circ\xi'$ is constant. 
By the last property $\xi'(t)$ is of the form $\rho(t)\lambda+a$ for some definable $C^1$ function $\rho$ on $I$ and $a\in R^n$. 
Hence $\Im d\xi_t$ is a line $R(\frac{d\rho}{dt}(t)\lambda,1)$ for each $t\in I$. 
It follows $(\frac{d\rho}{dt}(t)\lambda,1)\in T_{\xi(t)}A_j$. 
Therefore $(0,1)\in T_{\xi(t)}A_j$, which contradicts the property that $q|_{A_j}$ is a $C^1$ finite-to-one map. 
Thus $q_2|_{p_1(\Sing p_1|_{A_j})}$ is always a finite-to-one map, which completes the proof in the case of $C=\emptyset$. \par
  In place of the assumption $C=\emptyset$ we forget $\{X_i\}$ and assume $X$ is a PL manifold with boundary of the same dimension $n$ as the dimension of the 
ambient space of $X$, $f\ge 0$, $f^{-1}(0)=\partial X$ and $f$ is PL on a neighborhood of $\partial X$ in $X$. 
Let $\epsilon>0$ be such that $A\cap R^n\times[0,\,3\epsilon]$ is a polyhedron. 
Note that $X$ is the closure of the union of some connected components of $R^n-\partial X$ and $X$ is uniquely determined by 
$\partial X$, to be precise, if $X'$ is another compact PL manifold with the same boundary as $X$ then $X'=X$, which is clear 
for $R=\R$ and for general $R$ by remark 1.3. 
We wish to delete $C$. 
Go back to the definable $C^1$ stratification $p_1|_A:\{A_j\}_j\to\{p_1(A_j)\}_j$ such that $p_1|_{A'-R^n\times B}$ is a 
finite-to-one map and $p'$ on $q_1(A_j)$ is a finite-to-one map for $A_j$ of dimension $<n$. 
Approximate $p_1$ by $\tilde p_1$, in the compact-open $C^1$ topology as mentioned, and define $\hat p_1$ as in note 2. 
Assume that (i) $p_1|_{\hat p_1(A')}$ is a finite-to-one map, i.e., $B$ vanishes for $\hat p_1(A')$, (ii) $\tilde p_1=p_1$ on 
$A\cap R^n\times[0,\,\epsilon]$, (iii) the value of $\tilde p_1-p_1$ depends on only $x_1$ on $A-R^n\times[0,\,2\epsilon]$ and 
(iv) $p'_1$ on $q\circ\hat p_1(A_j-R^n\times[0,\,2\epsilon])$ is a finite-to-one map for $A_j$ of dimension $<n$. 
Then $q|_{\hat p_1(A)}$ is injective for the following reason. 
If we fix $x_1$, $\hat p_1$ is a parallel translation of $A\cap\{x_1\}\times R^{n-1}\times[2\epsilon,\,\infty)$ by (iii). 
Hence $q|_{\hat p_1(A)\cap R^n\times[2\epsilon,\,\infty)}$ is injective. 
On the other hand, $A\cap R^n\times[0,\,3\epsilon]$ is a compact polyhedron and $\hat p_1$ is close to the identity map in the 
compact-open $C^1$ topology. 
Hence $q\circ\hat p_1|_{A\cap R^n\times[0,\,3\epsilon]}$ is a $C^1$ imbedding since $q|_{A\cap R^n\times[0,\,3\epsilon]}$ is a 
$C^1$ imbedding (Lemma I.3.14). 
Therefore $q|_{\hat p_1(A)}$ is injective. 
By the same reason it follows from (iv) that $p'_1$ on $q\circ\hat p_1(A_j)$ is a finite-to-one map for $A_j$ of dimension $<n$. 
By (ii) we see $q(\hat p_1(A\cap R^n\times\{0\}))=f^{-1}(0)=\partial X$. 
Hence $q(\hat p_1(A))$ is a compact PL manifold in $R^n$ with the same boundary as $X$ and equals $X$ as noted above. 
Consequently $q|_{\hat p_1(A)}$ is a homeomorphism onto $X$. \par
  Substratify $\{\hat p_1(A_j)\}_j$ so that $p_1|_{\hat p_1(A)}:\{\hat p_1(A_j)\}_j\to\{p_1\circ\hat p_1(A_j)\}_j$ is a definable 
$C^1$ stratification and keep the same notation. 
In the same way we define $p_2,\tilde p_2,\hat p_2$ for $\cup\{p_1\circ\hat p_1(A_j):\dim p_1\circ\hat p_1(A_j)_t<n-1,t\in R\}_j$ 
with the properties corresponding to (i), (ii), (iii) and (iv). 
Repeat the same arguments and obtain $p_3,\tilde p_3,\hat p_3,...,p_{n-1},\tilde p_{n-1},\hat p_{n-1}$ and $\eta=\hat p_{n-1}
\circ\cdots\hat p_1$. 
Then $q|_{\eta(A)}$ is a homeomorphism onto $X$ by the same reason as for $q|_{\hat p_1(A)}$, and the arguments of triangulation 
of $\eta(A)$ via $R^n\times R\overset{p_1}\to\longrightarrow\cdots\overset{p_n}\to\longrightarrow R$ proceed with empty $C$ so 
that the condition in the above proof of injectivity of $q|_{|K_1|}$ is satisfied. 
Hence there exist a simplicial complex $K$ in $R^n\times R$ and a triangulation $\tau:|K|\to\eta(A)$ such that $q|_{|K|}$ is a homeomorphism onto $X$ and $\tau$ is of the form $\tau(x,t)=(\tau_1(x,t),t)$. 
Then $q\circ\eta^{-1}\circ\tau\circ(q|_{|K_1|})^{-1}:X\to X$ is the required definable triangulation of $f$. 
Thus what we need to prove is the following statement. \par
  {\it Statement.} Given $A,A',p_1$ and $\{A_j\}_j$ as above, then $p_1$ is approximated by $\tilde p_1$ so that (i), (ii), (iii) 
and (iv) are satisfied. \par
  {\it Proof of the statement.} We can forget (ii) and replace (iii) and (iv) with the conditions that (iii$)'$ the value of $\tilde p_1
-p_1$ depends on only $x_1$ and (iv$)'$ $p'_1$ on $q\circ\hat p_1(A_j)$ is a finite-to-one map for $A_j$ of dimension $<n$, 
respectively, for the following reason. 
Let $\zeta$ be a definable $C^1$ function on $R$ such that $\zeta=0$ on $(-\infty,\,\epsilon]$ and $\zeta=1$ on $[2\epsilon,\,
\infty)$. 
Assume we have $\tilde p_1$ sufficiently close to $p_1$ with (i), (iii$)'$ and (iv$)'$. 
Set 
$$
\Tilde{\Tilde p}_1(x,t)=(1-\zeta(t))p_1(x)+\zeta(t)\tilde p_1(x)\quad\text{for}\ (x,t)\in R^n\times R,
$$
and define $\Hat{\Hat p}_1$ by $\Tilde{\Tilde p}_1$ as in note 2. 
Then $\Tilde{\Tilde p}_1$ is close to $p_1$, $p_1$ is a finite-to-one map on $\Hat{\Hat p}_1(A')-R^n\times[0,\,2\epsilon]$ by (i) 
and on $\Hat{\Hat p}_1(A')\cap R^n\times[0,\,2\epsilon]$ by Lemma I.3.14 because we can choose $\{A_j\}_j$ so that $A'\cap R^n
\times[0,\,2\epsilon]$ is a polyhedron, hence, (i) holds for $\Hat{\Hat p}_1$, and (ii), (iii) and (iv) are clear for $\Hat{\Hat p}_1$ 
by the definition of $\Tilde{\Tilde p}_1$. 
Thus (i),..,(iv) are replaced with (i), (iii$)'$ and (iv$)'$. \par
  Moreover if (iii$)'$ is satisfied then (iv$)'$ is equivalent to that $p_1|_{\hat p_1(q^{-1}(q(A_j)))}$ is a finite-to-one map 
for $A_j$ of dimension $<n$. 
Hence, adding such $q^{-1}(q(A_j))$ to $A'$ and assuming only that $A'$ is a compact definable set in $R^n\times R$ with $\dim A'
_t<n$ for each $t\in R$ we will find $\tilde p_1$ so that (i) and (iii$)'$ are satisfied. \par
  For each $s=(s_1,...,s_{n+1})\in(\{0\}\times R^{n-1})^{n+1}\subset(R^n)^{n+1}$, define a definable $C^\infty$ submersion $r_s:
R^n\times R\to R^{n-1}\times R$ by 
$$
r_s(x,t)=p_1(x,t)-(\sum_{i=1}^{n+1}x_1^is_i,0)\quad\text{for}\ (x,t)=(x_1,...,x_n,t)\in R^n\times R.
$$
Then $r_0=p_1$, if $s$ is close to 0 in $(R^n)^{n+1}$ then $r_s$ is close to $p_1$ as a $C^1$ map, and if we choose $r_s$ as $\tilde p_1$ then (iii$)'$ is satisfied. 
Hence we will find $s$ close to 0 so that $\tilde p_1=r_s$ satisfies (i). 
For each $t\in R$, let $Y_t$ denote the points $s\in(\{0\}\times R^{n-1})^{n+1}$ such that $A'\cap r_s^{-1}(x_2,...,x_n,t)\ (=
A'_t\times\{t\}\cap r_s^{-1}(x_2,...,x_n,t))$ is of dimension 1 for some $(x_2,...,x_n,t)$. 
Then $Y_t$ and the union $Y$ of all $Y_t,\ t\in R$, are definable subsets of $(\{0\}\times R^{n-1})^{n+1}$ and $r_s$ satisfies 
(i) if and only if $s\not\in Y$. 
Hence it suffices to see $Y_t$ is of codimension $>1$ in $(\{0\}\times R^{n-1})^{n+1}$ for each $t$ because if so then $Y$ is of 
positive codimension and there exists a point outside of $Y$ and arbitrarily close to 0. 
Therefore we consider only $A'_0\subset R^n$ and $Y_0$. 
Moreover we restrict the problem at each point of $A'_0$ by the same reason. 
To be precise, let $x\in A'_0$ and let $Y_{x,0}$ denote the points $s$ of $Y_0$ such that the germ of $A'_0\cap r_s^{-1}(x'_2,...,
x'_n,0)$ at $x$ is of dimension 1 for some $(x'_2,...,x'_n)$. 
Then we only need to prove $\codim Y_{x,0}>n$. \par
  For simplicity of notation we restate the problem. 
Let $A$ be a definable set in $R^n$ of dimension $<n$ containing 0, for each $s=(s_1,...,s_{n+1})\in (\{0\}\times R^{n-1})^
{n+1}$, define $r_s:R^n\to R^{n-1}$ by 
$$
r_s(x)=(x_2,...,x_n)-\sum_{i=1}^{n+1}x_1^is_i\quad\text{for}\ x=(x_1,...,x_n)\in R^n,
$$
and let $Y$ denote the points $s\in(\{0\}\times R^{n-1})^{n+1}$ such that the germ of $A\cap r^{-1}_s(0)$ at 0 is of dimension 1. 
Then what we prove is that $Y$ is of codimension $>n$ in $(\{0\}\times R^{n-1})^{n+1}$. 
We write $Y$ as $Y(A)$ when we need to clarify that $Y$ is defined from $A$. \par
  We can assume $A$ is $C^{n+1}$ smooth because if $A$ is the union of definable sets $Z_1$ and $Z_2$ then $Y(A)=Y(Z_1)\cup Y(Z_2)$ 
and $A$ admits a finite stratification into definable $C^{n+1}$ manifolds (II.1.8, whose $R$-case is proved in the same way). 
For each integer $m$ with $1\le m\le n+1$, let $\theta_m:(\{0\}\times R^{n-1})^{n+1}\to\{0\}\times R^{n-1}$ denote the projection 
to the $m$-th factor. 
Then it suffices to see that $\theta_m\big(Y\cap(\theta_1\times\cdots\times\theta_{m-1})^{-1}(s_1,...,s_{m-1})\big)$ is of positive codimension in $\{0\}
\times R^{n-1}$ for each $(s_1,...,s_{m-1})\in\theta_1\times\cdots\times\theta_{m-1}(Y)\ (\subset(\{0\}\times R^{n-1}
)^{m-1})$. \par
  Case of $m=1$. 
Let $s_1=(0,a_2,...,a_n)\in\theta_1(Y)\subset\{0\}\times R^{n-1}$. 
Then the germ of $r_s^{-1}(0)$ at 0 is included in $A$ for some $s=(s_1,...,s_{n+1})\in(\{0\}\times R^{n-1})^{n+1}$, 
$$
r_s^{-1}(0)=\{(x_1,...,x_n)\in R^n:(x_2,...,x_n)-\sum_{i=1}^{n+1}x_1^is_i=0\},
$$
$r_s^{-1}(0)$ is a smooth curve and, hence, the tangent line of the curve at 0---$\{(x_1,...,x_n)\allowmathbreak\in R^n:(x_2,...,x_n)-x_1s_1=0\}\ (=
\{(x_1,...,x_n)\in R^n:x_i=a_ix_1,\ i=2,...,n\})$---is tangent to the tangent space $T_0A$ of $A$ at 0, namely, $(1,a_2,...,a_n)$ 
is a tangent vector of $A$ at 0. 
Therefore $\theta_1(Y)$ is of positive codimension. \par
  Case of $m=2$. 
Let $s_1\in\theta_1(Y)$. 
Translate linearly $R^{n-1}$ so that $s_1=0$. 
Then the $x_1$-axis is tangent to $T_0A$. 
Let $h:R^n\ni(x_1,...,x_n)\to(x_1,x_1x_2,...,x_1x_n)\in R^n$ be the blowing-up with center 0. 
Then $h^{-1}(A)$ is a definable $C^n$ manifold, and $Y\cap\theta^{-1}_1(0)$ consists of points $(0,s_2,...,s_{n+1})\in\{0\}\times
(\{0\}\times R^{n-1})^n$ such that the germ of $h^{-1}(A)\cap r_{(s_2,...,s_{n+1})}^{-1}(0)$ at 0 is of dimension 1, where 
$r_{(s_2,...,s_{n+1})}:R^n\to R^{n-1}$ is defined by 
$$
r_{(s_2,...,s_{n+1})}(x)=(x_2,...,x_n)-\sum_{i=1}^nx_1^is_{i+1}\quad\text{for}\ x=(x_1,...,x_n)\in R^n. 
$$
Hence we see $\codim\theta_2(Y\cap\theta^{-1}_1(0))>0$ in the same way as in the case of $m=1$. \par
  Repeat the same arguments. 
Then we see $\theta_m\big(Y\cap(\theta_1\times\cdots\times\theta_{m-1})^{-1}(s_1,...,\allowmathbreak s_{m-1})\big)$ is of positive 
codimension in $\{0\}\times R^{n-1}$. 
Thus we prove the statement and, hence, theorem 2.2,(2) in the case where $\{X_i\}=\emptyset$, $X$ is a PL manifold with boundary of dimension $n$, $f\ge 0$, 
$f^{-1}(0)=\partial X$ and $f$ is PL on a neighborhood of $\partial X$ in $X$. \qed
\enddemo
\proclaim{Remark 2.3} 
By note 1 in the above proof we have already shown locally a definable triangulation of $f$. 
To be precise, let $f$ be a definable $C^0$ function on a compact polyhedron $X$, $X_i$ a finite number of definable subsets of 
$X$ and $C$ a finite subset of $R$. 
Then there exist a closed definable neighborhood $J$ of $C$ in $R$ and a definable triangulation $\pi_J:Y_J\to f^{-1}(J)$ 
of $f|_{f^{-1}(J)}$ such that $Y_J$ is a neighborhood of $f^{-1}(C)$ in $X$ and $\pi_J^{-1}(X_i)$ are polyhedra. 
If $X_i$ are compact polyhedra, we can choose $\pi_J$ so that $\pi^{-1}_J(X_i)\subset X_i$. \par
  It has been also shown that $f$ admits a definable triangulation outside of $f^{-1}(C)$ for some finite set $C\subset\Im f$. 
To be precise, for any compact definable subset $I$ of $\Im f-C$ there exists a definable triangulation $\pi_I:Y_I\to f^{-1}(I)$ 
of $f|_{f^{-1}(I)}$ such that $Y_I\subset R^n$, $\pi^{-1}_I(X_i)$ are polyhedra and if $X$ is a polyhedron then $Y_I\subset X$. 
\endproclaim
By theorem 2.2,(2) in the above special case and remark 2.3 we will prove theorem 2.1. 
We proceed by induction on $\dim X_1$ together with the following complement. 
\proclaim{Complement of theorem 2.1}
A definable homeomorphism $(X_1,X_{1,i})_{i=1,...,k}\allowmathbreak\to(X_2,X_{2,i})_{i=1,...,k}$ is definably isotopic to a PL 
homeomorphism $(X_1,X_{1,i})_i\to(X_2,X_{2,i})_i$ through homeomorphisms. 
\endproclaim
Let theorem 2.1$_m$ denote theorem 2.1 for $X_1$ of dimension $\le m$, and define complement$_m$ in the same way. 
\demo{Proof that theorem 2.1$_m$ implies complement$_m$}
Assume $\dim X_1\le m$ and theorem 2.1$_m$ holds. 
Let $\pi:(X_1,X_{1,i})_i\to(X_2,X_{2,i})_i$ be a definable homeomorphism, and let $K_1$ and $K_2$ be simplicial decompositions 
of $X_1$ and $X_2$ compatible with $\{X_{1,i}\}$ and $\{X_{2,i}\}$ respectively. 
First we reduce the problem to the case of $\{X_{1,i}\}=K_1$. 
By the triangulation theorem of definable sets we have a simplicial subdivision $K'_2$ of $K_2$ and a definable homeomorphism $\tau$ of $X_2$ 
preserving $K_2$ such that $\{\tau(\sigma'):\sigma'\in K'_2\}$ is compatible with $\{\pi(\sigma_1):\sigma_1\in K_1\}$. 
Then $\tau^{-1}(\sigma_2)$ and $\tau^{-1}\circ\pi(\sigma_1)$ are polyhedra for $\sigma_2\in K_2$ and $\sigma_1\in K_1$. 
Consider two definable homeomorphisms $\tau^{-1}:(X_2,\sigma_2)_{\sigma_2\in K_2}\to(X_2,\tau^{-1}(\sigma_2))_{\sigma_2\in K_2}$ 
and $\tau^{-1}\circ\pi:(X_1,\sigma_1)_{\sigma_1\in K_1}\to(X_2,\tau^{-1}\circ\pi(\sigma_1))_{\sigma_1\in K_1}$. 
Then it suffices to see they are definably isotopic to PL homeomorphisms through homeomorphisms. 
Hence we can assume from the beginning $\{X_{1,i}\}_i=K_1$. \par
  Let $r\in\N$. 
We will construct a definable isotopy $\Pi_r:|K_1^r|\times[0,\,1]\to\pi(|K_1^r|)$ of $\pi|_{|K_1^r|}$ through homeomorphisms by 
induction on $r$ so that $\Pi_r|_{|K_1^r|\times\{1\}}$ is PL and $\Pi_r(\sigma_1\times[0,\,1])=\pi(\sigma_1)$ for $\sigma_1\in 
K_1^r$. 
We define naturally $\Pi_0$ by $\Pi_0(x,t)=x$ for $(x,t)\in|K_1^0|\times[0,\,1]$. 
Assume we have $\Pi_r$ and let $\sigma_1\in K_1^{r+1}-K_1^r$. 
Then it suffices to extend $\Pi_r|_{\partial\sigma_1\times[0,\,1]}$ to $\sigma_1\times[0,\,1]$. 
Assume $0\in\Int\sigma_1$. 
By theorem 2.1$_m$ there exists a PL homeomorphism $\theta:\pi(\sigma_1)\to\sigma_1$. 
Then $\theta\circ\Pi_r|_{\partial\sigma_1\times[0,\,1]}$ is a definable isotopy of $\theta\circ\pi|_{\partial\sigma_1}$ through 
homeomorphisms and $\theta\circ\Pi_r|_{\partial\sigma_1\times\{1\}}$ is PL. 
Set 
$$
\Theta(x,t)=\cases
\theta\circ\Pi_r(x,t)\quad&\text{for}\ (x,t)\in\partial\sigma_1\times[0,\,1]\\
\theta\circ\pi(x)\qquad&\text{for}\ (x,t)\in\sigma_1\times\{0\},
\endcases
$$
and extend it to $\sigma_1\times[0,\,1]$ by the Alexander trick, i.e., so that $\Theta(0,1)=0$ and $\Theta$ is linear on the 
segment joining $(0,1)$ to each point of $\partial\sigma_1\times[0,\,1]\cup\sigma_1\times\{0\}$. 
Then $\Theta$ is an extension of $\theta\circ\Pi_r|_{\partial\sigma_1\times[0,\,1]}$, a definable isotopy of $\theta\circ\pi
|_{\sigma_1}$ through homeomorphisms and of class PL on $\sigma_1\times\{1\}$. 
Hence, if we define $\Pi_{r+1}$ on $\sigma_1\times[0,\,1]$ to be $\theta^{-1}\circ\Theta$ then $\Pi_{r+1}$ fulfills the 
requirements. 
Thus induction process works and complement$_m$ is proved. 
\qed
\enddemo
\demo{Proof of theorem 2.1}
We prove theorem 2.1 by induction on $m=\dim X_1$. 
It is trivial if $m=0$. 
Hence assume $m>0$ and theorem 2.1 and the complement hold for smaller dimensional $X_1$. 
Let $\pi:(X_1,X_{1,i})_i\to(X_2,X_{2,i})_i$ be a definable homeomorphism. \par
  First we prove that $X_1$ is a PL manifold with boundary if so is $X_2$. 
Assume $X_2$ is a PL manifold with boundary. 
Note then $X_1$ is a polyhedral definable $C^0$ manifold with boundary. 
Let $x\in X_1$ and $f$ a non-negative PL function on $X_1$ with zero set $\{x\}$. 
Then it suffices to see $f^{-1}(\epsilon)$ is a PL ball or sphere of dimension $m-1$ for small $\epsilon>0\in R$ by invariance 
of links (Lemma 2.19 in [R-S]). 
Apply the former half of remark 2.3 to a definable $C^0$ function $f\circ\pi^{-1}$ on $X_2$, a finite subset $\{0\}$ of $R$ and compact polyhedra $(X_2,\pi(x))$. 
Then there exists a compact polyhedral neighborhood $U$ of $\pi(x)$ in $X_2$ and a definable $C^0$ imbedding $\tau:U\to X_2$ such that 
$\tau\circ\pi(x)=\pi(x)$, $\tau(U)$ is a neighborhood of $\pi(x)$ in $X_2$ and $f\circ\pi^{-1}\circ\tau$ is PL. 
Hence $(f\circ\pi^{-1}\circ\tau)^{-1}(\epsilon)$ is a PL ball or sphere of dimension $m-1$ for small $\epsilon>0\in R$ since a 
neighborhood of $\pi(x)$ in $U$ is a PL manifold and $f\circ\pi^{-1}\circ\tau$ is a non-negative PL function with zero set $\{\pi
(x)\}$. 
Choose $\epsilon$ so small that $\tau(U)$ includes $(f\circ\pi^{-1})^{-1}(\epsilon)$. 
Then $f^{-1}(\epsilon)$ and $(f\circ\pi^{-1}\circ\tau)^{-1}(\epsilon)$ are definably homeomorphic and, hence, PL homeomorphic by 
induction hypothesis. 
Therefore $f^{-1}(\epsilon)$ is a PL ball or sphere of dimension $m-1$. \par
  Next we reduce the problem to the case where $X_1$ is included in $R^m$ and $X_2$ is a simplex of dimension $m$. 
Let $K_1$ be a simplicial decomposition of $(X_1,X_{1,i})_i$, i.e., a simplicial decomposition of $X_1$ such that each $X_{1,i}$ 
is the union of some simplexes in $K_1$. 
Note each $X_{2,i}$ is the union of some $\pi(\sigma_1),\ \sigma_1\in K_1$. 
Apply the triangulation theorem of definable sets to $(X_2,\pi(\sigma_1))_{\sigma_1\in K_1}$. 
Then we have a simplicial decomposition $K_2$ of $(X_2,X_{2,i})_i$ and a definable homeomorphism $\tau$ of $X_2$ such that for each $\sigma_1\in K_1$, 
$\pi(\sigma_1)$ is the union of some $\tau(\sigma_2),\ \sigma_2\in K_2$. 
Replace $\pi$ with $\tau^{-1}\circ\pi$. 
Then we can assume from the beginning there is a simplicial decomposition $K_2$ of $(X_2,X_{2,i})_i$ such that for each $\sigma_1\in K_1$, $\pi(\sigma_1)$ is the union of some $\sigma_2\in K_2$. 
Moreover repeating the same arguments on $K_1$ we assume there exists a simplicial subdivision $K'_1$ of $K_1$ such that for each $\sigma_2\in K_2$, 
$\pi^{-1}(\sigma_2)$ is the union of some $\sigma'_1\in K'_1$. 
Then it suffices to see each $\pi^{-1}(\sigma_2)$ is PL homeomorphic to $\sigma_2$ by the following reason. \par
  Assume we have proved it and let $\alpha_{\sigma_2}:\sigma_2\to\pi^{-1}(\sigma_2)$ be PL homeomorphisms. 
We will construct a definable isotopy $\Pi:X_1\times[0,\,1]\to X_2$ through homeomorphisms so that $\Pi(\pi^{-1}(\sigma_2)\times
[0,\,1])=\sigma_2$ for $\sigma_2\in K_2$ and $\Pi|_{X_1\times\{1\}}$ is PL. 
Let $r\in\N$ and suppose by induction on $r$ we have defined $\Pi$ on $\pi^{-1}(|K_2^r|)\times[0,\,1]$. 
Let $\sigma_2\in K_2^{r+1}-K_2^r$ and fix a point $a$ in $\Int\sigma_2$. 
Set 
$$
\Pi_{\sigma_2}(x,t)=\cases
\Pi(\alpha_{\sigma_2}(x),t)\quad&\text{for}\ (x,t)\in\partial\sigma_2\times[0,\,1],\\
(\pi\circ\alpha_{\sigma_2}(x),0)\quad&\text{for}\ (x,t)\in\sigma_2\times\{0\}. 
\endcases
$$
Then $\Pi_{\sigma_2}$ is a definable homeomorphism of $\partial\sigma_2\times[0,\,1]\cup\sigma_2\times\{0\}$ and $\Pi_{\sigma_2}$ 
is PL on $\partial\sigma_2\times\{1\}$. 
Extend it to a definable homeomorphism of $\sigma_2\times[0,\,1]$ by the Alexander trick, i.e., so that $\Pi_{\sigma_2}(a,1)=(a,1)
$ and $\Pi_{\sigma_2}$ is linear on the segment joining $(a,1)$ to each point of $\partial\sigma_2\times[0,\,1]\cup\sigma_2\times
\{0\}$. 
Then $\Pi_{\sigma_2}$ is PL on $\sigma_2\times\{1\}$. 
Hence if we define $\Pi$ on $\pi^{-1}(\sigma_2)\times[0,\,1]$ by $(\Pi(x,t),t)=\Pi_{\sigma_2}(\alpha_{\sigma_2}^{-1}(x),t)$ then 
$\pi$ fulfills the requirements. \par
  By the same reason it also suffices to prove $\pi(\sigma_1)$ is PL homeomorphic to $\sigma_1$ for each $\sigma_1\in K_1$ and by 
induction hypothesis we consider only $\sigma_1$ of dimension $m$ and $\pi|_{\sigma_1}:\sigma_1\to\pi(\sigma_1)$. 
Namely we assume $X_1=\sigma_1$. 
Then we can replace the ambient space $R^n$ of $X_1$ with the linear space spanned by $\sigma_1$, which is of dimension $m$, and 
assume $X_1\subset R^m$. 
Under this assumption we will prove $\pi^{-1}(\sigma_2)$ is PL homeomorphic to $\sigma_2$ for each $\sigma_2\in K_2$. 
By induction hypothesis $\pi^{-1}(\sigma_2)$ is PL homeomorphic to $\sigma_2$ if $\dim\sigma_2<m$. 
Hence we consider the case of $\dim\sigma_2=m$. 
Then since $\pi^{-1}(\sigma_2)$ is included in $R^m$, what we prove is the following statement. \par
  {\it Statement.} Let $X$ be a compact polyhedron in $R^m$ which is definably homeomorphic to a simplex $\sigma$ of dimension 
$m$. 
Then $X$ is PL homeomorphic to $\sigma$. \par
  {\it Proof of the statement.} As shown already $X$ is a PL manifold with boundary and, hence, $\partial X$ is a PL manifold. 
Let $c_X:\partial X\times[0,\,1]\to X$ and $c_\sigma:\partial \sigma\times[0,\,1]\to \sigma$ be PL collars on $\partial X$ 
in $X$ and on $\partial\sigma$ in $\sigma$, respectively (p.~24, 25 and 26 in [R-S]), i.e., PL imbeddings such that $c_X(\cdot,
0)=c_\sigma(\cdot,0)=\id$ and $\Im c_X$ and $\Im c_\sigma$ are neighborhoods of $\partial X$ in $X$ and of $\partial\sigma$ in 
$\sigma$ respectively. 
Then $X$ and $\sigma$ are PL homeomorphic to $X-c_X(\partial X\times[0,\,1/2))$ and $\sigma-c_\sigma(\partial\sigma\times[0,\,
1/2))$ respectively. 
Hence we have a definable homeomorphism from $X-c_X(\partial X\times[0,\,1/2))$ to $\sigma-c_\sigma(\partial\sigma\times[0,\,1/2
))$. 
The homeomorphism is extensible to a definable homeomorphism $\pi$ from $X$ to $\sigma$ which carries each $c_X(\partial 
X\times\{t\})$, $t\in[0,\,1/2]$, to $c_\sigma(\partial\sigma\times\{t\})$ and such that $\pi|_{c_X(\partial X\times[0,\,1/4])}$ 
is of class PL because a definable homeomorphism from $\partial X$ to $\partial\sigma$ is definably isotopic to a PL 
homeomorphism from $\partial X$ to $\partial\sigma$ through homeomorphisms by induction hypothesis (complement$_m$). \par
  Let $0\in\Int\sigma$, and define a non-negative PL function $g$ on $\sigma$ and a non-negative definable $C^0$ function $f$ on 
$X$ by 
$$
g(tx)=1-t\quad\text{for}\ (x,t)\in\partial\sigma\times[0,\,1],
$$
and $f=g\circ\pi$. 
Then $g^{-1}(0)=\partial\sigma,\ g^{-1}(1)=\{0\},\ f^{-1}(0)=\partial X$, $g^{-1}([t,\,1])$ is PL 
homeomorphic to $\sigma$ for each $t\in[0,\,1)$, and $f$ is PL on a compact polyhedral neighborhood of $\partial X$ in $X$. \par
  Apply theorem 2.2,(2) in the special case to $X$ and $f$. 
Then there exist a definable homeomorphism $\xi$ of $X$ such that $f\circ\xi$ is PL. 
Let $K$ be a simplicial decomposition of $X$ such that $f\circ\xi$ is linear on each simplex in $K$. 
Since $\dim X=m>0$, $f\circ\xi(K^0)$ consists of at least two numbers. 
Set $l=\#f\circ\xi(K^0)$. 
We will prove that $X$ is PL homeomorphic to $\sigma$ by induction on $l$. \par
  Assume $l=2$. 
Then $K^0-\partial X$ consists of one point $\pi^{-1}(0)$. 
Hence $X$ is the cone with base $\partial X$ and vertex $\pi^{-1}(0)$. 
Therefore $X$ is PL homeomorphic to $\sigma$. 
Assume $l>2$ and set $t_0=\min(f\circ\xi(K^0)-\{0\})$. 
Then $t_0<1$, and $(f\circ\xi)^{-1}([0,\,t_0/2])$ is a regular neighborhood of $\partial X$ in $X$ and, hence, a PL collar 
(Corollary 3.9 in [R-S]). 
(Here we give the definition of a regular neighborhood. 
Let $Y$ be a compact subpolyhedron of a compact polyhedron $X$. 
A {\it regular neighborhood} $U$ of $Y$ in $X$ is a compact polyhedral neighborhood of $Y$ in $X$ such that there exist compact polyhedra $X_1
\supset U_1\supset Y_1$, a PL homeomorphism $h:(X,U,Y)\to(X_1,U_1,Y_1)$ and a simplicial decomposition $K_1$ of $X_1$ such 
that $K_1|_{Y_1}$ is a full subcomplex of $K_1$ and $U_1=\cup\{|\st(\sigma,K'_1)|:\sigma\in K'_1|_{Y_1}\}$, where $K'_1$ denotes 
the barycentric subdivision of $K_1$.) 
Hence $(f\circ\xi)^{-1}([0,\,t_0/2])$ is PL homeomorphic to $\partial X\times[0,\,1]$. 
On the other hand, since $g^{-1}([0,\,t_0])$ is a PL manifold with boundary, $(f\circ\xi)^{-1}([0,\,t_0])$ is a polyhedral 
definable $C^0$ manifold with boundary. 
Hence $(f\circ\xi)^{-1}([0,\,t_0])$ is a PL manifold with boundary $(f\circ\xi)^{-1}(t_0)\cup\partial X$. 
Therefore $(f\circ\xi)^{-1}([t_0/2,\,t_0])$ is PL homeomorphic to $(f\circ\xi)^{-1}(t_0)\times[0,\,1]$ by the same reason as 
above. 
Here $(f\circ\xi)^{-1}(t_0)$ is PL homeomorphic to $\partial\sigma$. 
Hence $(f\circ\xi)^{-1}([0,\,t_0])$ is PL homeomorphic to $\partial X\times[0,\,1]$ by elementary arguments of PL topology. 
Remember that $g^{-1}([t_0,\,1])$ is PL homeomorphic to $\sigma$. 
Then by induction hypothesis on $l$, $(f\circ\xi)^{-1}([t_0,\,1])$ is PL homeomorphic to $\sigma$. 
Consequently $X$ is PL homeomorphic to $\sigma$, which proves the statement, theorem 2.1 and, hence, its complement. 
\qed
\enddemo
\demo{Continued proof of theorem 2.2} {\it (1)} 
By the latter half of remark 2.3 $f$ admits a definable triangulation outside of $f^{-1}(C)$ for some finite set $C\subset\Im f$. 
On the other hand, by the former half of remark 2.3 we have a closed definable neighborhood $J$ of $C$ in $\Im f$ and a definable triangulation 
$\pi_J:Y_J\to f^{-1}(J)$ of $f|_{f^{-1}(J)}$ such that $\pi^{-1}_J(X_j)$ are polyhedra. 
Set $I=\overline{\Im f-J}$. 
Note $I\cap J$ is a finite set. 
Apply the latter half to $f|_{f^{-1}(I)}$. 
Then there exists a definable triangulation $\pi_I:Y_I\to f^{-1}(I)$ 
of $f|_{f^{-1}(I)}$ such that $\pi^{-1}_I(X_i)$ are polyhedra. 
\par
  We will paste $Y_I$ and $Y_J$ at $(f\circ\pi_I)^{-1}(I\cap J)$ and $(f\circ\pi_J)^{-1}(I\cap J)$ respectively. 
Set $Y_{\partial I}=(f\circ\pi_I)^{-1}(I\cap J),\ \pi_{\partial I}=\pi_I|_{Y_{\partial I}},\ Y_{\partial J}=(f\circ\pi_J)^{-1}
(I\cap J)$ and $\pi_{\partial J}=\pi_J|_{Y_{\partial J}}$. 
Then $\pi_{\partial J}^{-1}\circ\pi_{\partial I}:Y_{\partial I}\to Y_{\partial J}$ is a definable homeomorphism. 
If it is of class PL, pasting $Y_I$ and $Y_J$ by the PL homeomorphism we obtain a polyhedron $Y$. 
Define a definable homeomorphism $\pi:Y\to X$ to be $\pi_I$ on $Y_I$ and $\pi_J$ on $Y_J$. 
Then $\pi:Y\to X$ is a definable triangulation of $f$ and $\pi^{-1}(X_i)$ are polyhedra. 
Hence it suffices to modify $\pi_I:Y_I\to X$ and $\pi_J:Y_J\to X$ so that $\pi_{\partial J}^{-1}\circ\pi_{\partial I}$ is PL. \par
  Let $K_I$ and $K_J$ be simplicial decompositions of $Y_I$ and $Y_J$, respectively, such that $f\circ\pi_I$ and $f\circ\pi_J$ 
are linear on each simplex in $K_I$ and $K_J$, respectively, and each $\pi^{-1}_I(X_i)$ and $\pi^{-1}_J(X_i)$ are the unions of 
some simplexes in $K_I$ and $K_J$ respectively. 
Set $K_{\partial I}=K_I|_{Y_{\partial I}}$ and $K_{\partial J}=K_J|_{Y_{\partial J}}$. 
First we modify $\pi_J:Y_J\to X$. 
By the triangulation theorem of definable sets there exists a definable homeomorphism $\theta$ of $Y_{\partial J}$ preserving $K_{\partial J}$ 
such that $\theta\circ\pi^{-1}_{\partial J}\circ\pi_{\partial I}(\sigma),\ \sigma\in K_{\partial I}$, are polyhedra. 
Extend $\theta$ to a definable homeomorphism $\Theta$ of $Y_J$ by the Alexander trick. 
Then $f\circ\pi_J\circ\Theta=f\circ\pi_J$ and $\Theta$ is preserving $K_J$ and, hence, $\{\pi_J^{-1}(X_i)\}_i$. 
Hence we can replace $\pi_J$ with $\pi_J\circ\Theta^{-1}$. 
If we replace, we can assume from the beginning that $\pi^{-1}_{\partial J}\circ\pi_{\partial I}(\sigma)$ are polyhedra for $\sigma\in K_{\partial 
I}$. \par
  Secondly, we modify $\pi_I:Y_I\to X$. 
Apply theorem 2.1 to the definable homeomorphism $\pi^{-1}_{\partial J}\circ\pi_{\partial I}:(Y_{\partial I},K_{\partial I})\to(
Y_{\partial J},\pi^{-1}_{\partial J}\circ\pi_{\partial I}(\sigma))_{\sigma\in K_{\partial I}}$. 
Then we have a PL homeomorphism $\xi:(Y_{\partial I},K_{\partial I})\to(Y_{\partial J},\pi^{-1}_{\partial J}\circ\pi_{\partial I
}(\sigma))_{\sigma\in K_{\partial I}}$ such that $\xi(\sigma)=\pi^{-1}_{\partial J}\circ\pi_{\partial I}(\sigma)$ for each 
$\sigma\in K_{\partial I}$. 
Repeat the same arguments as the above first modification to the definable homeomorphism $\pi^{-1}_{\partial I}\circ\pi_{
\partial J}\circ\xi:Y_{\partial I}\to Y_{\partial I}$, which is preserving $K_{\partial I}$. 
Then we have its definable homeomorphism extension $\Xi:Y_I\to Y_I$ preserving $K_I$ such that $f\circ\pi_I\circ\Xi=f\circ\pi_
I$. 
Hence we can replace $\pi_I$ with $\pi_I\circ\Xi$. 
Carry out the replacement. 
Then, since $\pi^{-1}_J\circ(\pi_I\circ\Xi)=\pi^{-1}_J\circ\pi_I\circ\pi^{-1}_{\partial I}\circ\pi_{\partial J}\circ\xi=\xi$ on 
$Y_{\partial I}$, $\pi^{-1}_{\partial J}\circ\pi_{\partial I}$ becomes of class PL. 
Thus theorem 2.2,(1) is proved. \par
  {\it (2)} 
Let $X$ be the underlying polyhedron to a finite simplicial complex $P$, $f$ a definable $C^0$ function on $X$, and $X_i$ a 
finite number of compact definable subsets of $X$. 
As above we have a definable triangulation $\pi:Y\to X$ of $f$ such that $\pi^{-1}(X_i)$ and $\pi^{-1}(\sigma),\ \sigma\in P$, 
are polyhedra. 
Let $\theta:(X,\sigma)_{\sigma\in P}\to(Y,\pi^{-1}(\sigma))_{\sigma\in P}$ be a PL homeomorphism such that $\theta(\sigma)=\pi
^{-1}(\sigma)$ for $\sigma\in P$, which exists by theorem 2.1. 
Then $\pi\circ\theta:X\to X$ is a definable triangulation of $f$ such that $\pi\circ\theta(\sigma)=\sigma$ for $\sigma\in P$ and 
$(\pi\circ\theta)^{-1}(X_i)$ are polyhedra. 
By the complement of theorem 2.1 we can choose $\theta$ so that there exists a definable isotopy from id to $\pi$ through 
homeomorphisms and preserving $P$, which proves theorem 2.2,(2). 
\qed
\enddemo
\demo{Proof of the complement of theorem 2.2}
Lemma II.3.10 states that there exists a PL homeomorphism $\omega:Y'\to Y$ such that $f\circ\pi\circ\omega=f\circ\pi'$ 
and $\omega(\pi^{\prime-1}(X_i))=\pi^{-1}(X_i)$ for each $i$ in the real number case. 
Its proof uses only theorem 2.1, its complement and theorems 2.2,(1) and 2.2,(2), and works 
for general $R$. 
In the proof we obtain $\omega$ as the finishing homeomorphism of some definable isotopy $\omega_t:Y'\to Y,\ 0\le t\le 1$, 
through homeomorphisms such that $\omega_0=\pi^{-1}\circ\pi'$, $f\circ\pi\circ\omega_t=f\circ\pi'$ for $t\in[0,\,1]$ 
and $\omega_t(\pi^{\prime-1}(X_i))=\pi^{-1}(X_i)$ for each $i$ and $t$. 
Thus the complement holds. 
We do not repeat the proof. 
\qed
\enddemo
\demo{Proof of theorem 2.2,(3)}
Assume $X$ is the underlying polyhedron to a finite simplicial complex $P$ and $f$ is simplicial on $P$. 
Let $\tau:X\to X$ be a definable triangulation of $f$ preserving $P$ such that $\tau^{-1}(X_i)$ are polyhedra. 
Regard $\id:X\to X$ as a definable triangulation of $f$ preserving $P$, and apply the complement of theorem 2.2 to $\tau$, id and $P$. 
Then there exists a definable isotopy $\omega_t,\ 0\le t\le 1$, of $\tau^{-1}$ through homeomorphisms such that $\omega_1$ is PL, $f\circ\tau\circ\omega_t=f$ for $t\in[0,\,1]$ and $\omega_t(\sigma)=\tau^{-1}(\sigma)=\sigma$ for each $\sigma\in P$ and $t$. 
Set $\pi_t=\tau\circ\omega_t$ for $t\in[0,\,1]$. 
Then $\pi_0=\id$, $\pi^{-1}_t(\sigma)=\omega^{-1}_t(\tau^{-1}(\sigma))=\sigma$ for $\sigma\in P$, hence, $\pi_t,\ 0\le t\le 1$, is a definable isotopy of $X$ preserving $P$, $\pi^{-1}_1(X_i)=\omega_1^{-1}(\tau^{-1}(X_i))$ are polyhedra, and $f\circ\pi_t=f\circ\tau\circ\omega_t=f$ for $t\in[0,\,1]$. 
Therefore it suffices to modify $\pi_t$ so that $\pi_t=\id$ on $f^{-1}(R-(s'_1,\,s'_2))$. \par
  Let $P'$ be a simplicial subdivision of $P$ such that $\pi^{-1}_1(X_i)$ are the unions of some simplexes of $P'$. 
Choose $s''_1<s'''_1<s'''_2<s''_2$ in $R$ so that $s'_1<s''_1<s'''_1<s'''_2<s''_2<s'_2$ and $f(P^{\prime 0})\cap([s''_1,\,s'''_1]\cup[s'''_2,\,s''_2])=\emptyset$, and set 
$$
Q=\{\sigma_P\cap f^{-1}(\sigma):\sigma_P\in P',\,\sigma\in\{s''_1,s'''_1,s''_2,s'''_2,[s''_1,\,s'''_1],[s'''_2,\,s''_2]\}\}.
$$
Then $Q$ is a cell complex and $f|_{|Q|}:Q\to\{s''_1,s'''_1,s''_2,s'''_2,[s''_1,\,s'''_1],[s'''_2,\,s''_2]\}$ is a trivial cellular map. 
Hence by the Alexander trick there exists a definable isotopy $\pi'_t,\ 0\le t\le 1$, of $X$ such that for any $t\in[0,\,1]$, $\pi'_t=\pi_t$ on $f^{-1}([s'''_1,\,s'''_2])$, $\pi'_t=\id$ on $f^{-1}(R-(s''_1,\,s''_2))$, $f\circ\pi'_t=f$ on $f^{-1}(R-(s'''_1,\,s'''_2))$ and $\pi_t|_{|Q|}$ is preserving $Q$. 
Then $\pi'_t,\ 0\le t\le 1$, clearly fulfills requirements except that $\pi^{\prime-1}_1(X_i)$ are polyhedra. 
Divide each $\pi^{\prime-1}_1(X_i)$ to three $\pi^{\prime-1}_1(X_i)\cap f^{-1}([s'''_1,\,s'''_2]),\ \pi^{\prime-1}_1(X_i)\cap|Q|$ and $\pi^{\prime-1}_1(X_i)\cap f^{-1}(R-(s''_1,\,s''_2))$. 
The first set is $\pi^{-1}_1(X_i)\cap f^{-1}([s'''_1,\,s'''_2])$, the second the union of some cells of $Q$, and the third $X_i\cap f^{-1}(R-(s''_1,\,s''_2))$. 
These are all polyhedra. 
Hence $\pi^{\prime-1}_1(X_i)$ is a polyhedron. 
Thus theorem 2.2,(3) is proved. 
\qed
\enddemo
Consider the non-compact case of theorems 2.1 and 2.2. 
We use the notion {\bf semi-linear}. 
A {\it semi-linear} set or map is a definable one in the o-minimal structure $(R,<,0,+,c\cdot:c\in R)$, where $c\cdot$ denote 
the function $R\ni x\to c x\in R$. 
(For simplicity of notation in the following proofs we use this o-minimal structure for the definition of a semi-linear set or map. 
The results below, however, hold clearly in the o-minimal structure $(R,<,0,+)$.) 
Let $(X,X_i)_i$ be a definable set and a finite number of definable subsets. 
Assume $X$ is bounded in $R^n$ and consider a definable triangulation of $(\overline X,X,X_i)_i$. 
Then it follows that there exist semi-linear sets $Y_i\subset Y$ in $R^n$ and a definable homeomorphism $\pi:(Y,Y_i)_i\to(X,X_i
)_i$. 
We call $\pi:(Y,Y_i)_i\to(X,X_i)_i$ a {\it definable semi-linearization} of $(X,X_i)_i$. 
We wish to find conditions under which uniqueness of a definable semi-linearization holds. 
First, since $(0,\,1)$ and $R$ are definably homeomorphic but not semi-linearly homeomorphic, we treat only bounded semi-linear 
subsets of $R^n$. 
Boundedness condition is, however, not sufficient. 
For example, let $\sigma$ be a 2-simplex in $R^2$ with $0\in\Int\sigma$. 
Then $\sigma-\{0\}$ and $\sigma-\sigma/2$ are definably homeomorphic but not semi-linearly homeomorphic. 
We treat semi-linear sets of the latter form. 
To be precise, a semi-linear subset $Y$ of $R^n$ with a finite number of semi-linear subsets $Y_i$ is called {\it standard} if 
$Y$ is bounded and each point of $\overline Y-Y$ has a semi-linear neighborhood $U$ in $\overline Y$ such that $(U,U\cap Y,U
\cap Y_i)_i$ is semi-linearly homeomorphic to $(V\times[0,\,1),V\times(0,\,1),V_i\times(0,\,1))_i$ for some semi-linear sets 
$(V,V_i)_i$. 
We call $U$, or the semi-linear homeomorphism, a {\it local collar} on $\overline Y-Y$ in $(\overline Y,Y,Y_i)_i$ (see p.~24 in 
[R-S]). 
As Theorem 2.25 in [R-S], which states the existence of local collars implies that of a global collar, we see the semi-linear case. 
\proclaim{Lemma 2.4}
Let $(Y,Y_i)_i$ be a standard family of semi-linear sets. 
Then there exist a family of bounded semi-linear sets $(V,V_i)_i$ and a semi-linear homeomorphism of $(V\times[0,\,1),V\times(0,\,1),V_i\times
(0,\,1))_i$ to $(U,U\cap Y,U\cap Y_i)_i$ for some semi-linear neighborhood $U$ of $\overline Y-Y$ in $\overline Y$. 
\endproclaim
We call $U$, or the homeomorphism, a {\it collar} on $\overline Y-Y$ in $(\overline Y,Y,Y_i)_i$. 
\demo{Proof}
We proceed as in the proof of Theorem 2.25 in [R-S]. 
Let $K$ be a simplicial decomposition of $\overline Y$ such that $Y$ and $Y_i$ are the unions of some open simplexes in $K$, $K'$ and 
$K''$ the barycentric subdivisions of $K$ and $K'$, respectively, and $v_\sigma$ the barycenter of each $\sigma\in K'$, set $V=\overline Y-Y$, and 
let $Y\subset R^n$. 
Identify $Y$ with $Y\times\{0\}\subset Y\times R$, set $W=\overline Y\times\{0\}\cup\overline V\times[0,\,2]$. 
For each $\sigma\in K'$ with $\Int\sigma\subset V$, let $U_\sigma$ be a local collar at $v_\sigma$ on $V$ in $(\overline Y,Y,Y_i
)_i$, and $\phi_\sigma:K''\to\{0,1,[0,\,1]\}$ the simplicial map such that $\phi_\sigma=1$ at $v_\sigma$ and $\phi_\sigma=0$ at 
the other vertices. 
For simplicity of notation, assume $v_\sigma=0$, and choose $\epsilon>0\in R$ so that $\epsilon|\st(v_\sigma,K'')|\subset\Int U_
\sigma$. 
Then by properties of a local collar, there exists a semi-linear homeomorphism $h_{\epsilon\sigma}$ of $W$ such that 
$$
\gather
h_{\epsilon\sigma}(h_{\epsilon\sigma}^{-1}(Y_i\times\{0\}))\subset Y_i\times\{0\},\quad h_{\epsilon\sigma}(v_\sigma,1)=(v_\sigma,
0),\\
h_{\epsilon\sigma}=\id\quad\text{on}\ W-(V\cap\Int\epsilon|\st(v_\sigma,K'')|)\times[0,\,2)-(\Int\epsilon|\st(v_\sigma,K'')|)
\times\{0\},
\endgather
$$
and $h_{\epsilon\sigma}$ carries linearly each segment joining $(v_\sigma,1)$ and a point of $\partial(V\cap\epsilon|\st(v_\sigma,
\allowmathbreak K'')|)\times[0,\,2]\cup(V\cap\epsilon|\st(v_\sigma,K'')|)\times\{2\}$. 
Note 
$$
h^{-1}_{\epsilon\sigma}(\overline V\times[0,\,2])-\{(v_\sigma,0)\}=\overline V\times[0,\,2]-(v_\sigma,1)*(V\cap\Int\epsilon|\st(
v_\sigma,K'')|). 
$$
Hence by the star property of $K''$ at $v_\sigma$ we have a semi-linear homeomorphism $h_\sigma$ of $W$ such that 
$$
\gather
h_\sigma(h_\sigma^{-1}(Y_i\times\{0\}))\subset Y_i\times\{0\},\quad h_\sigma(v_\sigma,1)=(v_\sigma,0),\\
h_\sigma=\id\quad\text{on}\ W-(V\cap\Int|\st(v_\sigma,K'')|)\times[0,\,2)-(\Int|\st(v_\sigma,K'')|)\times\{0\},
\endgather
$$
and $h_\sigma$ carries linearly each segment joining $(v_\sigma,1)$ and a point of $\partial(V\cap|\st(v_\sigma,K'')|)
\allowmathbreak\times[0,\,2]\cup(V\cap|\st(v_\sigma,K'')|)\times\{2\}$. 
Note
$$
\gather
h_\sigma(\{(y,t)\in \overline V\times R:t=\phi_\sigma(y)\})=\overline V\times\{0\}. \\
\{\sigma\in K':\Int\sigma\in V\}=\{\sigma_1,...,\sigma_k\},\quad\phi=\sum_{i=1}^k\phi_{\sigma_i},\tag"\ \ \ \ Set"\\
h=h_{\sigma_1}\circ\cdots\circ h_{\sigma_k},\quad Z=\overline Y\times\{0\}\cup\{(y,t)\in\overline V\times R:0\le t\le\phi(y)\}. 
\endgather
$$
Then $\phi$ is a simplicial map from $K''$ to $\{0,1,[0,\,1]\}$ with $\phi^{-1}(0)\cap K^{\prime\prime 0}=\{v_{\sigma_1},...,
\allowmathbreak v_{\sigma_k}\}$, 
$h|_Z$ is a semi-linear homeomorphism to $\overline Y\times\{0\}$ such that 
$$
h(h^{-1}(Y\times\{0\}))=Y\times\{0\}\quad\text{and}\ h(h^{-1}(Y_i\times\{0\}))=Y_i\times\{0\},
$$
and there exists a semi-linear homeomorphism $\tau:Z\to Y\times\{0\}\cup V\times[0,\,1]$ of the form $\tau(y,t)=(\tau'(y,t),t)$ 
for $(y,t)\in Z$ such that $\tau=\id$ on $Y\times\{0\}$ and $\tau(\sigma\times R\cap Z)\subset\sigma$ for $\sigma\in K'$. 
Therefore $\big(Z-\overline Y\times\{0\},h^{-1}(Y\times\{0\})-\overline Y\times\{0\},h^{-1}(Y_i\times\{0\})-\overline Y\times\{0\}\big)
_i$ and, hence, $\big(Z,h^{-1}(Y\times\{0\}),h^{-1}(Y_i\times\{0\})\big)_i$ admit a collar on $Z-h^{-1}(Y\times\{0\})$. 
Thus the lemma is proved. 
\qed
\enddemo
  A typical example of a standard semi-linear set is the interior of a compact PL manifold with boundary. 
We see that a definable $C^0$ manifold is definably homeomorphic to the interior of some compact PL manifold possibly with 
boundary and the compact manifold possibly with boundary is unique up to PL homeomorphisms in the same way as Theorem V.2.1 in 
[S$_1$]. 
More generally, 
\proclaim{Theorem 2.5}
A definable set with a finite number of definable subsets admits uniquely a standard semi-linearization. 
Moreover a definable homeomorphism between standard semi-linear set families is definably isotopic to a semi-linear 
homeomorphism through homeomorphisms. 
\endproclaim
\demo{Proof}
Let $(X,X_i)_i$ be the given definable sets. 
As shown above we can assume $X$ and $X_i$ are the unions of some open simplexes in some finite simplicial complex $K$ in $R^n$ 
with $|K|=\overline X$. 
Let $K'$ and $K''$ denote the barycentric subdivisions of $K$ and $K'$, respectively, and $v_\sigma$ the barycenter of $\sigma\in K$. 
(We require $K'$ to be only a simplicial subdivision of $K$ such that $K'|_{|K^k|},\ k=1,...,$ are full subcomplexes of $K'$.) 
Define a definable $C^0$ imbedding $\tau_\sigma:\overline X-\Int\sigma\to\overline X-\Int\sigma$ by $\tau_\sigma=\id$ outside of 
$\Int|\st(v_\sigma,K')|-\Int\sigma$ and 
$$
\gather
\tau_\sigma(t_0v_\sigma+\sum_{j=1}^l t_j v_j)=t_0v_\sigma/2+\sum_{j=1}^kt_j v_j+(\sum_{j=k+1}^lt_j+t_0/2)\sum_{j=k+1}^lt_j v_j/
\sum_{j=k+1}^lt_j\\
\qquad\text{for}\ t_0,...,t_l\in[0,\,1]\ \text{with}\ \sum_{j=0}^kt_j<1\ \text{and}\ \sum_{j=0}^lt_j=1,
\endgather
$$
where $v_1,...,v_k,v_\sigma,v_{k+1},...,v_l$ are both the vertices of a simplex in $\st(v_\sigma,K')$ and the respective barycenters 
of increasing simplexes $\sigma_1\subset\cdots\subset\sigma_k\subset\sigma\subset\sigma_{k+1}\subset\cdots\subset\sigma_l$ in $K$. 
Then $\tau_\sigma$ is preserving $K-\{\sigma\}$, 
$$
\gather
\overline X-\overline{\Im\tau_\sigma}=\cup\{\Int|\st(\sigma'',K'')|:\sigma''\in K'',\Int\sigma''\subset\Int\sigma\},\\
\Im\tau_{\sigma_1}\circ\tau_{\sigma_2}=\Im\tau_{\sigma_1}\cap\Im\tau_{\sigma_2}\quad\text{for}\ \sigma_1\not=\sigma_2\in K. 
\endgather
$$\par
  Set $\{\sigma_1,...,\sigma_k\}=\{\sigma\in K:X\cap\Int\sigma=\emptyset\}$ and $\tau=\tau_{\sigma_1}\circ\cdots\circ\tau_{\sigma
_k}$. 
Then $\tau$ is a semialgebraic $C^0$ imbedding of $(X,X_i)_i$ into $(X,X_i)_i$, $\overline{\Im\tau}\subset X$, 
$$
\overline X-\overline{\Im\tau}=\cup\{\Int|\st(\sigma'',K'')|:\sigma''\in K'',X\cap\Int\sigma''=\emptyset\}, 
$$
$(\Im\tau,\tau(X_i))_i$ is a family of semi-linear sets, and it it standard for the following reason. \par
  Set $V=\overline{\Im\tau}-\Im\tau$. 
Note $V$ is the union of some open simplexes in $K''$. 
Let the above $\tau_\sigma$ be rewritten as $\tau_{\sigma 1/2}$ and define $C^0$ imbeddings $\tau_{\sigma t}:\overline X-\Int
\sigma\to\overline X-\Int\sigma$ for $t\in(0,\,1]$ by replacing $t_0/2$ in the definition of $\tau_{\sigma 1/2}$ with $t_0t$. 
Define $\tau_t$ to be $\tau_{\sigma_1 t}\circ\cdots\circ\tau_{\sigma_k t}$ for $t\in(0,\,1]$. 
Then $\tau_t$ are semialgebraic $C^0$ 
imbeddings of $(X,X_i)_i$ into $(X,X_i)_i$, $\overline{\Im\tau_t}\subset\Im\tau_{t'}$ for $t<t'$, $\tau_1=\id$, and the map 
$\pi:V\times[1/2,\,1]\ni(x,t)\to\tau_t(x)\in\overline{\Im\tau_{1/2}}$ is a definable $C^0$ imbedding whose image---$\overline{
\Im\tau_{1/2}}-\Im\tau_{1/4}$---is a semi-linear neighborhood of $V$ in $\overline{\Im\tau_{1/2}}$ and such that 
$$
\gather
\pi^{-1}(\Im\tau_{1/2})=V\times[1/2,\,1),\\
\pi^{-1}(\tau_{1/2}(X_i))=(V\cap X_i)\times[1/2,\,1),\\
\pi^{-1}(\sigma)=(V\cap\sigma)\times[1/2,\,1]\quad\text{for}\ \sigma\in K'.
\endgather
$$
It also follows that $\pi(\cdot,t)$ is extended to a cellular map $\{\sigma\cap\overline V:\sigma\in 
K''\}\to\{\sigma\cap(\overline{\Im\tau_{t/2}}-\Im\tau_{t/2}\overline):\sigma\in K''\}$ for each $t\in[1/2,\,1]$, where $(\ \ \overline)$ denotes the 
closure of a set $(\ \ )$. 
Let $\sigma\in K''$ with $\sigma\cap\overline V\not=\emptyset$. 
If $\sigma\cap(\overline V-V)=\emptyset$, then $\sigma\cap\Im\pi$ is both a cell and the disjoint union of the cells $\sigma\cap\pi(V
\times\{t\}),\ t\in[1/2,1]$. 
If $\sigma\cap(\overline V-V)\not=\emptyset$, then $\sigma$ is described as $\sigma_1*\sigma_2$ for some $\sigma_1,\sigma_2\in K''$ with $\sigma_1
\subset\overline V-V$ and $\sigma_2\cap(\overline V-V)=\emptyset$, $\sigma\cap\overline{\pi(V\times\{t\})}$ is a cell and equal to $\sigma_1*(\sigma_2\cap\pi(V\times\{t\}))$ for each $t\in[1/2,\,1]$, and 
$$
\sigma\cap\overline{\pi(V\times\{t_1\})}\cap\overline{\pi(V\times\{t_2\})}=\sigma_1\quad\text{for}\ t_1\not= t_2\in[1/2,\,1], 
$$
where $\sigma_1*\sigma_2$ denotes the join of $\sigma_1$ and $\sigma_2$---the smallest cell including  $\sigma_1$ and $\sigma_2$ 
under the condition that any two distinct segments with ends in $\sigma_1$ and $\sigma_2$ do not meet except in $\sigma_1\cup\sigma_2$. 
In particular, $\sigma\cap\overline V=\sigma_1*(\sigma_2\cap\overline V)$ in the latter case. 
Hence there exists a cellular map $\theta:\{\sigma\cap\overline V:\sigma\in K''\}\to\{0,1,[0,\,1]\}$ such that $\theta=0$ on 
$K^{\prime\prime 0}\cap(\overline V-V)$ and $\theta=1$ at the other vertices. 
Set 
$$
\Xi=\{(x,t)\in\overline V\times[0,\,1]:0\le t\le\theta(x)\}. 
$$
Then we can construct a PL homeomorphism $\xi:\Xi\to\overline{\Im\pi}$ inductively on $\Xi\cap|K^{\prime\prime l}|\times[0,\,1
],\ l=0,1,...,$ by the Alexander trick so that $\xi(\cdot,0)=\id$ on $\overline V$, $\pi^{-1}\circ\xi(x,t)$ is of the form 
$(\pi'(x,t),1-t/(2\theta(x)))$ for $(x,t)\in\Xi\cap V\times[0,\,1]$ and 
$$
\xi(\Xi\cap\sigma\times[0,\,1])=\overline{\pi((V\cap\sigma)\times[1/2,\,1])}\quad\text{for}\ \sigma\in K'\ \text{(not}\ K'')\ 
\text{with}\ V\cap\Int\sigma\not=\emptyset. 
$$
Clearly there exists a semi-linear homeomorphism $\xi':\overline V\times[0,\,1]\to\Xi$ of the form $\xi'(x,t)=(\xi''(x,t),t)$ 
for $(x,t)\in\overline V\times[0,\,1]$ such that 
$$
\xi'((\overline V\cap\sigma)\times[0,\,1])=\Xi\cap(\overline V\cap\sigma)\times[0,\,1]\quad\text{for}\ \sigma\in K'\ \text{(not}
\ K''). 
$$
Hence $\xi\circ(\xi'|_{V\times[0,1]}):V\times[0,\,1]\to\overline{\Im\tau}$ is a collar on $V$ in $(\overline{\Im\tau},\Im\tau,
Y_i\cap\Im\tau)_i$. \par
  We have shown a standard semi-linearization of $(X,X_i)_i$. 
We prepare for proving uniqueness. \par
  Set $W=\overline X-X$ and $U_K=\cup\{\Int|\st(\sigma,K'')|:\sigma\in K'',\sigma\subset W\}$, let $L$ be another simplicial 
decomposition of $\overline X$ with the same properties as $K$ and such that $K|_{(\overline W-W\overline)}=L|_{(\overline W-W
\overline)}$, 
and define $L',L''$ and $U_L$ in the same way. 
Then $(U_K,U_K\cap X,U_K\cap X_i)_i$ and $(U_L,U_L\cap X,U_L\cap X_i)_i$ are semi-linearly homeomorphic for the following reason. 
This is called Regular Neighborhood Theorem in the case of compact $W$ (Theorem 3.24 in [R-S]). 
As the proof in our case is the same, we rapidly repeat the proof in [R-S]. 
We can assume $L$ is a subdivision of $K$ by replacing $L$ with a simplicial subdivision of the cell complex $\{\sigma_K\cap
\sigma_L:\sigma_K\in K,\sigma_L\in L\}$ without new vertices. 
Let $K'_1$ and $K''_1$ be any derived subdivisions of $K$ and $K'_1$, respectively, (see p.~20 in [R-S]) and define $U_{K_1}$ by $K''_1$ and 
$L'_1,\ L''_1$ and $U_{L_1}$ in the same way. 
Then by 3.6 in [R-S], $(U_K,U_K\cap X,U_K\cap X_i)_i$ and $(U_L,U_L\cap X,U_L\cap X_i)_i$ are semi-linearly homeomorphic to 
$(U_{K_1},U_{K_1}\cap X,U_{K_1}\cap X_i)_i$ and $(U_{L_1},U_{L_1}\cap X,U_{L_1}\cap X_i)_i$, respectively, and by Lemma 3.7 in 
[R-S] there exist $K_1$ and $L_1$ such that 
$$
(U_{K_1},U_{K_1}\cap X,U_{K_1}\cap X_i)_i=(U_{L_1},U_{L_1}\cap X,U_{L_1}\cap X_i)_i. 
$$
Hence $(U_K,U_K\cap X,U_K\cap X_i)_i$ and $(U_L,U_L\cap X,U_L\cap X_i)_i$ are semi-linearly homeomorphic. 
Note the homeomorphism is extensible to a semi-linear homeomorphism of $(\overline X,X,X_i)_i$. \par
  We use this invariance property to show that if $(X,X_i)_i$ is standard then $U_K$ is a collar on $W$ in $(\overline X,X,X_i)
_i$ as follows. 
We are interested in only a semi-linear neighborhood of $W$ in $\overline X$ and by the invariance property we can modify $X$ by a simplicial 
map fixing $(\overline W-W\overline)$. 
Hence we can replace $(X,X_i,W)_i$ with $\big(X\times(0,\,1]\cup W\times[0,\,1],X_i\times(0,\,1]\cup W\times[0,\,1],X\times\{0\}\big)_i$. 
Fix a simplicial decomposition $K$ of $\overline X$. 
Let $K'_\times$ be a simplicial subdivision of the cell complex $\{\sigma\times\{0\},\sigma\times\{1\},\sigma\times[0,\,1]:
\sigma\in K'\}$ without new vertices such that for $\sigma_1,\sigma_2\in K'$ with $\sigma_1*\sigma_2\in K',\ \sigma_1\subset W$ 
and $\sigma_2\subset X$, $(\sigma_1\times\{0\})*(\sigma_2\times[0,\,1])$ is the union of some simplexes in $K'_\times$. 
(Proposition 2.9 in [R-S] states a simplicial subdivision without new vertices only. 
Its proof, however, shows the above last additional condition may be satisfied.) 
Define $U_{K_\times}$ by $K'_\times$ as $U_K$ by $K'$. 
Then by the invariance property it suffices to see $U_{K_\times}$ is a collar on $X\times\{0\}$ in $\big(\overline X\times[0,\,1],
X\times(0,\,1]\cup W\times[0,\,1],X_i\times(0,\,1]\cup W\times[0,\,1]\big)_i$. 
That was already shown. 
Indeed, $\big(U_{K_\times},U_{K_\times}\cap X\times(0,\,1],U_{K_\times}\cap X_i\times(0,\,1]\big)_i$ is semi-linearly homeomorphic to 
$\big(\Phi,\Phi\cap X\times(0,\,1],\Phi\cap X_i\times(0,\,1]\big)_i$, where $\Phi=\{(x,t)\in X\times[0,\,1):0\le t\le\psi(x)\}$ and 
$\psi$ is the simplicial map from $K'$ to $\{0,1,[0,\,1]\}$ such that $\psi=0$ on $K^{\prime 0}\cap W$ and $\psi=1$ at the 
other vertices. \par
  We wish to replace $U_K$ with a set according to the following proof of uniqueness. 
Define a definable $C^0$ function $f_K$ on $\overline X-\cup\{\sigma\in K':\sigma\subset\overline W-W\}$ by 
$$
\gather
f_K=\cases
0\quad&\text{on}\ K^{\prime 0}\cap W\\
1\quad&\text{on}\ K^{\prime 0}-\overline W,
\endcases\\
f_K(\sum_{j=1}^{l_3}t_j v_j)=\sum_{j=l_1+1}^{l_2}t_j/\sum_{j=1}^{l_2}t_j 
\endgather
$$
for $t_1,...,t_{l_3}\in[0,\,1]$ with $\sum_{j=1}^{l_3} t_j=1$ and $\sum_{j=1}^{l_2}t_j>0$ and for a simplex in $K'$ with vertices 
$v_1,...,v_{l_1}$ in $W$, $v_{l_1+1},...,v_{l_2}$ outside of $\overline W$ and $v_{l_2+1},...,v_{l_3}$ in $\overline W-W$. 
Then by the above arguments $(\overline X,U_K,X,X_i)_i$ is semi-linearly homeomorphic to $(\overline X,f^{-1}_K([0,1/2)),X,
\allowmathbreak X_i)_i$ and, hence, we can replace $U_K$ with $f^{-1}_K([0,\,1/2))$. 
We keep the notation $U_K$. \par
  Let $K''_1$ and $K''_2$ be any simplicial subdivisions of the barycentric subdivision of $K'$, and set $U'_i=\cup\{\Int|\st(\sigma,K''_i)|:\sigma\in K''_i,\sigma\subset W\},\ i=1,2$. 
Then $U'_i$ are regular neighborhoods of $\cup\{\sigma\in K':\sigma\subset W\}$ in $\overline X$. 
Hence by Regular Neighborhood Theorem there exists a semi-linear isotopy $\beta_t,\ 0\le t\le 1$, of $\overline X$ preserving 
$K'$ such that $\beta_1(U'_1)=U'_2$. 
Therefore, if we define $f'_1$ and $f'_2$ from $K''_1$ and $K''_2$, respectively, as $f_K$ from $K'$ then the following note holds. \par
  {\it Note 1.} There exists a semi-linear isotopy $\gamma_t,\ 0\le t\le 1$, of $\overline X$ preserving $K'$ such that $\gamma_1
(\Dom f'_1)=\Dom f'_2$, where $\Dom$ denotes the domain of a map. \par
  {\it Uniqueness.} 
Now we begin to prove uniqueness. 
Let $\eta:(Y,Y_i)_i\to(Z,Z_i)_i$ be a definable homeomorphism between standard semi-linear set families in $R^n$. 
Note $\eta$ is not necessarily extensible to a definable $C^0$ map from $\overline Y$ to $\overline Z$, but if $\eta$ is 
semi-linear then $\eta$ is extensible to a semi-linear homeomorphism from $\overline Y$ to $\overline Z$ because a semi-linear 
$C^0$ imbedding of an open simplex to $R^n$ is extensible to a semi-linear $C^0$ immersion of the closure. 
First we reduce the problem to the case where $\eta$ is extensible to a definable $C^0$ map $\overline\eta:\overline Y\to\overline Z$. 
Set $X=\graph\eta$ and $X_i=\graph\eta|_{Y_i}$ and let $p_Y:(X,X_i)_i\to(Y,Y_i)_i$ and $p_Z:(X,X_i)_i\to(Z,Z_i)_i$ denote the 
projections. 
Then $p_Y$ and $p_Z$ are extensible to definable $C^0$ maps $\overline{p_Y}:\overline X\to\overline Y$ and $\overline{p_Z}:
\overline X\to\overline Z$. 
By the definable triangulation theorem of definable sets we regard $\overline X$ as the underlying polyhedron to a finite simplicial complex $K$ such 
that $X$ and $X_i$ are the unions of some open simplexes in $K$. 
Define $K',\,K''$ and $\tau:(X,X_i)_i\to(X,X_i)_i$ as above. 
Then $(\Im\tau,\tau(X_i))_i$ is a standard family of semi-linear sets, and $\tau^{-1}|_{\Im\tau}:\Im\tau\to X$ is extensible to a 
semi-algebraic $C^0$ map from $\overline{\Im\tau}$ to $\overline X$ by the definition of $\tau$. 
Hence $p_Y\circ(\tau^{-1}|_{\Im\tau}):(\Im\tau,\tau(X_i))_i\to(Y,Y_i)_i$ and $p_Z\circ(\tau^{-1}|_{\Im\tau}):(\Im\tau,\tau(X_i))_i
\to(Z,Z_i)_i$ are definable homeomorphisms between standard semi-linear set families and extensible to definable $C^0$ maps 
$\overline{\Im\tau}\to\overline Y$ and $\overline{\Im\tau}\to\overline Z$ respectively. 
Thus replacing $\eta$ with $p_Y\circ(\tau^{-1}|_{\Im\tau})$ and $p_Z\circ(\tau^{-1}|_{\Im\tau})$ we assume $\eta$ is extensible to 
$\overline\eta$. \par
  Let $K_Y$ be a simplicial decomposition of $\overline Y$ such that $Y$ and $Y_i$ are the unions of some open simplexes in $K_Y$, 
and $K'_Y$ the barycentric subdivision of $K_Y$. 
Set $W_Y=\overline Y-Y$ and define  a definable $C^0$ function $f_Y$ and a definable 
set $U_Y$ as $f_K$ and $U_K$ from $K$. 
Let $K_Z,\,K'_Z,\,W_Z,\,f_Z$ and $U_Z$ be given for $(Z,Z_i)_i$ in the same way. 
Then there exists a definable isotopy $\alpha_{Zt},\ 0\le t\le 1$, of $\overline Z$ preserving $K'_Z$ such that $\alpha_{Z1}\circ
\overline\eta(\sigma)$ is semi-linear for each $\sigma\in K'_Y$. 
Hence subdividing $K'_Z$ we assume $\overline\eta(\sigma)$ is the union of some simplexes in $K'_Z$. \par
  We modify $K'_Y$ and $K'_Z$ to some cell (not simplicial in general) complexes. 
Set
$$
\gather
\tilde Y=(\overline Y-\Dom f_Y)\cup f^{-1}_Y((1/2,\,1]),\ \tilde Y_i=\tilde Y\cap Y_i,\\
\tilde K'_Y=\{\sigma\cap(\tilde Y\overline),\,\big(\sigma\cap f^{-1}(1/2)\overline{\big)}:\sigma\in K'_Y\}\ \text{and}\ \tilde f_Y=2f_Y-1. 
\endgather
$$
Then $\tilde K'_Y$ is a cell complex, $((\tilde Y\overline),\tilde Y,\tilde Y_i)_i$ and $(\overline Y,Y,Y_i)_i$ are semi-linearly 
homeomorphic as shown already, each cell in $\tilde K'_Y$ is of the form $\sigma_1*(\sigma_2\cap f^{-1}_Y(1/2))$ or $\sigma_1*(\sigma_2\cap f^{-1
}_Y([1/2,\,1]))$ for some $\sigma_1,\sigma_2\in K'_Y$ with $\sigma_1\subset\overline Y-\Dom f_Y$ and $\sigma_2\subset\Dom f_Y$, and there exists a definable homotopy $\alpha_{Yt}:(\tilde Y\overline)\to\overline Y,\ 0\le t\le 1$, 
of the homeomorphism such that the following three conditions are satisfied. 
(1) $\alpha_{Yt}$ is a homeomorphism to $\overline Y$ for each $t\in[0,\,1)$ but not for $t=1$, 
(2) $\alpha_{Y1}|_{\tilde Y}$ is a homeomorphism to $Y$, and (3) $\alpha_{Y1}\big(\sigma\cap(\tilde Y\overline)\big)=\sigma$ for $\sigma\in K'_Y$ 
with $\Int\sigma\not\subset W_Y$. 
It follows from (1), (2) and (3) that $\alpha_{Y1}\big(\sigma\cap((\tilde Y\overline)-\tilde Y\overline)\big)=\sigma\cap\overline{W_Y}$ 
for the same $\sigma$ as in (3), and $\sigma\cap((\tilde Y\overline)-\tilde Y\overline)$ and $\sigma\cap\overline {W_Y}$ may be of different dimension. 
Hence, by replacing $(Y,Y_i)_i,\ K'_Y$ and $f_Y$ with $(\tilde Y,\tilde Y_i)_i,\ \tilde K'_Y$ and $\tilde f_Y$, respectively, and keeping the 
notations, we assume $K'_Y$ is a cell complex, for any $t\in[0,\,1)$ there exists a cellular isomorphism from $K'_Y$ to $\{\sigma
\cap(f_Y^{-1}(t)\overline),\sigma\cap(f^{-1}_Y([t,\,1])\overline):\sigma\in K'_Y\}$, the restriction of $f_Y$ to $\{\sigma\in K'_
Y:\sigma\subset\Dom f_Y\}$ is a cellular map to $\{0,1,[0,\,1]\}$, each cell in $K'_Y$ is uniquely described as of the form 
$\sigma_1*\sigma_2$ for $\sigma_1,\sigma_2\in K'_Y$ with $\sigma_1\subset\overline Y-\Dom f_Y$ and $\sigma_2\subset\Dom f_Y$, 
and for such $\sigma_1$ and $\sigma_2$ 
$$
f_Y(ty_1+(1-t)y_2)=f_Y(y_2)\quad\text{for}\ (y_1,y_2,t)\in\sigma_1\times\sigma_2\times[0,\,1). 
$$
Note $(Y,Y\cap\sigma)_{\sigma\in K'_Y}$ is standard. 
In the same way we modify $K'_Z$. \par
  Under these assumptions we will find a definable isotopy $\alpha_t,\ 0\le t\le 1$, of $Y$ preserving $\{\Int\sigma:\sigma\in K'_Y,
\Int\sigma\subset Y\}$ by induction on $\dim Y$ such that $\eta\circ\alpha_1$ is semi-linear, $\alpha_t$ is extensible to a definable homeomorphism $\overline\alpha_t$ of $\overline Y$ for each $t\in[0,\,1)$ but not for $t=1$, and the map $\overline Y\times[0,
\,1)\ni(y,t)\to\overline\alpha_t(y)\in\overline Y$ is continuous. 
Since $W_Y\cap Y=\emptyset$ and $\overline {W_Y}-Y=W_Y$ we have 
$$
\overline{W_Y}-W_Y\subset Y,\ \overline{W_Z}-W_Z\subset Z\ \text{and}\ \eta(\overline{W_Y}-W_Y)=\overline{W_Z}-W_Z. 
$$\par
  First we reduce the problem to the case where $\overline\eta$ is semi-linear on $(\overline{W_Y}-W_Y\overline)$. 
Clearly $(\overline{W_Y}-W_Y,Y_i\cap\overline{W_Y}-W_Y)_i,\,K'_Y|_{(\overline{W_Y}-W_Y\overline)},\,\eta|_{\overline{W_Y}-W_Y},\,
(\overline{W_Z}-W_Z,Z_i\cap\overline{W_Z}-W_Z)_i$ and $K'_Z|_{(\overline{W_Z}-W_Z\overline)}$ satisfy the assumptions on $(Y,Y_i)_i,
\,K'_Y,\,\eta,\,\allowmathbreak(Z,Z_i)_i$ and $K'_Z$. 
Hence by induction hypothesis there exists a definable isotopy $\alpha_{Wt},\ 0\le t\le 1$, of $\overline{W_Y}-W_Y$ preserving 
$\{\Int\sigma:\sigma\in K'_Y,\Int\sigma\subset\overline{W_Y}-W_Y\}$ such that $\eta\circ\alpha_{W1}$ is semi-linear, $\alpha_{Wt}$ is extensible to a definable homeomorphism $\overline{\alpha_{Wt}}$ of $(\overline{W_Y}-W_Y\overline)$ for each $t\in[0,\,1)$, 
and the map $(\overline{W_Y}-W_Y\overline)\times[0,\,1)\ni(y,t)\to\overline{\alpha_{Wt}}(y)\in(\overline{W_Y}-W_Y\overline)$ is 
continuous. 
Then we can extend $\alpha_{Wt}$ to a definable isotopy $\tilde\alpha_{Wt},\ 0\le t\le 1$, of $\overline Y-\big((\overline{W_Y}-W_Y
\overline)-(\overline{W_Y}-W_Y)\big)$ preserving $\{\Int\sigma:\sigma\in K'_Y,\Int\sigma\subset\overline Y-\big((\overline{W_Y}-W_Y
\overline)-(\overline{W_Y}-W_Y)\big)\}$, for which it suffices to see the following statement. 
Given a cell $\sigma$ in $R^n$, a union $\sigma_1$ of open faces of $\sigma$ of codimension $\ge 2$ and a definable isotopy 
$\alpha_{\partial\sigma t},\ 0\le t\le 1$, of $\partial\sigma-\sigma_1$ such that $\alpha_{\partial\sigma 
t}$ is extensible to a definable homeomorphism $\overline{\alpha_{\partial\sigma t}}$ of $\partial\sigma$ for each $t\in[0,\,1)$ and the map $\partial
\sigma\times[0,\,1)\ni(y,t)\to\overline{\alpha_{\partial\sigma t}}(y)\in\partial\sigma$ is continuous, then $\alpha_{\partial
\sigma t}$ is extensible to a definable isotopy $\alpha_{\sigma t},\ 0\le t\le 1$, of $\sigma-\sigma_1$ so that $\alpha_{\sigma t}$ is extensible to a definable homeomorphism $\overline{\alpha_{\sigma t}}$ of $\sigma$ for each $t\in
[0,\,1)$ and the map 
$\sigma\times[0,\,1]-\sigma_1\times\{1\}\ni(y,t)\to\alpha_{\partial\sigma_t(y)}$, or $\overline{\alpha_{\sigma t}}(y),\in\sigma$ 
is continuous. 
We prove as usual by the Alexander trick. 
Assume $0\in\Int\sigma$, and set 
$$
\alpha_{\sigma t}(sy)=
\cases
\alpha_{\partial\sigma 1}(y)\quad&\text{for}\ (y,s,t)\in(\partial\sigma-\sigma_1)\times\{1\}^2,\\
s\overline{\alpha_{\partial\sigma st}}(y)\quad&\text{for}\ (y,s,t)\in\sigma\times([0,\,1]^2-\{1\}^2). 
\endcases
$$
Then $\alpha_{\sigma t}$ fulfills the requirements. 
Hence, replacing $\eta$ with $\eta\circ\tilde\alpha_{W1}|_Y$ we assume $\overline\eta$ is semi-linear on $(\overline{W_Y}-W_Y
\overline)$. \par
  Secondly, we reduce to the case where $\overline\eta(\Dom f_Y)=\Dom f_Z$. 
In the above arguments of uniqueness we start with the barycentric subdivisions of $K_Y$ and $K_Z$ in place of $K_Y$ and $K_Z$ and 
use the same notation. 
Let $K''_Y$ and $K''_Z$ be simplicial subdivisions of the barycentric subdivisions $K'_Y$ and $K'_Z$ of the new $K_Y$ and $K_Z$, respectively, 
such that $\overline\eta|_{(\overline{W_Y}-W_Y\overline)}:K''_Y|_{(\overline{W_Y}-W_Y\overline)}\to K''_Z|_{(\overline{W_Z}-W_Z
\overline)}$ is an isomorphism, and define $f'_Y$ and $f'_Z$ from $K''_Y$ and $K''_Z$ as $f_Y$ and $f_Z$ from $K'_Y$ and $K'_Z$. 
Then $\overline\eta(\Dom f'_Y)=\Dom f'_Z$, $\overline\eta^{-1}(\Dom f'_Z)=\Dom f'_Y$ and by note 1 there exist semi-linear isotopies 
$\gamma_{Yt}$ and $\gamma_{Zt},\ 0\le t\le 1$, of $\overline Y$ and $\overline Z$ preserving $K_Y$ and $K_Z$, respectively, such 
that $\gamma_{Y1}(\Dom f_Y)=\Dom f'_Y$ and $\gamma_{Z1}(\Dom f_Z)=\Dom f'_Z$. 
Define a semi-linear isotopy $\gamma'_{Yt},\ 0\le t\le 1$, of $(\overline{W_Y}-W_Y\overline)$ to be $\overline\eta^{-1}\circ\gamma
_{Zt}\circ\overline\eta$ and that $\gamma''_{Yt},\ 0\le t\le 1$, of $(\overline{W_Y}-W_Y\overline)$ to be $\gamma_{Yt}^{\prime-1}
\circ\gamma_{Yt}$, which are preserving $K_Y|_{(\overline{W_Y}-W_Y\overline)}$, and extend $\gamma''_{Yt},\ 0\le t\le 1$, to a 
semi-linear isotopy $\overline\gamma''_{Yt},\ 0\le t\le 1$, of $\overline Y$ preserving $K_Y$ by the Alexander trick. 
Then $\overline\eta\circ\overline\gamma''_{Y1}(\Dom f_Y)=\Dom f_Z$. 
Hence we assume $\overline\eta(\Dom f_Y)=\Dom f_Z$ from the 
beginning. \par
  Define a cellular map $\phi_Y:K'_Y\to\{0,1,[0,\,1]\}$ by 
$$
\phi_Y=
\cases
0\quad&\text{on}\ K^{\prime 0}_Y\cap(\overline{W_Y}-W_Y)\\
1\quad&\text{on}\ K^{\prime 0}_Y-(\overline{W_Y}-W_Y),
\endcases
$$
which is well-defined because for non-empty $\sigma\in K'_Y$ there exist uniquely $\sigma_1,\sigma_2\in K'_Y$ with $\sigma_1*
\sigma_2=\sigma,\ \sigma_1\subset\overline Y-\Dom f_Y$ and $\sigma_2\subset\Dom f_Y$, hence, $\phi_Y(\sigma)=\{0\}$ if $\sigma_2=
\emptyset,\ \phi_Y(\sigma)=\{1\}$ if $\sigma_1=\emptyset$ and $\phi_Y((1-t)y_1+ty_2)=t$ for $(y_1,y_2,t)\in\sigma_1\times\sigma
_2\times[0,\,1]$. 
Then $\phi_Y^{-1}(0)=\overline Y-\Dom f_Y,\ \phi_Yf_Y$ is extensible to a PL function $\overline{\phi_Y f_Y}$ on $\overline Y$, 
and $\overline{\phi_Yf_Y}^{-1}(0)=\overline{W_Y}$. 
We also define $\phi_Z$ for $K'_Z$ in the same way. 
Note $\phi_Z\circ\overline\eta$ and $\overline{\phi_Zf_Z}\circ\overline\eta$ are definable $C^0$ functions on $\overline Y$ with 
zero sets $\phi_Y^{-1}(0)$ and $\overline{W_Y}$, respectively, and 
$$
\Im(\phi_Y,\overline{\phi_Yf_Y})=\Im(\phi_Z,\overline{\phi_Zf_Z})=\Im(\phi_Z\circ\overline\eta,\overline{\phi_Zf_Z}\circ
\overline\eta)=\{(r,s)\in[0,\,1]^2:s\le r\}. 
$$\par
  Set $F=(F_1,F_2)=(\phi_Z\circ\overline\eta,\overline{\phi_Zf_Z}\circ\overline\eta)$. 
We will, finally, reduce to the case where there exists a definable neighborhood $N$ of $(0,\,\epsilon]\times\{0\}$ in $(0,\,\epsilon]\times
[0,\,\epsilon]$ for some $\epsilon\in(0,\,1)$ such that $F|_{F^{-1}(\overline N)}$ is extensible to a PL map $\tilde F=(\tilde F_1,
\tilde F_2)$ on $\overline Y$ with $\tilde F_1\ge 0,\ \tilde F_2\ge 0,\ \tilde F_1^{-1}(0)=F_1^{-1}(0)$ and $\tilde F_2^{-1}(0)=
F_2^{-1}(0)$. 
Before carrying out reduction we prove the uniqueness under this condition. \par
  Compare $\tilde F$ with the PL map $G=(G_1,G_2)=(\phi_Y,\overline{\phi_Yf_Y})$ around $\phi^{-1}_Y(0)$. 
Then we can assume $\tilde F=G$ there for the following reason. 
There exists a semi-linear isotopy $\delta_t,\ 0\le t\le 1$, of $\overline Y$ preserving $K'_Y$ such that $\tilde F_1\circ\delta_1=
\phi_Y$ on a neighborhood of $\phi^{-1}_Y(0)$ in $\overline Y$ because $\tilde F_1$ and $\phi_Y$ are non-negative PL functions on 
$\overline Y$ with the same zero set. 
(Existence of $\delta_1$ under these conditions is stated in Lemma I.3.10 in the real case, and existence of $\delta_t$ in the 
general case is clear by its proof and remark 1.3.) 
Hence we assume $\tilde F_1=\phi_Y$ from the beginning. 
Let $K''_Y$ be a simplicial subdivision of $K'_Y$, where $\tilde F$ and $G$ are simplicial, $0<\epsilon\in R$ so small that $\phi_Y
(K^{\prime\prime 0}_Y)\cap(0,\,\epsilon]=\emptyset$, and $K'''_Y$ a simplicial subdivision of the cell complex $\{\sigma''\cap
\phi^{-1}_Y(\sigma_\epsilon):\sigma''\in K''_Y,\sigma_\epsilon\in\{0,\epsilon,1,[0,\,\epsilon],[\epsilon,\,1]\}\}$ without new 
vertices. 
Set $K'''_{Y\epsilon}=K'''_Y|_{\phi^{-1}_Y(\epsilon)}$ and compare $\tilde F_2$ and $G_2$ on $\phi^{-1}_Y(\epsilon)$. 
Since both are PL and non-negative and have the same zero set there exists a semi-linear isotopy $\delta'_{\epsilon t},\ 0\le 
t\le 1$, of $\phi^{-1}_Y(\epsilon)$ preserving $K'''_{Y\epsilon}$ such that $\tilde F_2\circ\delta'_{\epsilon 1}=G_2$ on a 
neighborhood of $\phi^{-1}_Y(\epsilon)\cap G_2^{-1}(0)$ in $\phi^{-1}_Y(\epsilon)$, say, $\tilde F^{-1}(N')$ for a semi-linear 
neighborhood $N'$ of the point $(\epsilon,0)$ in $\{\epsilon\}\times[0,\,\epsilon]$. 
We need to extend $\delta'_{\epsilon t}$ to a semi-linear isotopy $\delta'_t,\ 0\le t\le 1$, of $\overline Y$ preserving $K''_Y$. 
That is possible on $\phi^{-1}_Y([\epsilon,\,1])$ because the restriction of $\phi_Y$ to $\{\sigma\cap\phi_Y^{-1}(\epsilon),\sigma
\cap\phi_Y^{-1}(\epsilon'),\sigma\cap\phi_Y^{-1}([\epsilon,\,\epsilon']):\sigma\in K''_Y\}$ is a cellular map to $\{\{\epsilon\},
\{\epsilon'\},[\epsilon,\,\epsilon']\}$ and trivial for some $\epsilon'\,(>\epsilon)\in R$. 
We define $\delta'_t$ on $\phi^{-1}_Y([0,\,\epsilon])$ by 
$$
\gather
\delta'_t(ry_1+(1-r)y_2)=ry_1+(1-r)\delta'_{\epsilon t}(y_2)\\
\text{for}\ (y_1,y_2,r)\in\sigma_1\times\sigma_2\times[0,\,1],\ \sigma_1\in K'''_Y|_{\phi^{-1}_Y(0)},\ \sigma_2\in K'''_{Y\epsilon}
\ \text{with}\ \sigma_1*\sigma_2\in K'''_Y. 
\endgather
$$
Then $\tilde F_2\circ\delta'_1=G$ on $\tilde F^{-1}((0,0)*N')$, where the cone $(0,0)*N'$ is a semi-linear neighborhood of $(0,\,
\epsilon]\times\{0\}$ in $(0,\,\epsilon]\times[0,\,\epsilon]$. 
Hence we can assume $\tilde F=G$ on a neighborhood of $\phi_Y^{-1}(0)$ in $\overline Y$ and $F=G$ on $F^{-1}(\overline N)$. \par
  Moreover we can assume $f_Z\circ\overline\eta=f_Y$ on $F^{-1}(N'')$ for some neighborhood $N''$ of $[\epsilon,\,1]\times\{0\}$ 
in $[\epsilon,\,1]\times[0,\,1]$ for the following reason. 
Note $f_Z\circ\overline\eta|_{F^{-1}([\epsilon,1]\times[0,1])}$ and $f_Y|_{F^{-1}([\epsilon,1]\times[0,1])}$ are definable $C^0$ 
and PL, respectively, functions with zero set $F^{-1}([\epsilon,1]\times\{0\})$, and $f_Z\circ\overline\eta$ and $f_Y$ are PL 
and coincide each other on $F^{-1}\big([\epsilon/2,\,\epsilon]\times[0,\,\epsilon'']\big)$ for some $\epsilon''>0\in R$. 
Choose $\epsilon''$ so small that $K_Y^{\prime\prime\prime 0}\cap F^{-1}\big([\epsilon/2,\,\epsilon]\times[0,\,\epsilon'']\big)\subset F^{-1}(
\epsilon,0)$, where $K'''_Y$ was defined above. 
By theorem 2.2,(2) there exists a definable homeomorphism $\pi_Y$ of $F^{-1}\big([\epsilon,\,1]\times[0,\,1]\big)$ preserving $\{\sigma\cap 
F^{-1}(\{\epsilon\}\times[0,\,\epsilon'']),\sigma\cap F^{-1}(\{\epsilon\}\times[\epsilon'',\,\epsilon]),\sigma\cap F^{-1}\big([
\epsilon,\,1]\times[0,\,1]\big):\sigma\in K'''_Y\}$ such that $f_Z\circ\overline\eta\circ\pi_Y$ is PL. 
Apply the complement of theorem 2.2 to $f_Z\circ\overline\eta$ and $f_Z\circ\overline\eta\circ\pi_Y$ on $F^{-1}(\{\epsilon\}
\times[0,\,\epsilon''])$. 
Then there is a PL homeomorphism $\pi'_Y$ of $F^{-1}(\{\epsilon\}\times[0,\,\epsilon''])$ preserving $\{\sigma\cap F^{-1}(\{
\epsilon\}\times[0,\,\epsilon'']):\sigma\in K'''_Y\}$ such that $f_Z\circ\overline\eta=f_Z\circ\overline\eta\circ\pi_Y\circ\pi'
_Y$ on $F^{-1}(\{\epsilon\}\times[0,\,\epsilon''])$. 
Extend $\pi'_Y$ to a PL homeomorphism of $F^{-1}\big([\epsilon,\,1]\times[0,\,1]\big)$ preserving $\{\sigma\cap F^{-1}(\{\epsilon\}\times[
0,\,\epsilon'']),\sigma\cap F^{-1}(\{\epsilon\}\times[\epsilon'',\epsilon]),\sigma\cap F^{-1}\big([\epsilon,\,1]\times[0,\,1]\big):\sigma
\in K'''_Y\}$ by the Alexander trick, and replace $\pi_Y$ with $\pi_Y\circ\pi'_Y$. 
Then we can assume $f_Z\circ\overline\eta\circ\pi_Y=f_Z\circ\overline\eta=f_Y$ on $F^{-1}(\{\epsilon\}\times[0,\,\epsilon''])$ 
from the beginning. 
Next extend $\pi_Y$ to a definable homeomorphism of $\overline Y$ preserving $K'''_Y$ by the Alexander trick so that $f_Z\circ
\overline\eta\circ\pi_Y=f_Z\circ\overline\eta$ on $F^{-1}\big([\epsilon/2,\,\epsilon]\times[0,\,\epsilon'']\big)$ and $\pi_Y=\id$ on 
$F^{-1}\big([0,\,\epsilon/2]\times[0,\,1]\big)$. 
Furthermore by the method of construction of $\pi_Y$ and $\pi'_Y$ there exists a definable isotopy of $\overline Y$ to $\pi_Y$ 
preserving $K'''_Y$. 
Thus forgetting $\pi_Y$ we assume $f_Z\circ\overline\eta$ is PL on $F^{-1}\big([\epsilon,\,1]\times[0,\,1]\cup[\epsilon/2,\,\epsilon]
\times[0,\,\epsilon'']\big)$. 
Then we have a PL homotopy $\pi'_{Yt},\ 0\le t\le 1$, of $\overline Y$ preserving $K'''_Y$ such that $f_Z\circ\overline\eta\circ
\pi'_{Y1}=f_Y$ on $F^{-1}(N'')$ for some neighborhood $N''$ of $[\epsilon/2,\,1]\times\{0\}$ in $[\epsilon/2,\,1]\times[0,\,1]$ 
(Lemma I.3.10). 
Moreover the proof of Lemma I.3.10 says that $\pi'_{Yt}$ is chosen to be the identity map on $F^{-1}\big([0,\,\epsilon/2]\times[0,
\,1]\big)$. 
Consequently $F\circ\pi'_{Yt}=G$ on $F^{-1}\big(N\cap[0,\,\epsilon/2]\times[0,\,1]\big)$ and we can assume $f_Z\circ\overline\eta=f_Y$ 
on $F^{-1}(N'')$. \par
  Let $K''_Y$ and $L$ be a simplicial subdivision of $K'_Y$ and a simplicial decomposition of $\{(r,s)\in[0,\,1]^2:s\le r\}$, 
respectively, such that $G:K''_Y\to L$ is simplicial. 
Let $\psi_1$ and $\psi_2$ be definable non-negative $C^0$ functions on $[0,\,1]$ such that $\psi_1(r)<r$ and $\psi_2(r)<r$ for $r
\in(0,\,1]$, $\psi^{-1}_i(0)=\{0\}$ and 
$$
\gather
\cup\{\Int|\st(\sigma,L)|\!:\!\sigma\in L,\Int\sigma\!\subset\!(0,\,1]\times\{0\}\}\!\supset\!\{(r,s)\in(0,\,1]^2\!:\!s\le\psi_i(r)\},\,i=1,
2. \\
\Psi_i=\{(r,s)\in(0,\,1]^2:\psi_i(r)\le s\le r\},\ i=1,2. \tag"Set"
\endgather
$$
{\it Note 2.} Then there exists a definable isotopy $\delta''_t,\ 0\le t\le 1$, of $\overline Y$ preserving $K'_Y$ such that 
$\delta''_1(\overline{G^{-1}(\Psi_1)})=\overline{G^{-1}(\Psi_2)}$. \par
  The reason is the following. 
Let $\psi_0$ and $\psi_3$ be functions with the same properties as $\psi_1$ and $\psi_2$ and such that $\psi_0<\psi_i<
\psi_3,\ i=1,2$, on $(0,\,1]$. 
First define a definable isotopy $\delta''_{Lt},\ 0\le t\le 1$, of $|L|$ preserving $L$ so that for each $(r,s)\in|L|$, 
$$
\delta''_{Lt}(r,s)=
\cases
(r,s)\quad&\text{for}\ s\in[0,\,\psi_0(r)]\cup[\psi_3(r),\,r]\\
\big(r,t\psi_2(r)+(1-t)\psi_1(r)\big)\quad&\text{for}\ s=\psi_1(r),
\endcases
$$
and $\delta''_{Lt}(r,\cdot)$ is linear on the segments $[\psi_0(r),\,\psi_1(r)]$ and $[\psi_1(r),\,\psi_3(r)]$. 
Next lift $\delta''_{Lt}$ to a definable isotopy of $\overline Y$ by the following note. 
Then the isotopy is the required one. \par
  {\it Note 3.} 
Let $h:K\to L$ be a cellular map between finite cell complexes and $L_1$ a subcomplex of $L$. 
Let $\delta_L$ be a definable homeomorphism of $|L|$ preserving $L$ such that $\delta_L=\id$ on $|L_1|$. 
Then there exists a definable homeomorphism $\delta_K$ of $|K|$ preserving $K$ such that $\delta_K=\id$ on $h^{-1}(|L_1|)$ and 
$\delta_L\circ h=h\circ\delta_K$. \par
  {\it Proof of note 3.} 
By induction on $l\in\N$ we assume $\delta_K$ is already constructed on $|K^l|$ and let $\sigma\in K^{l+1}-K^l$. 
Then it suffices to extend $\delta_K|_{\partial\sigma}$ to $\sigma$ so that $\delta_L\circ h=h\circ\delta_K$ on $\sigma$. 
Since the extension is trivial if $\dim\sigma=\dim h(\sigma)$ or $h(\sigma)\in L_1$, we suppose $\dim\sigma>\dim h(\sigma)$ and 
$h(\sigma)\not\in L_1$. 
Let $c:h(\sigma)\to\sigma$ be a definable $C^0$ cross-section of $h|_\sigma$ such that $c(x)\in\Int (h|_{\sigma})^{-1}(x)$ for 
$x\in h(\sigma)$. 
Set 
$$
\delta'_K\big(ty+(1-t)c\circ h(y)\big)=t\delta_K(y)+(1-t)c\circ\delta_L\circ h(y)\quad\text{for}\ (y,t)\in\partial\sigma\times[0,\,1]. 
$$
Then $\delta'_K$ is a definable homeomorphism of $\sigma$ preserving $K|_\sigma$, $h\circ\delta'_K=\delta_L\circ h$ and $\delta'
_K=\delta_K$ on $\partial(h|_\sigma)^{-1}(x)$ for each $x\in h(\sigma)$. 
The last equality, however, does not necessarily hold on $\cup\{\Int(h|_\sigma)^{-1}(x):x\in\partial h(\sigma),\dim(h|_\sigma)^{-1}(
x)>0\}\,(\subset\partial\sigma)$. 
We need to modify $\delta'_K$. 
Consider $\delta_K^{\prime -1}\circ\delta_K$ on $\partial\sigma$ and the identity map on $h(\sigma)$ in place of $\delta_K$ on 
$\partial\sigma$ and $\delta_L$ on $h(\sigma)$. 
Then we only need to find a definable homeomorphism $\delta''_K$ of $\sigma$ preserving $K|_\sigma$ such that $\delta''_K=\delta_K
^{\prime -1}\circ\delta_K$ on $\partial\sigma$ and $h\circ\delta''_K=h$ because if such $\delta''_K$ exists then $\delta'_K\circ\delta''
_K$ is what we want. 
Let $a\in\Int\sigma$ and set 
$$
\delta''_K\big(ty+(1-t)a\big)=t\delta_K^{\prime -1}\circ\delta_K(y)+(1-t)a\quad\text{for}\ (y,t)\in\partial\sigma\times[0,\,1]. 
$$
Then $\delta''_K$ fulfills the requirements. 
Thus notes 2 and 3 are proved. \par
  Choose $\psi_1$ and $\psi_2$ in note 2 so that $\psi_1(r)=d_1r$ for some small $d_1>0\in R$, $\psi_2=\psi_1$ on $[\epsilon,\,1]$ 
and $\Psi_2\cup N\cup N''\subset\{(r,s)\in(0,\epsilon]\times[0,\,1]:s\le r\}$, and set $G_Z=(\phi_Z,\overline{\phi_Z f_Z})$. 
Then $F=G_Z\circ\overline\eta$, $G^{-1}(\Psi_1)=f_Y([d_1,\,1])$ by definition of $G$, $\overline{G^{-1}(\Psi_1)}$ is identified with 
$\overline{G^{-1}(\Psi_2)}$ through a definable homeomorphism of $\overline Y$ by note 2, and $\overline\eta|_{\overline{G^{-1}(\Psi_2)}
}$ is a homeomorphism onto $\overline{G^{-1}_Z(\Psi_2)}$ for the following reason. 
$$
\gather
\overline{G^{-1}(\Psi_2)}\cap(\overline Y-G^{-1}(\Psi_2)\overline)=\{y\in\overline Y:\overline{\phi_Yf_Y}(y)=\psi_2\circ\phi_Y(y)\}
\subset\\
\qquad\qquad\qquad\qquad\qquad\qquad\qquad\qquad\qquad\qquad F^{-1}(0,0)\cup F^{-1}(N\cup N''),\\
\{y\in F^{-1}(N):\overline{\phi_Yf_Y}(y)=\psi_2\circ\phi_Y(y)\}\qquad\qquad\qquad\qquad\qquad\qquad\qquad\qquad\\
\qquad=\{y\in F^{-1}(N):\overline{\phi_Zf_Z}\circ\overline\eta(y)=\psi_2\circ\phi_Z\circ\overline\eta(y)\}
\quad\text{by}\ F=G\ \text{on}\ F^{-1}(N),\\
\{y\in F^{-1}(N''):\overline{\phi_Yf_Y}(y)=\psi_2\circ\phi_Y(y)\}\qquad\qquad\qquad\qquad\qquad\qquad\qquad\qquad\\
=\{y\in F^{-1}(N''):f_Y(y)=d_1\}\quad\text{by}\ \psi_2(r)=d_1r\ 
\text{for}\ r\in[\epsilon,\,1]\\
=\{y\in F^{-1}(N''):f_Z\circ\overline\eta(y)=d_1\}\quad\text{by}\ f_Z\circ\overline\eta=f_Y\text{on}\ F^{-1}(N'')\\
=\{y\in F^{-1}(N''):\overline{\phi_Zf_Z}\circ\overline\eta(y)=\psi_2\circ\phi_Z\circ
\overline\eta(y)\}\quad\text{by}\ \psi_2(r)=d_1r\ \text{for}\ r\in[\epsilon,\,1],\\
\{y\in\overline Y:\overline{\phi_Yf_Y}(y)=\psi_2\circ\phi_Y(y)\}=\{y\in\overline Y:\overline{\phi_Zf_Z}\circ\overline\eta(y)=\psi_
2\circ\phi_Z\circ\overline\eta(y)\}. \tag"and, hence,"
\endgather
$$\par
  Define $\psi_{it}$ and $\Psi_{it}$ to be $t\psi_i$ and $\{(r,s)\in(0,\,1]^2:\psi_{it}(r)\le s\le r\}$, respectively, for each $t\in(0,\,1]$, and repeat the above arguments for $\psi_{it}$ and $\Psi_{it}$. 
Then by the above proof we have a definable isotopy $\omega_{Yt},\ 0\le t\le 1$, of $\overline Y$ preserving $K''_Y$ such that 
$\omega_{Yt}(\overline{G^{-1}(\Psi_{1t})})=\overline{G^{-1}(\Psi_{2t})}$ for $t\in(0,\,1]$ and $\overline\eta|_{\overline G^{-1}
(\Psi_{2t})}$ is a homeomorphism onto $\overline{G^{-1}_Z(\Psi_{2t})}$. 
Consider $Z$. 
Then $G^{-1}_Z(\Psi_{1t})=f^{-1}_Z([td_1,\,1])$ for each $t\in(0,\,1]$ and there exists a definable isotopy $\omega_{Zt},\ 0\le t
\le 1$, of $\overline Z$ preserving $K''_Z$ such that $\omega_{Zt}(\overline{G^{-1}_Z(\Psi_{1t})})=\overline{G^{-1}_Z(\Psi_{2t})}$ 
for $t\in(0,\,1]$. 
Remember that there are definable homeomorphisms $\kappa_{Yt}:\overline Y\to\overline{f_Y([td_1,\,1])}$ and $\kappa_{Zt}:
\overline Z\to\overline{f_Z([td_1,\,1])}$ for each $t\in[0,\,1]$ preserving $\{\Int\sigma:\sigma\in K'_Y,\Int\sigma\subset Y\}$ 
and $\{\Int\sigma:\sigma\in K'_Y,\Int\sigma\subset Y\}$, respectively, such that $\kappa_{Y0}=\kappa_{Z0}=\id$ and the maps 
$\overline Y\times[0,\,1]\ni(y,t)\to\kappa_{Yt}(y)\in\overline Y$ and $\overline Z\times[0,\,1]\ni(z,t)\to\kappa_{Zt}(z)\in
\overline Z$ are continuous. \par
  In conclusion, set $\kappa'_{Yt}=\omega_{Yt}\circ\kappa_{Yt}$ and $\kappa'_{Zt}=\omega_{Zt}\circ\kappa_{Zt}$. 
Then $\kappa'_{Yt}:\overline Y\to\overline{G^{-1}(\Psi_{2t})}$ and $\kappa'_{Zt}:\overline Z\to\overline{G^{-1}_Z(\Psi_{2t})},\ 0
\le t\le 1$, are definable homeomorphisms preserving $\{\Int\sigma:\sigma\in K'_Y,\Int\sigma\subset Y\}$ and $\{\Int\sigma:\sigma
\in K'_Z,\Int\sigma\subset Z\}$, respectively, such that $\kappa'_{Y0}=\kappa'_{Z0}=\id$ and the maps $\overline Y\times[0,\,1]
\ni(y,t)\to \kappa'_{Yt}(y)\in\overline Y$ and $\overline Z\times[0,\,1]\ni(z,t)\to \kappa'_{Zt}(z)\in\overline Z$ are continuous. 
Note $\Im\overline\eta\circ\kappa'_{Yt}=\Im\kappa'_{Zt}$ for each $t\in[0,\,1]$. 
Here we can assume $\kappa_{Z1}^{\prime -1}\circ\eta\circ\kappa'_{Y1}$ is semi-linear by theorem 2.1 through some definable 
isotopy of $\overline Y$ because $\kappa^{\prime -1}_{Z1}\circ\eta\circ\kappa'_{Y1}:\overline Y\to\overline Z$ is a definable 
homeomorphism. 
Now we define $\alpha_t$ by $\kappa'_{Yt}$ and $\kappa'_{Zt}$ as follows
$$
\alpha_t=
\cases
\eta^{-1}\circ\kappa_{Z1}^{\prime -1}\circ\eta\circ\kappa'_{Y1}&\text{for}\ t=1\\
\kappa_{Y1-t}^{\prime -1}\circ\eta^{-1}\circ\kappa'_{Z1-t}\circ\kappa_{Z1}^{\prime -1}\circ\eta\circ\kappa'_{Y1}&\text
{otherwise}.
\endcases
$$
Then $\alpha_t$ fulfills the requirements that $\alpha_t,\ 0\le t\le 1$, is a definable isotopy of $Y$ preserving $\{\Int\sigma:
\sigma\in K'_Y,\Int\sigma\subset Y\}$, $\eta\circ\alpha_1$ is semi-linear, $\alpha_t$ is extensible to a 
definable homeomorphism $\overline\alpha_t$ of $\overline Y$ for each $t\in[0,\,1)$, and the map $\overline Y\times[0,\,1)\ni(y,t)\to\overline\alpha_t
(y)\in\overline Y$ is continuous. 
Thus the uniqueness is proved under the assumption that $F$ is triangulable on $F^{-1}(\overline N)$. \par
  It remains to triangulate $F$ on $F^{-1}(\overline N)$. 
To be precise, what we prove is the following statement. \par
  Let $K$ be a finite simplicial complex in $R^n$ with underlying polyhedron $X$ and $H=(H_1,H_2):X\to R^2$ a definable $C^0$ map 
with $H_1\ge 0$ and $H_2\ge 0$. 
Then there exist a definable isotopy $\zeta_t,\ 0\le t\le 1$, of $X$ preserving $K$ and a definable neighborhood $N$ of $(0,\,
\epsilon]\times\{0\}$ in $(0,\,\epsilon]\times[0,\,\epsilon]$ for some $\epsilon>0\in R$ such that $H\circ\zeta_1|_{(H\circ\zeta_
1)^{-1}(\overline N)}$ is extensible to a PL map $X\to R^2$. \par
  Note a triangulation of $H$ is impossible in general, e.g., the blowing-up $[0,\,1]^2\ni(x_1,x_2)\to(x_1,x_1x_2)\in R^2$. \par
  {\it Proof of the statement.} 
We argue as in the proof of theorem 2.2. 
Let $p:R^n\times R^2\to R^2$ denote the projection. 
Set $\{X_i\}=\{\Int\sigma:\sigma\in K\}$, $A=\graph H$ and $B_{(s,t)}=\{x\in R^n:(x,s,t)\in B\}$ for $(s,t)\in R^2$ and any subset 
$B$ of $R^n\times R^2$. 
Let $\{A_j\}_j$ be a finite stratification of $A$ into definable connected $C^1$ manifolds compatible with $\{X_i\times 
R^2\}_i$, with the frontier condition and such that $p|_A:\{A_j\}_j\to\{ p(A_j)\}_j$ is a definable $C^1$ stratification 
of $p|_A$. 
Let $A'$ denote the union of $A_j$ with $\dim A_{j(s,t)}<n$ for each $(s,t)\in R^2$. 
Set $S^{n-1}=\{\lambda\in R^n:|\lambda|=1\}$, let $T_{(s,t)}\subset S^{n-1}$ denote the closure of the set of singular directions 
$\lambda$ for $A'_{(s,t)}$ and set $T=\{(\lambda,s,t)\in S^{n-1}\times R^2:\lambda\in T_{(s,t)}\}$ and $\tilde T=(\overline T-T
\overline)$. 
Then $T_{(s,t)}$ is a definable set of dimension $<n-1$ for each $(s,t)$, and $T,\overline T,\overline T-T,\tilde T$ are definable 
sets of dimension smaller than $n+1,n+1,n,n$ respectively. \par
  We will find a definable subset of $S^{n-1}$ of dimension $<n-1$ any definable neighborhood of which includes $(\cup_{(s,t)\in
\overline N}T_{(s,t)}\overline)$ for some $\epsilon$ and $N$ in the statement. 
Consider a definable $C^1$ stratification of $p|_{\overline T}$ such that the stratification of $\overline T$ is compatible with 
$T$, and choose $\epsilon$ and $N$ so small that $\overline N-[0,\,\epsilon]\times\{0\}$ is included in a stratum of the 
stratification of $p(\overline T)$. 
Then 
$$
(\cup_{(s,t)\in\overline N}T_{(s,t)}\overline)\subset(\cup_{(s,t)\in\overline N-\{(0,0)\}}\tilde T_{(s,t)}\overline)\cup\cup_{(s,t)
\in\overline N}T_{(s,t)}. 
$$\par
  First we consider $(\cup_{(s,t)\in\overline N-\{(0,0)\}}\tilde T_{(s,t)}\overline)$. 
Let $p|_{\tilde T}:\{\tilde T_j\}_j\to\{p(\tilde T_j)\}_j$ be a definable $C^1$ stratification of $p|_{\tilde T}:\tilde T\to p(
\tilde T)$. 
It suffices to treat only $T_j$ such that $p(T_j)$ is of dimension 1 and includes $(0,\,\epsilon]\times\{0\}$ for some $\epsilon>0$ 
or $p(T_j)\cup(0,\,\epsilon]\times\{0\}$ is a neighborhood of $(0,\,\epsilon]\times\{0\}$ in $(0,\,\epsilon]\times[0,\,\epsilon]$. 
In the former case, $\tilde T_{j(s,t)}$ is of dimension $<n-1$ for each $(s,t)\in(0,\,\epsilon]\times\{0\}$ and, hence, $(\tilde T_j
\overline)-\tilde T_j$ is of dimension $<n-1$. 
Let $O_j$ be any definable neighborhood of $q\big(((\tilde T_j\overline)-\tilde T_j\overline)\big)$ in $S^{n-1}$, where $q:R^n\times R^2
\to R^n$ denotes the projection. 
Then there exists $\epsilon>0$ such that $(\cup_{(s,t)\in(0,\,\epsilon]}\tilde T_{j(s,t)}\overline)\subset O_j$. 
In the latter case, $\tilde T_{j(s,t)}$ is of dimension $<n-2$, $\cup_{s^2+t^2=c}\tilde T_{j(s,t)}$ is of dimension $<n-1$ for 
each $c>0\in R$, and $(\tilde T_j\overline)-\tilde T_j$ is of dimension $<n-1$. 
Hence for any definable neighborhood $O_j$ of $q\big(((\tilde T_j\overline)-\tilde T_j\overline)\big)$ in $S^{n-1}$ there exist $\epsilon$ 
and $N$ such that $(\cup_{(s,t)\in\overline N}T_{j(s,t)}\overline)\subset O_j$. 
Let $J$ denote the family of such $j'$s. 
Then, for any definable neighborhood $O$ of $q\big(\cup_{j\in J}((\tilde T\overline)-\tilde T\overline)\big)$ in $S^{n-1}$ there exist 
$\epsilon$ and $N$ such that $(\cup_{(s,t)\in\overline N-\{(0,0)\}}\tilde T_{(s,t)}\overline)\subset O$. \par
  Next we consider $\cup_{(s,t)\in\overline N}T_{(s,t)}$. 
Set 
$$
\gather
S(s)=\overline{T\cap R^n\times\{s\}\times R}\cap R^n\times\{(s,0)\}\quad\text{for each}\ s\ge 0\in R,\\
S=\overline{\cup_{s\ge 0\in R}S(s)}\cap R^n\times\{(0,0)\}. 
\endgather
$$
Then $S(s)$ and $S$ are definable sets of dimension $<n-1$. 
Let $O'$ be a definable neighborhood of $q(S)$ in $S^{n-1}$. 
Then $q\big((\cup_{s\in[0,\,\epsilon]}S(s)\overline)\big)\subset O'$ for some $\epsilon>0\in R$, and there exists a definable non-negative 
$C^0$ function $\delta$ on $[0,\,\epsilon]$ such that $\delta^{-1}(0)=\{0\}$ and $T_{(s,t)}\subset O'$ for $s\in[0,\,\epsilon]$ 
and $0\le t\le\delta(\epsilon)$. 
Hence $\cup_{(s,t)\in\overline N}T_{(s,t)}\subset O'$ for $N=\{(s,t)\in(0,\,\epsilon]\times R:0\le t\le\delta(\epsilon)\}$. \par
  Thus we obtain the required definable subset of $S^{n-1}$ of dimension $<n-1$. 
Choose $k$ and $N$ so that $\big(\cup_{(s,t)\in\overline N}T_{(s,t)}\overline{\big)}$ is included in a small definable neighborhood of the set. 
Then $\cup_{(s,t)\in\overline N}T_{(s,t)}\not= S^{n-1}$. 
Hence there exists a non-singular direction for $A_{(s,t)}$ for any $(s,t)\in\overline N$. \par
  After changing linearly the coordinate system of $R^n$ we can assume $(1,0,...,0)\in R^n$ is not a singular direction for $A'_
{(s,t)}$ for any $(s,t)\in\overline N$. 
Repeat the same arguments for the image of $\cup_{(s,t)\in\overline N}A'_{(s,t)}$ under the projection $R^n\times R^2\to R^{n-1}
\times R^2$ forgetting the first factor of $R^n$ and so on. 
Then in the same way as in the proof of theorem 2.2 we obtain $N$ and a definable isotopy $\zeta_t,\ 0\le t\le 1$, of $X$ such 
that $H\circ\zeta_1|_{(H\circ\zeta_1)^{-1}(\overline N)}$ is PL. 
Here we can choose $\zeta_t$ preserving $K$. 
Thus the statement is proved. 
\qed
\enddemo
  Standard semi-linearization of a bounded definable $C^0$ function is possible as follows. 
\proclaim{Theorem 2.6}
(1) Let $(X,X_i)_i$ be a finite family of definable sets and $f:X\to R$ a bounded definable $C^0$ function. 
Then there exist a standard family of semi-linear sets $(Y,Y_i)_i$ and a definable homeomorphism $\pi:(Y,Y_i)_i\to(X,X_i)_i$ such 
that $f\circ\pi$ is semi-linear. \par
  (2) If $X$ is a standard semi-linear set and $P$ is a cell decomposition of $\overline X$ such that $X$ is the union of some open cells in $P$ and $(X,X\cap\sigma)_{\sigma\in P}$ is standard then we can choose $Y$ and $\pi$ in (1) so that $Y=X$ and there exists a definable isotopy $\pi_t,\ 0\le t\le 1$, of $X$ from id to $\pi$ preserving $P$.  
\endproclaim
\demo{Proof}
{\it (1).} Let $X$ be contained and bounded in $R^n$. 
Replacing $X$ with the graph of $f$ we assume $f$ is extensible to a definable $C^0$ function $\overline f$ on $\overline X$. 
By theorem 2.5 we have a standard semi-linearization $\pi:(Y,Y_i)_i\to(X,X_i)_i$ of $(X,X_i)$. 
In that proof $\pi$ is constructed so as to be extensible to a definable $C^0$ map $\overline\pi:\overline Y\to\overline X$. 
Let $K$ be a simplicial decomposition of $\overline Y$ such that $Y$ and $Y_i$ are the unions of some open simplexes in $K$. 
Apply theorem 2.2,(2) to $\overline f\circ\overline\pi$ and $K$. 
Then there exists a definable homeomorphism $\tau$ of $\overline Y$ preserving $K$ and, hence, $Y_i$ such that $\overline f\circ\overline\pi\circ\tau$ is PL. 
Hence $\pi\circ\tau|_Y:(Y,Y_i)_i\to(X,X_i)_i$ is the required semi-linearization of $f$. \par
  {\it (2).} By (1) we have a standard family of semi-linear sets $(Y,Y_\sigma,Y_{i,\sigma})_{i,\sigma\in P}$ and a definable homeomorphism $\tau:(Y,Y_\sigma,Y_{i,\sigma})_{i,\sigma\in P}\to(X,X\cap\sigma,X_i\cap\sigma)_{i,\sigma\in P}$ such that $f\circ\tau$ is semi-linear. 
Apply theorem 2.5 to a definable homeomorphism $\tau:(Y,Y_\sigma)_{\sigma\in P}\to(X,X\cap\sigma)_{\sigma\in P}$ between standard semi-linear set families. 
Then there exists a definable isotopy $\tau_t:(Y,Y_\sigma)_\sigma\to(X,X\cap\sigma)_\sigma,\ 0\le t\le 1$, from $\tau$ to a semi-linear homeomorphism through homeomorphisms. 
Set $\pi_t=\tau\circ\tau_t^{-1}$. 
Then $\pi_t,\ 0\le t\le 1$, is a definable isotopy of $(X,X\cap\sigma)_\sigma$ preserving $P$; $f\circ\pi_1$ is semi-linear because $f\circ\pi_1=f\circ\tau\circ\tau_1^{-1}$; each $\pi^{-1}_1(X_i)$ is semi-linear because $\pi^{-1}_1(X_i)=\tau_1(\tau^{-1}(X_i))=\tau_1(\cup_{\sigma\in P}Y_{i,\sigma})$; $(X,\pi_1^{-1}(X_i))$ is standard because $(Y,\cup_{\sigma\in P}Y_{i,\sigma})_i$ is standard and $\tau_1:(Y,\cup_{\sigma\in P}Y_{i,\sigma})_i\to(X,\pi^{-1}_1(X_i))_i$ is a semi-linear homeomorphism. 
Hence $\pi_1:(X,\pi_1^{-1}(X_i))\to(X,X_i)$ is the required standard semi-linearization of $f$. \qed
\enddemo
  Uniqueness of the semi-linearization in (1) does not necessarily hold. 
A counter-example can be constructed as the above-mentioned counter-example to uniqueness of a definable semi-linearization of a family of definable sets. \medskip
  It seems very possible that any definable fiber bundle with compact polyhedral base and total space is definably isomorphic 
to some PL fiber bundle. 
We give a partial answer as follows. 
We define a {\it definable microbundle}, a {\it definable isomorphism} of two definable microbundles and a {\it definable 
s-isomorphism} as in the PL case over $\R$ in [M]. 
See [M] for properties of PL microbundles. 
\proclaim{Remark 2.7}
A definable microbundle $\xi:B\overset i\to\to E\overset j\to\to B$ is definably $s$-isomorphic to some 
PL microbundle if $B$ is a compact polyhedron. 
\endproclaim
\demo{Proof of remark 2.7}
Let $B\subset R^n$, $B'$ a regular neighborhood of $B$ in $R^n$, and $p:B'\to B$ a PL retraction. 
Let $E'\,(\subset B'\times E)$ denote the fiber product of $p:B'\to B$ and $j:E\to B$, and define $i':B'\to E'$ to be $(\id,i
\circ p)$ and $j':E'\to B'$ to be the projection. 
Then $\xi':B'\overset i'\to\longrightarrow E'\overset j'\to\longrightarrow B'$ is a definable microbundle and we regard $\xi$ 
as the restriction of $\xi'$ to $B$. 
Hence it suffices to see $\xi'$ is definably $s$-isomorphic to some PL microbundle. 
Therefore we assume from the beginning $B$ is a compact PL manifold with boundary of dimension $n$ and $E$ is a definable $C^0$ 
manifold with boundary. 
Let $\tilde B$ be a small definable open neighborhood of $B$ in $R^n$ and $\tilde E$ a definable $C^0$ manifold including  $E$ 
defined by some definable retraction $\tilde B\to B$ as above. \par
  Let $E\subset R^m$. 
We can replace $E$ with the graph of $j$. 
Hence we assume $E\subset B\times R^m$ and $j$ is the restriction to $E$ of the projection $B\times R^m\to B$. 
Set $E_b=\{x\in R^m:(b,x)\in E\}$ for each $b\in B$. 
Then $E=\cup_{b\in B}\{b\}\times E_b$, and moving parallel each $E_b$ in $R^m$ we suppose $\Im i=B\times\{0\}$. \par
  Let us naturally extend the definition of the tangent microbundle of a $C^0$ manifold to that of a definable $C^0$ manifold 
with boundary. 
The tangent microbundle of $E$ $\eta_1:E\overset i_1\to\longrightarrow E\times\tilde E\overset j_1\to\longrightarrow E$ is 
defined by $i_1(b,x)=(b,x,b,x)$ for $(b,x)\in E$ and $j_1(b,x,b',x')=(b,x)$ for $(b,x,b',x')\in E\times\tilde E$. 
Set $V_2=\cup_{b\in B}\{b\}\times E_b\times\tilde B\times E_b$ and $V_3=\cup_{b\in B}\{b\}\times\{0\}\times\tilde B\times E_
b$, and define definable microbundles $\eta_2:E\overset i_2\to\longrightarrow V_2\overset j_2\to\longrightarrow E$ and $\eta_3
:B\times\{0\}\overset i_3\to\longrightarrow V_3\overset j_3\to\longrightarrow B\times\{0\}$ by $i_2(b,x)=(b,x,b,x)$ for $(b,x)
\in E$, $j_2(b,x,b',x')=(b,x)$ for $(b,x,b',x')\in V_2$, $i_3=i_2|_{B\times\{0\}}$ and $j_3=j_2|_{V_3}$. 
Then $\eta_1$ is definably isomorphic to $\eta_2$, the restriction of $\eta_2$ to $B\times\{0\}$ is $\eta_3$, and $\eta_3$ is 
regarded as definably isomorphic to the definable microbundle $B\ni b\to(b,0,b)\in E\times R^n\ni(b,x,b')\to b\in B$, which 
is the Whitney sum of $\xi$ and the trivial microbundle $B\to B\times R^n\to B$. 
Thus $\xi$ is regarded as definably $s$-isomorphic to $\eta_1|_{B\times\{0\}}$. \par
  Let $U$ be the closed $\epsilon$-neighborhood of $B\times\{0\}$ in $E$ for small $\epsilon>0$. 
Then $U$ is a compact definable $C^0$ manifold with boundary. 
Let $\pi:X\to U$ be a definable triangulation such that $\pi^{-1}(B\times\{0\})$ is a polyhedron (remember that a definable 
polyhedral $C^0$ manifold possibly with boundary is a PL manifold possibly with boundary), and $\tau_X$ denote the tangent 
microbundle of $X$. 
Then $\tau|_{\pi^{-1}(B\times\{0\})}$ is a PL microbundle and definably isomorphic to the induced microbundle $(\pi|_{\pi
^{-1}(B\times\{0\})}\allowmathbreak)^*\eta_1|_{B\times\{0\}}$ by invariance of tangent microbundles. 
Let $\Pi:\pi^{-1}(B\times\{0\})\times[0,\,1]\to B\times\{0\}$ be a definable isotopy of $\pi|_{\pi^{-1}(B\times\{0\})}$ through 
homeomorphisms such that $\Pi(\cdot,1)$ is PL (the complement of theorem 2.1). 
Then by the covering homotopy property of topological microbundles, which is proved in the definable $R$-case in the same way, 
$(\pi|_{\pi^{-1}(B\times\{0\})})^*\eta_1|_{B\times\{0\}}$ and $\Pi(\cdot,1)^*\eta_1|_{B\times\{0\}}$ are definably isomorphic. 
Hence $\Pi(\cdot,1)^*\eta_1|_{B\times\{0\}}$ is definably isomorphic to the PL microbundle $\tau_X|_{\pi^{-1}(B\times\{0\})}$. 
Therefore $\eta_1|_{B\times\{0\}}$ is definably isomorphic to the PL microbundle $(\Pi(\cdot,1)^{-1})^*\tau_X|_{\pi^{-1}(B\times
\{0\})}$. 
\qed
\enddemo
\head \S 3. differential topology in o-minimal structure \endhead
In this section we consider differential topology in o-minimal structure over $R$. 
We fix an o-minimal structure over $R$ and one over $\R$, which expand the semialgebraic structure,  such that for any definable set 
$X$ in $\R^n$ there exists a definable set $Y$ in $R^n$ such that $Y\cap\R^n=X$. 
We call the smallest $Y$ the $R$-{\it extension} of $X$ and denote by $X_R$. 
Given a definable $C^0$ map $f:X\to Y$ between definable sets in $\R^n$, then its $R$-{\it extension} $f_R:X_R\to Y_R$ is naturally 
defined. 
An example of such an o-minimal structure over $\R$ is the semialgebraic structure. 
\proclaim{Remark 3.1} 
Let $r$ be $0,...$ or $\infty$. 
A semialgebraic $C^r$ version of lemma 1.1 is that given an $R$-semialgebraic $C^r$ map $f:X\to Y$ between $R$-semialgebraic 
$C^r$ manifolds then there exist an $\R$-semialgebraic $C^r$ map $g:U\to V$ between $\R$-semialgebraic $C^r$ manifolds and 
$R$-semialgebraic $C^r$ diffeomorphisms (homeomorphisms if $r=0$) $\pi:X\to U_R$ and $\tau:Y\to V_R$ such that $\tau\circ f=g_R
\circ\pi$. 
This is not correct. 
Indeed the cardinal number of the $R$-semialgebraic $C^r$ right-left equivalence classes of all $R$-semialgebraic $C^r$ maps 
between $R^3$ is $\# R$ by [T], which treated only the case of $\R$ but whose proof works for any $R$. 
\endproclaim
If we treat only manifolds but not maps, a reduction to the $\R$-case is possible as follows, which was proved by [C-S$_1$] in the semialgebraic structure and by [E] in general. 
A {\it Nash manifold} is a semi-algebraic $C^\infty$ manifold, and a {\it Nash map} between Nash manifolds is a semi-algebraic $C^\infty$ map. 
\proclaim{Lemma 3.2}
Given a definable $C^r$ manifold $M$ in $R^n$, $r>0\in\N$, then there exists uniquely a Nash manifold $N$ in $\R^n$ such that $M$ and 
$N_R$ are definably $C^r$ diffeomorphic. 
Here uniqueness means that if $N'$ is another Nash manifold with the same properties then $N$ and $N'$ are Nash diffeomorphic.
\endproclaim
As in [S$_1$], [C-S$_{1,2}$] and [E] we see 
\proclaim{Lemma 3.2$'$}
Assume the above $M$ is non-compact. 
Then there exists uniquely a compact definable $C^r$ manifold with boundary $N$ in $R^n$ such that $M$ and $\Int N$ are definably 
$C^r$ diffeomorphic. 
\endproclaim
Let $X$ and $Y$ be definable sets in $R^n$, $X',X_i\subset X$ and $Y_i\subset Y$ a finite number of definable subsets and 
$\phi:X'\to Y$ a definable $C^0$ map. 
Let $\Hom_R(X,Y)$ denote the definable homotopy classes of definable $C^0$ maps from $X$ to $Y$. 
We also define $\Hom_R(X,Y)_\phi,\Iso_R(X,Y),$ $\Iso_R(X,Y)_\phi,\ \Aiso_R(X,Y)$, $\Aiso_R(X,Y)_\phi$ and the relative classes 
$\Hom_R(X,X_i;Y,Y_i),...$ in the same way as in the PL case. 
Then we have the following lemma, which generalizes the main theorems on homotopy in [D-K]. 
We apply theorem 2.1 and its complement many times in the proof in order to reduce the problem to the PL case. 
\proclaim{Lemma 3.3} 
Let $X,X',X_i,Y,Y_i\subset\R^n$. 
The following former four natural maps are bijective, and if $Y$ is locally closed in $R^n$ then the latter two ones are injective. 
$$
\gather
\Hom_\R(X,Y)\to\Hom_R(X_R,Y_R),\ \Hom_\R(X,Y)_\phi\to\Hom_R(X_R,Y_R)_{\phi_R},\\
\Iso_\R(X,Y)\to\Iso_R(X_R,Y_R),\ \Iso_\R(X,Y)_\phi\to\Iso_R(X_R,Y_R)_{\phi_R},\\
\Aiso_\R(X,Y)\to\Aiso_R(X_R,Y_R),\ \Aiso_\R(X,Y)_\phi\to\Aiso_R(X_R,Y_R)_{\phi_R}.
\endgather
$$
The relative case also holds. 
\endproclaim
\demo{Proof}
{\it Surjectivity of $\Hom_\R(X,Y)\to\Hom_R(X_R,Y_R)$.} 
Let $f:X_R\to Y_R$ be a definable $C^0$ map. 
By the triangulation theorem of definable sets we can assume $X$ and $Y$ are the unions of some open simplexes in some finite simplicial 
complexes $K$ and $L$ in $\R^n$, respectively, such that $\overline X=|K|$ and $\overline Y=|L|$. 
First we reduce the problem to the case where $X$ and $Y$ are compact. \par
  As in the proof of theorem 2.5, for each $\sigma\in K$ with $X\cap\Int\sigma=\emptyset$, let $\pi_{\sigma,t}:\overline X-\Int\sigma\to\overline X-\Int\sigma,
\ 0\le t\le1$, be a semialgebraic  isotopy of the identity map of $\overline X-\Int\sigma$ such that $\pi_{\sigma,t}=\id$ outside 
of $\Int|\st(\sigma,K)|$ for each $t\in[0,\,1]$ and $\overline X-\overline{\Im\pi_{\sigma,1}}=\Int|\st(v_\sigma,K')|$, where 
$K'$ is the barycentric subdivision of $K$ and $v_\sigma$ is the barycenter of $\sigma$. 
Set $\pi_t=\pi_{\sigma_1,t}\circ\cdots\circ\pi_{\sigma_k,t}$ on $X$ for $\{\sigma_1,...,\sigma_k\}=\{\sigma\in K:X\cap\Int\sigma
=\emptyset\}$. 
Then $\pi_t:X\to X$ is a semialgebraic isotopy of the identity map of $X$, $\overline{\Im\pi_t}\subset X$ for $t\in(0,\,1]$, and 
$f\circ\pi_{tR},\ 0\le t\le1$, is a definable homotopy of $f$. 
Hence we can consider $f|_{\overline{\Im\pi_{1R}}}:\overline{\Im\pi_{1R}}\to Y_R$ in place of $f:X_R\to Y_R$. 
Then $f$ is extensible to a definable $C^0$ map $\overline X_R\to Y_R$ and we can assume $X$ is compact. 
By the same reason we can suppose $\Im f$ does not intersect with some semialgebraic neighborhood of $\overline{Y_R}-Y_R$ in 
$\overline{Y_R}$ defined by polynomials with coefficients in $\R$ and, hence, $Y$ is compact. \par
  The simplicial approximation theorem does not holds over general $R$. 
Its proof, e.g., in [H], however, is available. 
We can prove by a similar method that $f$ is definably homotopic to a PL map as follows, which implies surjectivity of the map $\Hom_\R
(X,Y)\to\Hom_R(X_R,Y_R)$ by lemma 1.2. \par 
  Let $\tau:Z\to X_R$ be a definable triangulation of $X_R$ such that $(f\circ\tau)^{-1}(\sigma)$ are polyhedra for $\sigma\in L_R$ 
and $P$ a simplicial decomposition of $Z$ such that each $(f\circ\tau)^{-1}(\sigma)$ is the union of some simplexes in $P$. 
Then we can replace $X_R,K_R$ and $f$ with $Z,P$ and $f\circ\tau$ because $\tau$ is definably isotopic to a PL homeomorphism 
through homeomorphisms (theorem 2.1 and its complement), $P$ is simplicially isomorphic to $K_{1R}$ for some simplicial complex 
$K_1$ over $\R$ and $|K_1|$ is PL homeomorphic to $X$. 
Hence we assume from the beginning for each $\sigma\in L_R$, $f^{-1}(\sigma)$ is the union of some simplexes in $K_R$. 
Then for each $\sigma\in L_R$, $f^{-1}(\Int\sigma)$ is the union of some open simplexes in $K_R$. 
Hence for each $\sigma\in K_R$ there exists uniquely $\delta_\sigma\in L_R$ such that $f(\Int\sigma)\subset\Int\delta_\sigma$. 
For each $v\in K^0$, let $v_L$ denote the barycenter of $\delta_v$. 
Let $K'$ denote the barycentric subdivision of $K$ and define $\delta_\sigma$ and $v_L$ for $\sigma\in K'_R$ and $v\in K^
{\prime 0}_R$ in the same way. 
Given a simplex $v_1\cdots v_k$ in $K'_R$ spanned by vertices $v_1,...,v_k$, then $f(\Int v_1\cdots v_k)$ is included in $\Int
\delta_{v_1\cdots v_k}$ and, hence, the map carrying $v_1,...,v_k$ to $v_{1L},...,v_{kL}$ in order is linearly extended to a map from $v_1\cdots 
v_k$ to $\delta_{v_1\cdots v_k}$. 
Here the restriction of the extension to a face of $v_1\cdots v_k$ coincides with the linear map defined by the face in the 
same way. 
Hence a PL map $f_1:X_R\to Y_R$ is well-defined so that $f_1$ is linear on each simplex in $K'_R$ and $f_1(v)=v_L$ for $v\in K
^{\prime 0}_R$, and, moreover, a definable homotopy $f_t:X_R\to Y_R,\ 0\le t\le1$, of $f$ is defined by 
$$
f_t(x)=tf_1(x)+(1-t)f(x)\quad\text{for}\ (x,t)\in X_R\times[0,\,1]. 
$$
Thus $f$ is definably homotopic to a PL map, and surjectivity is proved. 
Note that each vertex $v$ of $K'_R$ is the barycenter of some $\sigma\in K_R$, $v_L$ is the barycenter of $\delta_\sigma$ and, 
hence, $f_t(\Int\sigma)\subset\Int\delta_\sigma$. 
In particular, if the restriction of $f$ to a subcomplex $K_1$ of $K_R$ is a simplicial map to some subcomplex of $L_R$ then 
$f=f_t$ on $|K_1|$. 
(Here we need the assumption that for each $\sigma\in L_R$, $f^{-1}(\sigma)$ is the union of some simplexes in $K_R$.) \par
  Moreover we generalize the last statement as follows for subsequence applications. 
If the restriction of $f$ to a closed subpolyhedron $Z$ of $X_R$ is PL and if $f^{-1}(\sigma)$ is a 
polyhedron for each $\sigma\in L_R$ then we can choose $f_t$ so that $f=f_t$ on $Z$. 
Indeed, subdivide $K_R$ and assume the restriction of $f$ to each simplex in a subcomplex $K_1$ of $K_R$ with $|K_1|=Z$ is a 
linear map into some simplex in $L_R$ and for each $\sigma\in L_R$, $f^{-1}(\sigma)$ is the union of some simplexes in $K_R$. 
Then we can define a PL map $f_1:X_R\to Y_R$ so that $f_1$ is linear on each simplex in $K'_R$, $f_1(v)=v_L$ for $v\in K_R^
{\prime 0}-K'_1$ and $f_1(v)=f(v)$ for $v\in K_R^{\prime 0}\cap K'_1$. 
For such $f_1$ we define $f_t,\ 0\le t\le 1$, in the same way as above. \par
  {\it Surjectivity of $\Hom_\R(X,X_i;Y,Y_i)\to\Hom_R(X_R,X_{iR};Y_R,Y_{iR})$.} 
Let $f$ be a definable $C^0$ map from $(X_R,X_{iR})_i$ to $(Y_R,Y_{iR})_i$. 
We repeat the above proof. 
First we assume $X,X_i,Y,Y_i$ are the unions of some open simplexes in finte simplicial complexes $K$ and $L$. 
We reduce the problem, secondly, to the case where $X$ and $Y$ are compact and coincide with $|K|$ and $|L|$, respectively, and, 
thirdly, to the case where for each $\sigma\in L_R$, $f^{-1}(\sigma)$ is the union of some simplexes in $K_R$. 
Lastly, we define a definable homotopy $f_t:X_R\to Y_R,\ 0\le t\le1$, of $f$ to a PL map. 
Then by the above method of construction of $f_t$, $f_t(X_{iR})\subset Y_{iR}$, i.e., $f_t,\ 0\le t\le1$, is a definable 
homotopy from $(X_R,X_{iR})_i$ to $(Y_R,Y_{iR})_i$. 
Thus we can assume $f$ is PL from the beginning. \par
  Let $K_1$ and $L_1$ be simplicial subdivisions of $K_R$ and $L_R$, respectively, such that $f:K_1\to L_1$ is simplicial. 
As in the proof of lemma 1.2 we find PL isotopies $\pi_t$ of $X_R$ and $\tau_t$ of $Y_R,\ 0\le t\le1$, through homeomorphisms 
preserving $K_R$ and $L_R$, respectively, such that $\pi_1$ and $\tau_1$ are isomorphisms from $K_1$ and $L_1$ to the $R$-extensions 
of some simplicial subdivisions of $K$ and $L$ respectively. 
Then $\tau_t\circ f\circ\pi_t^{-1},\ 0\le t\le1$, is a PL homotopy of $f$ from $(X_R,X_{iR})_i$ to $(Y_R,Y_{iR})_i$, and $\tau
_1\circ f\circ\pi^{-1}_1$ is the $R$-extension of some PL map from $(X,X_i)_i$ to $(Y,Y_i)_i$. \par
  {\it Injectivity of $\Hom_\R(X,Y)\to\Hom_R(X_R,Y_R)$.}
Let $f,g:X\to Y$ be definable $C^0$ maps such that $f_R,g_R:X_R\to Y_R$ are definably homotopic. 
Let $F:X_R\times[0,\,1]\to Y_R$ be a definable $C^0$ map such that $F(\cdot,t)=f_R(\cdot)$ for $t\in[0,\,1/4]$ and $F(\cdot,t)=
g_R(\cdot)$ for $t\in[3/4,\,1]$. 
We will see then $f$ and $g$ are definably homotopic. 
By the same reason as above we reduce the problem to the case where $X$ and $Y$ are compact polyhedra, and by the simplicial 
approximation theorem over $\R$ we assume $f$ and $g$ are PL. 
Then by lemma 1.2 it suffices to find a definable homotopy $F_s:X_R\times[0,\,1]\to Y_R$, $0\le s\le 1$, of $F$ such that $F_1$ 
is PL and $F_s=F$ on $X_R\times\{0,1\}$ for any $s$. 
Thus what we prove is the following statement. \par
$(*)$ Let $U_1\subset U_2\subset U$ and $V$ be compact polyhedra over $R$ such that $U_2$ is a neighborhood of $U_1$ in $U$ and 
$h:U\to V$ a definable $C^0$ map such that $h|_{U_2}$ is PL. 
Then there exists a definable homotopy $h_t:U\to V,\ 0\le t\le 1$, of $h$ such that $h_1$ is PL and $h_t=h$ on $U_1$ for any $t$. 
\par
  Existence of such a homotopy without the last condition is already shown. 
We need some additional arguments for this condition. 
Let $K$ and $L$ be simplicial decompositions of $U$ and $V$, respectively, such that $U_1,U_2,h(U_1)$ and $h(U_2)$ are the 
unions of some simplexes in $K$ and $L$, respectively, and $K|_{U_1}$ is full in $K$. 
Here we can choose $K$ and $L$ so that $h|_{U_2}:K|_{U_2}\to L|_{h(U_2)}$ is simplicial for the following reason. 
First subdivide $K$ so that $h$ is linear on each simplex in $K|_{U_2}$. 
Secondly, subdivide $L$ so that for each $\sigma\in K|_{U_2}$, $h(\sigma)$ is the union of some simplexes in $L$. 
Then $\{\sigma\in K:\Int\sigma\cap U_2=\emptyset\}\cup\{\sigma_K\cap h^{-1}(\sigma_L):\sigma_K\in K|_{U_2},\sigma_L\in L\}$ is 
a cell complex, whose underlying polyhedron is $U$. 
Lastly, subdivide the cell complex to a simplicial complex without introducing new vertices (Proposition 2.9 in [R-S]) and 
replace $K$ with the subdivision. 
Then $h|_{U_2}:K|_{U_2}\to L|_{h(U_2)}$ is simplicial. \par
  Next we shrink $U_2$. 
Let $\xi:K\to\{0,1,[0,\,1]\}$ be the simplicial map defined by $\xi(v)=0$ for $v\in(K|_{U_1})^0$ and $\xi(v)=1$ for $v\in K^0-U
_1$. 
Then $\xi^{-1}(0)=U_1$ by fullness of $K|_{U_1}$. 
(Note fullness continues to hold after subdivision.) 
Set $\epsilon=1/2$. 
Then $\xi^{-1}([0,\,\epsilon])\subset U_2$, and $(\xi^{-1}([\epsilon/2,\,\epsilon]),\xi^{-1}([\epsilon/2,\,\epsilon])\cap\sigma)
_{\sigma\in K}$ is PL homeomorphic to $[0,\,1]\times(\xi^{-1}(\epsilon),\xi^{-1}(\epsilon)\cap\sigma)_{\sigma\in K}$ by the 
usual arguments of PL topology. 
We consider $\xi^{-1}([0,\,\epsilon])$ in place of $U_2$ and use the last triviality. \par
  We will construct a definable triangulation $\tau:Z\to U$ and a simplicial decomposition $P$ of $Z$ such that for $\sigma_L\in L$ and $\sigma_K\in K$, $(h\circ\tau)^
{-1}(\sigma_L)$ and $\tau^{-1}(\sigma_K\cap U_1)$ are the unions of some simplexes in 
$P$ and $\tau|_{\tau^{-1}(U_1)}:P|_{\tau^{-1}(U_1)}\to K|_{U_1}$ is an isomorphism. 
Let $\tau_\epsilon:Z_\epsilon\to\xi^{-1}([\epsilon,\,1])$ be a definable triangulation such that $\tau^{-1}_\epsilon(\xi^{-1}(
\epsilon)\cap\sigma_K)$ and $(h\circ\tau_\epsilon)^{-1}(\sigma_L)$ are polyhedra for $\sigma_K\in K$ and $\sigma_L\in L$. 
Apply theorem 2.1 and its complement to $\tau_\epsilon|_{(\xi\circ\tau_\epsilon)^{-1}(\epsilon)}:\big((\xi\circ\tau_\epsilon)^{-1}
(\epsilon),\tau^{-1}_\epsilon(\xi^{-1}(\epsilon)\cap\sigma)\big)_{\sigma\in K}\to(\xi^{-1}(\epsilon),\xi^{-1}(\epsilon)\cap\sigma)
_{\sigma\in K}$. 
Then $\tau_\epsilon$ is extended to a definable triangulation $\tau_{\epsilon/2}:Z_{\epsilon/2}\to\xi^{-1}([\epsilon/2,\,1])$ 
so that $\tau_{\epsilon/2}|_{(\xi\circ\tau_{\epsilon/2})^{-1}(\epsilon/2)}:(\xi\circ\tau_{\epsilon/2})^{-1}(\epsilon/2)\to\xi
^{-1}(\epsilon/2)$ is PL and $\tau^{-1}_{\epsilon/2}(\xi^{-1}([\epsilon/2,\,\epsilon])\cap\sigma)$ are 
polyhedra for $\sigma\in K$. 
Then $(h\circ\tau_{\epsilon/2})^{-1}(\sigma)$ are polyhedra for $\sigma\in L$ because for $\sigma_L\in L$, $\tau^{-1}_{\epsilon/2}(\xi^{-1}([
\epsilon/2,\,\epsilon])\cap h^{-1}(\sigma_L))$ are the unions of $\tau^{-1}_{\epsilon/2}(\xi^{-1}([
\epsilon/2,\,\epsilon])\cap\sigma_K)$ for some $\sigma_K\in K$. 
Let $P_{\epsilon/2}$ be a simplicial decomposition of $Z_{\epsilon/2}$ such that for $\sigma_L\in L$ and $\sigma_K\in K$, $(h\circ\tau_{\epsilon/2})^{-1}(\sigma_L)$ 
and $\tau^{-1}_{\epsilon/2}(\sigma_K\cap\xi^{-1}(\epsilon/2))$ are the unions of some 
simplexes in $P_{\epsilon/2}$. 
Subdivide the cell complex $\{\sigma_K\cap\xi^{-1}(0),\sigma_K\cap\xi^{-1}([0,\,\epsilon/2]),\tau_{\epsilon/2}(\sigma_P)\cap
\xi^{-1}(\epsilon/2):\sigma_K\in K,\sigma_P\in P_{\epsilon/2}\}$ to a simplicial complex $\tilde K$ without introducing new 
vertices. 
Note $|\tilde K|=\xi^{-1}([0,\,\epsilon/2])$, $\tilde K|_{U_1}=K|_{U_1}$ and $\tilde K|_{\xi^{-1}(\epsilon/2)}=\{\tau_{\epsilon
/2}(\sigma)\cap\xi^{-1}(\epsilon/2):\sigma\in P_{\epsilon/2}\}$. 
Paste $P_{\epsilon/2}$ with $\tilde K$ through $\tau_{\epsilon/2}|_{(\xi\circ\tau_{\epsilon/2})^{-1}(\epsilon/2)}$. 
Then we have a finite simplicial complex $P$ and a definable triangulation $\tau:Z=|P|\to U$ such that for $\sigma_L\in L$ and $\sigma_K\in K$, $(h\circ\tau)^{-1}(
\sigma_L)$ and $\tau^{-1}(\sigma_K\cap U_1)$ are the unions of some simplexes in $P$ 
and $\tau|_{\tau^{-1}(U_1)}:P|_{\tau^{-1}(U_1)}\to K|_{U_1}$ is an isomorphism. \par
  The rest of proof is the same as in the proof of surjectivity. 
First replacing $U$ and $K$ with $Z$ and $P$ we reduce the problem to the case where for each $\sigma\in L$, $h^{-1}(\sigma)$ is 
the union of some simplexes in $K$ and $h|_{U_1}:K|_{U_1}\to L|_{h(U_1)}$ is simplicial. 
Next define a definable homotopy $h_t:U\to V$, $0\le t\le1$, of $h$ so that $h_1$ is PL in the same way. 
Then $h_t=h$ on $U_1$ by the simplicial map property of $h|_{U_1}$ and the definition of $h_t$. 
Thus injectivity is proved. \par
  {\it Injectivity of $\Hom_\R(X,X_i;Y,Y_i)\to\Hom_R(X_R,X_{iR};Y_R,Y_{iR})$.}
We only need to generalize $(*)$ to the relative case as follows. 
Given $W_i\subset U$ and $V_i\subset V$ a finite number of unions of a finite number of open simplexes such that $h(W_i)\subset 
V_i$, then we can choose $h_t$ so that $h_t(W_i)\subset V_i$ for any $t$. 
That is clear as shown in the proof of surjectivity of $\Hom_\R(X,X_i;Y,Y_i)\to\Hom_R(X_R,X_{iR};Y_R,Y_{iR})$. \par
  {\it Bijectivity of $\Hom_\R(X,X_i;Y,Y_i)_\phi\to\Hom_R(X_R,X_{iR};Y_R,Y_{iR})_{\phi_R}$.}
We can assume $X$ and $Y$ are bounded in $\R^n$ and, moreover, $X'$ is closed in $X$ and for the following reason. 
Since we are interested in the case of $\Hom_R(X_R,X_{iR};Y_R,Y_{iR})_{\phi_R}\not=\emptyset$, we assume $\phi_R$ is extensible to a definable $C^0$ map $(X_R,X_{iR})_i\to(Y_R,Y_{iR})_i$. 
Then $\phi_R$ is uniquely extensible to the closure of $X'_R$ in $X_R$. 
Let $\overline\phi_R$ denote the extension. 
Clearly $\graph\overline\phi_R$ is closed in $X_R\times Y_R$ and the closure of $X'_R$ in $X_R$ includes the closure of $X'$ in $X$. 
Therefore $\graph\overline\phi_R$ includes the closure of $\graph\phi$ in $X\times Y$. 
Namely $\phi$ is extensible to the closure of $X'$ in $X$. 
Therefore we can assume $X'$ is closed in $X$. \par
  We reduce the problem to the compact PL case one by one. 
First we reduce to the case where $\phi$ is an imbedding as follows. 
Let $\tilde\phi:X\to\R^n$ be a bounded definable $C^0$ extension of $\phi$. 
Replace $X$ with $\graph\tilde\phi$. 
Then we assume $\tilde\phi$ is extensible to a definable $C^0$ map $\hat\phi:\overline X\to\R^n$. 
Using Thom's transversality theorem (Theorem II.5.4), by usual arguments of singularity theory (see \S II.5 and \S II.6 in [S$_3$]) we 
find a definable $C^1$ map $\alpha':\overline X\to\R^N$ (i.e., a map extensible to a definable $C^1$ map from a definable 
open neighborhood of $\overline X$ in $\R^n$ to $\R^N$) for some large integer $N$ such that $\alpha'=0$ on $\overline{X'}$ 
and $\alpha'|_{\overline X-\overline{X'}}$ is a $C^0$ imbedding into $\R^N-\{0\}$. 
Set $\alpha=(\hat\phi,\alpha')|_X,\ Z=\Im\alpha,\ Z'=\alpha(X')$ and $Z_i=\alpha(X_i)$, and define a definable $C^0$ 
map $\psi:Z'\to Y$ to be the restriction to $Z'$ of the projection $\R^n\times\R^N\to\R^n$. 
Then $\psi$ is an imbedding and, moreover, extensible to a definable $C^0$ imbedding of $\overline{Z'}$ into $\overline Y$, and 
the map $\Hom_R(Z_R,Z_{iR};Y_R,Y_{iR})_{\psi_R}\ni G\to G\circ\alpha_R\in\Hom_R(X_R,X_{iR};Y_R,Y_{iR})_{\phi_R}$ is bijective. 
Hence we can replace $(X,X',X_i)_i$ with $(Z,Z',Z_i)_i$ and assume $\phi$ and its $C^0$ extension to $\overline{X'}$ are 
imbeddings into $Y$ and $\overline Y$ respectively. \par
  Secondly, by the triangulation theorem of definable sets, theorem 2.1 and its complement we suppose $X,X',X_i,Y$ and $Y_i$ are the 
unions of some open simplexes in finite simplicial complexes $K$ and $L$ such that $|K|=\overline X$, $|L|=\overline Y$ and $K|
_{\overline{X'}}$ is a full subcomplex of $K$, and the extension of $\phi$ to $\overline{X'}$ is a PL imbedding into 
$\overline Y$. 
(Note the original $\phi$ is not necessarily triangulable.) \par
  Let $\pi_{\sigma,t}:\overline X-\Int\sigma\to\overline X-\Int\sigma,\ 0\le t\le 1$, be the semialgebraic isotopy defined in 
the proof of surjectivity of $\Hom_{\R}(X,Y)\to\Hom_R(X_R,Y_R)$ and $K'$ the barycentric subdivision of $K$. 
Set $\{\sigma_1,...,
\sigma_k\}=\{\sigma\in K:X\cap\Int\sigma=\emptyset\}$ ordered so that $\dim\sigma_i\le\dim\sigma_{i+1}$. 
Then we can choose $\pi_{\sigma_i,t}$ by downward induction on $i$ by fullness of $K|_{\overline{X'}}$ in $K$ so that 
$$
\overline{X-\Im\pi_{\sigma_i,1}\circ\cdots\circ\pi_{\sigma_k,1}}=\cup_{j=i}^k|\st(v_{\sigma_j},K')|,
$$
where $v_\sigma$ denotes the barycenter of $\sigma$. 
Hence if we set $\pi_1=\pi_{\sigma_1,1}\circ\cdots\circ\pi_{\sigma_k,1}$ as before then $\pi_1$ is a semialgebraic $C^0$ imbedding of 
$X$ into $X$,  
$$
\overline{X-\Im\pi_1}=\cup\{|\st(\sigma,K')|:\sigma\in K',\,\Int\sigma\cap X=\emptyset\} 
$$
and $\overline{\Im\pi_1}$ is the union of some simplexes in $K'$. 
It follows that $\cup\{|\st(\sigma',K')|\cap\overline{\Im\pi_1}:\sigma'\in K',\,\Int\sigma'\subset\sigma\cap\Im\pi_1\}$ is a neighborhood of $\sigma\cap\overline{\Im\pi_1}$ in $\overline{\Im\pi_1}$ for any $\sigma\in K$ with $\Int\sigma\subset X$. 
Therefore replacing $X,\,\phi$ and $K$ with $\Im\pi_1,\,\phi\circ\pi_1^{-1}$ and $K'|_{\overline{\Im\pi_1}}$ and keeping the same 
notation we, thirdly, assume 
$$
\sigma\cap(\overline{X'}-X)\subset\overline{\sigma\cap X'}\quad\text{for}\ \sigma\in K\ \text{with}\ \Int\sigma\subset X,\tag$**$
$$
namely, there does not exist a simplex in $K$ which touches $\overline{X'}$ but not $X'$ and whose interior is included in $X$. 
\par
  Though reduction to the compact PL case is not yet complete, we consider surjectivity. \par
  {\it Proof of surjectivity.}
Let $f:(X_R,X_{iR})_i\to(Y_R,Y_{iR})_i$ be a definable $C^0$ map with $f|_{X'_R}=\phi_R$. 
Then $f$ is definably homotopic to a definable $C^0$ map extensible to a definable $C^0$ map $\overline X_R\to\overline Y_R$ 
through a homotopy fixing $X'_R$ for the following reason. 
Let $U$ be an open definable neighborhood of $X'_R$ in $X_R$ such that $f|_U$ is extensible to a definable $C^0$ map from $\overline U$ 
to $\overline Y$, define $K'$ as above and let $\xi$ be a positive definable $C^0$ function on $X'_R$ with $\xi\le 1$. \par
  We will define a definable homeomorphism $\tau_\xi:X_R\to X_R$ by which $\overline U$ becomes a neighborhood of $\overline{X'}$ in $\overline X$. 
Let $N(\overline{X'},K)$ denote the simplicial neighborhood of $\overline{X'}$ in $K$---the subcomplex of $K$ generated by 
$\sigma\in K$ with $\sigma\cap\overline{X'}\not=\emptyset$---and $\partial N(\overline{X'},K)$ the subcomplex of $N(
\overline{X'},K)$ consisting of $\sigma$ with $\sigma\cap\overline{X'}=\emptyset$. 
Set $\tau_\xi=\id$ on $X'_R\cup(X_R\cap|\partial N(\overline{X'},K)|_R)\cup(X_R-|N(\overline{X'},K)|_R)$. 
By $(**)$ $X_R\cap|N(\overline{X'},K)|_R-X'_R-|\partial N(\overline{X'},K)|_R$ is the disjoint union of the interiors of joins 
$\sigma*\sigma'$ such that $\sigma\in\partial N(\overline{X'},K)_R,\,\sigma'\in K_R|_{\overline{X'_R}},\ \sigma*\sigma'\in K_R$ 
and $\Int\sigma'\subset X'_R$. 
Hence it suffices to define $\tau_\xi$ on such $\Int\sigma*\sigma'$. 
Clearly $\Int\sigma*\sigma'$ is the disjoint union of open segments joining points in $\Int\sigma$ and points in $\Int\sigma'$. 
Let $l(t),\ 0\le t\le 1$, denote the points of the closed segment linearly parameterized by $t$ so that $l(0)\in\Int\sigma'$ 
and $l(1)\in\Int\sigma$. 
Define $\tau_\xi$ to carry linearly the segments $\overline{l(0)l(1/2)}-\{l(0)\}$ and $\overline{l(1/2)l(1)}-\{l(1)\}$ to 
$\overline{l(0)l(\xi(l(0))/2)}-\{l(0)\}$ and $\overline{l(\xi(l(0))/2)l(1)}-\{l(1)\}$ respectively. 
Then $\tau_\xi$ is a well-defined definable homeomorphism of $X_R$ and preserving $\{\Int\sigma:\sigma\in K_R,\,\Int\sigma
\subset X_R\}$, and for sufficiently small $\xi$ 
$$
\tau_\xi(X_R\cap|N(\overline{X'},K')|_R)\subset U.
$$
Fix such a $\xi$ and set $\tau'_t=\tau_{1-t+t\xi},\ 0\le t\le 1$. 
Then $\tau'_t,\ 0\le t\le 1$, is a definable isotopy of $X_R$ fixing $X'_R$ such that 
$$
\tau'_1(X_R\cap|N(\overline{X'},K')|_R)\subset U.
$$
Hence $f\circ\tau'_1$ is extensible to $|N(\overline{X'},K')|_R$. 
After then we can see $f\circ\tau'_1$ is definably homotopic to a definable $C^0$ map extensible to a definable $C^0$ map 
$\overline X_R\to\overline Y_R$ through a homotopy fixing $X'_R$ in the same way as in the proof of surjectivity of $\Hom_\R
(X,Y)\to\Hom_R(X_R,Y_R)$. \par
  Thus we have reduced surjectivity problem to the case where $X$ and $Y$ are compact. 
It suffices then to show $f$ is definably homotopic to a PL map through a homotopy fixing $X'_R$ by lemma 1.2. 
Subdivide $K$ and $L$ so that $\phi:K|_{X'}\to L|_{\phi(X')}$ is simplicial. 
By the triangulation theorem of definable sets, theorem 2.1 and its complement there exists a definable isotopy $\eta_t,\ 0\le t\le 1$, of $X_R$ 
preserving $K_R$ such that $\eta_1(f^{-1}(\sigma)),\ \sigma\in L_R$, are polyhedra. 
Since $f\circ\eta_t\circ f^{-1}|_{f(X'_R)},\ 0\le t\le 1$, is a definable isotopy of $f(X'_R)$ preserving $L_R|_{f(X'_R)}$ we 
can extend it to a definable isotopy $\rho_t,\ 0\le t\le 1$, of $Y$ preserving $L_R$. 
Then $\rho_t\circ f\circ\eta_t^{-1},\ 0\le t\le 1$, is a definable homotopy fixing $X'$ such that $(\rho_1\circ f\circ\eta_1^{-1}
)^{-1}(\sigma)=\eta_1(f^{-1}(\sigma))$ are polyhedra for $\sigma\in L_R$. 
Hence we assume $f^{-1}(\sigma)$ are polyhedra for $\sigma\in L_R$. 
Then as shown at the end of the proof of surjectivity of $\Hom_\R(X,Y)\to\Hom_R(X_R,Y_R)$ there exists a definable homotopy 
$f_t:(X_R,X_{iR})_i\to(Y_R,Y_{iR})_i,\ 0\le t\le 1$, from $f$ to a PL map $f_1$ fixing $X'_R$. \par
  {\it Proof of injectivity.} 
Let $f,g:(X,X_i)_i\to(Y,Y_i)_i$ be  definable $C^0$ maps with $f=g=\phi$ on $X'$ such that $f_R$ and $g_R$ are definably 
homotopic through a homotopy fixing $X'_R$. 
Let $F:(X_R, X_{iR})_i\times[0,\,1]\to(Y_R,Y_{iR})_i\times[0,\,1]$ be a definable $C^0$ map of the form $F(x,t)=(F'(x,t),t)$ 
such that $F'(\cdot,t)=f_R(\cdot)$ for $t\in[0,\,1/4],\ F'(\cdot,t)=g_R(\cdot)$ for $t\in[3/4,\,1]$ and $F'(\cdot,t)=\phi_R(
\cdot)$ on $X'_R$ for any $t$. 
Then by the above proof we can assume $X$ and $Y$ are compact polyhedra and $f$ and $g$ are of class PL, and it suffices to find 
a definable homotopy $F_s:(X_R,X_{iR})_i\times[0,\,1]\to(Y_R,Y_{iR})_i\times[0,\,1],\ 0\le s\le 1$, of $F$ such that $F_1$ is PL 
and $F_s=F$ on $X_R\times\{0,1\}\cup X'_R\times[0,\,1]$. 
Hence we need only to prove the following statement. \par 
  $(*)'$ Let $U_1,U_2,U,W_i,V,V_i$ and $h$ be the same as in $(*)$ in the proof of injectivity of $\Hom_\R(X,Y)\to\Hom_R(X_R,Y_R)$ or its relative case in the proof of injectivity of $\Hom_\R(X,X_i;Y,Y_i)\to\Hom_R(X_R,X_{iR};Y_R,Y_{iR})$. 
Let $U'$ be a closed subpolyhedron of $U$ such that $h|_{U'}$ is a PL imbedding. 
Let $K$ and $L$ be simplicial decompositions of $U$ and $V$, respectively, such that $U_1,U_2,U',h(U_1),h(U_2)$ and $h(U')$ are the unions 
of some simplexes in $K$ and $L$ and $h|_{U_2}:K|_{U_2}\to L|_{h(U_2)}$ and $h|_{U'}:K_{U'}\to L|_{h(U')}$ are simplicial. 
Let $\xi:K\to\{0,1,[0,\,1]\}$ be the simplicial map with $\xi^{-1}(0)=U_1$. 
Then we can choose the homotopy $h_t,\ 0\le t\le 1$, of $h$ in $(*)$ or its relative case so that $h_t=h$ on $U'$. \par
  {\it Proof of $(*)'$.} 
We proceed to prove in a similar way. 
Shrink $U_2$ to $\xi^{-1}([0,\,\epsilon])$, where $\epsilon=1/2$, and set $U_3=\xi^{-1}([\epsilon,\,1])$. 
Let $\tau_\epsilon:Z_\epsilon\to U_3\cup U'$ be a definable triangulation such that $(h\circ\tau_\epsilon)^{-1}(\sigma_L)$ and $\tau^{-1}_\epsilon(\sigma_K)$ are polyhedra for 
$\sigma_L\in L$ and $\sigma_K\in K$. 
Here by theorem 2.1 and its complement we can assume $\tau_\epsilon$ is PL on $\tau_\epsilon^{-1}(U'\cap\xi^{-1}([0,\,\epsilon
/2]))$ since for each $\sigma\in K$, $\sigma\cap U'\cap\xi^{-1}([\epsilon/2,\,\epsilon])$ is a PL ball to which the restriction 
of $\xi$ is trivial. 
Moreover we suppose $\tau_\epsilon$ is globally PL for the following reason. \par
  By theorem 2.1 and its complement we have a definable isotopy $\tau_{\epsilon t},\ 0\le t\le 1$, of $\tau_\epsilon$ through 
homeomorphisms such that $\tau^{-1}_{\epsilon t}(\sigma)=\tau^{-1}_\epsilon(\sigma)$ for $\sigma\in K$, $\tau_{\epsilon 1}$ is 
PL and $\tau_{\epsilon t}=\tau_\epsilon$ on $\tau^{-1}_\epsilon(U'\cap\xi^{-1}([0,\,\epsilon/2]))$. 
(Here the last condition is not stated in the complement. 
It is, however, clearly satisfied by the proof of the complement.) 
Consider two definable isotopies $\tau_\epsilon\circ\tau_{\epsilon t}^{-1},\ 0\le t\le 1$, of $U_3\cap U'$ and $h\circ\tau_\epsilon
\circ\tau^{-1}_{\epsilon t}\circ h^{-1}|_{h(U')},\ 0\le t\le 1$, of $h(U')$, which are preserving $K|_{U_3\cap U'}$ and $L|_{h
(U')}$ respectively. 
Extend them to definable isotopies $\eta_t,\ 0\le t\le 1$, of $U$ preserving $K$ and $\zeta_t,\ 0\le t\le 1$, of $V$ preserving 
$L$ so that $\eta_t=\id$ on $\xi^{-1}([0,\,\epsilon/2])$ and $\zeta_t=\id$ on $h(\xi^{-1}([0,\,\epsilon/2]))$. 
Define $g_t$ to be $\zeta^{-1}_t\circ h\circ\eta_t$ for $t\in[0,\,1]$. 
Then $g_t$, $0\le t\le 1$, is a definable homotopy of $h$, 
$$
\gather
g_t=(h\circ\tau_\epsilon\circ\tau_{\epsilon t}^{-1}\circ h^{-1})^{-1}\circ h\circ\tau_\epsilon\circ\tau^{-1}_{\epsilon t}=h 
\quad\text{on}\ U'\ \text{and}\\
g_t=\id\circ h\circ\id=h\quad\text{on}\ \xi^{-1}([0,\,\epsilon/2]). 
\endgather
$$
Hence we can replace $h$ with $g_1$ though the equality $g^{-1}(g(\xi^{-1}(s)))=\xi^{-1}(s)$ holds for only $s\in[0,\,\epsilon]$. 
Replace $\tau_\epsilon$ with the triangulation $\tau_{\epsilon 1}:Z_\epsilon\to U_3\cup U'$. 
Then $(g_1\circ\tau_{\epsilon 1})^{-1}(\sigma_L)$ and $\tau^{-1}_{\epsilon 1}(\sigma_K)$ 
are polyhedra for $\sigma_L\in L$ and $\sigma_K\in K$ because 
$$
\gather
(g_1\!\circ\!\tau_{\epsilon 1})^{-1}(\sigma_L)\!=\!\tau^{-1}_{\epsilon 1}\!\circ\!\eta_1^{-1}\!\circ\! h^{-1}\!\circ\!\zeta_1(
\sigma_L)\!=\!\tau_{\epsilon 1}^{-1}\!\circ\!\tau_{\epsilon 1}\!\circ\!\tau_\epsilon^{-1}\!\circ\! h^{-1}(\sigma_L)\!=\!(h\!\circ
\!\tau_\epsilon)^{-1}(\sigma_L),\\\tau^{-1}_{\epsilon 1}(\sigma_K)=\tau^{-1}_\epsilon(\sigma_K). 
\endgather
$$
Hence, since $\tau_{\epsilon 1}$ is PL we can suppose $\tau_\epsilon$ is PL from the beginning. \par
  We have simplicial decompositions $P_\epsilon$ of $Z_\epsilon$ and $\tilde K$ of $\xi^{-1}([0,\,\epsilon])$ such that the following map is an isomorphism. 
$$\tau_\epsilon|_{(\xi\circ\tau_\epsilon)^{-1}([0,\,\epsilon])}:P_\epsilon|_{(\xi\circ\tau_\epsilon)^{-1}
([0,\,\epsilon])}\to\tilde K|_{\xi^{-1}(\epsilon)\cup\big(U'\cap\xi^{-1}([0,\,\epsilon])\big)}.
$$
(Note $\tau_\epsilon\big((\xi\circ\tau_\epsilon)^{-1}([0,\,\epsilon])\big)=\xi^{-1}(\epsilon)\cup\big(U'\cap\xi^{-1}([0,\,\epsilon])\big)$.)
Paste $P_\epsilon$ and $\tilde K$ through $\tau_\epsilon|_{(\xi\circ\tau_\epsilon)^{-1}([0,\,\epsilon])}$. 
Then we obtain a definable triangulation $\tau:Z\to U$ of class PL and a simplicial decomposition $P$ of $Z$ such that for $\sigma_L\in L$ and $\sigma_K\in K$, $(h\circ
\tau)^{-1}(\sigma_L)$ and $\tau^{-1}(\sigma_K)$ are the unions of some simplexes in $P$. \par
  In conclusion, replacing $U$ with $Z$ we assume from the beginning $h^{-1}(\sigma)$ is a polyhedron for each $\sigma\in L$. 
Then the required $h_t,\ 0\le t\le 1$, is constructed as stated at the end of proof of surjectivity of $\Hom_\R(X,Y)\to\Hom_R(
X_R,Y_R)$. \par
  {\it Bijectivity of $\Iso_\R(X,Y)\to\Iso_R(X_R,Y_R),\ \Iso_\R(X,Y)_\phi\to\Iso_R(X_R,Y_R)_{\phi_R}$ and their relative maps.}
Consider $\Iso_\R(X,X_i;Y,Y_i)_\phi\to\Iso_R(X_R,X_{iR};Y_R,Y_{iR})_{\phi_R}$ only because the map is the most general. 
As above we can reduce the problem to the case where $X,X'$ and $Y$ are compact polyhedra, $X_i$ and $Y_i$ are the unions of a 
finite number of some open simplexes and $\phi$ is PL. \par
  {\it Proof of surjectivity.}
Let $f:(X_R,X_{iR})\to(Y_R,Y_{iR})$ be a definable $C^0$ imbedding such that $f=\phi_R$ on $X'_R$. 
By lemma 1.2 it suffices to see $f$ is definably isotopic to some PL imbedding through an isotopy fixing $X'_R$. 
Let $K$ and $L$ be simplicial decompositions of $X_R$ and $Y_R$, respectively, such that $X_{iR},\ X'_R,\ Y_{iR}$ and $\phi_R(X'_R)$ are the unions of some open simplexes in $K$ or $L$ and $\phi_R:K|_{X'_R}\to L|_{\phi_R(X'_R)}$ is simplicial. 
First we reduce the problem to the case where $\{f(\sigma):\sigma\in K\}$ are polyhedra. \par
  By the triangulation theorem of definable sets we have a definable isotopy $\pi_t,\ 0\le t\le 1$, of $Y_R$ preserving $L$ such that $\{\pi_1\circ f(\sigma):\sigma\in K\}$ are polyhedra. 
Then $\pi_t\circ f:(X_R,X_{iR})\to(Y_R,Y_{iR}),\ 0\le t\le 1$, is a definable isotopy from $f$ to $\pi_1\circ f$ but not necessarily fixing $X'_{iR}$. 
We need to modify $\pi_t\circ f$. 
Set $\tau_t=f^{-1}\circ\pi_t\circ f$ on $X'_R$ for $t\in[0,\,1]$, which is a well-defined definable isotopy of $X'_R$ preserving $K|_{X'_R}$. 
By the Alexander trick $\tau_t$ is extended to a definable isotopy $\tilde\tau_t,\ 0\le t\le 1$, of $X_R$ preserving $K$ such that $\tilde\tau_0=\id$. 
Then $\pi_t\circ f\circ\tilde\tau_t^{-1},\ 0\le t\le 1$, is a definable isotopy from $f$ through an isotopy fixing $X'_R$ and such that $\{\pi_1\circ f\circ\tilde\tau_t^{-1}(\sigma):\sigma\in K\}$ are polyhedra by the following respective equalities. 
$$
\gather
\pi_0\circ f\circ\tilde\tau_0^{-1}=\id\circ f\circ\id=f,\\
\pi_t\circ f\circ\tilde\tau_t^{-1}=\pi_t\circ f\circ f^{-1}\circ\pi^{-1}_t\circ f=f\quad\text{on}\ X'_R,\\
\pi_1\circ f\circ\tilde\tau_1^{-1}(\sigma)=\pi_1\circ f(\sigma)\quad\text{for}\ \sigma\in K.
\endgather
$$
Hence we can assume each $f(\sigma)$ is a polyhedron. \par
  Apply theorem 2.1 and its complement to $f:(X_R,\sigma)_{\sigma\in K}\to(f(X_R),f(\sigma))_{\sigma\in K}$. 
Then there exists a definable isotopy $\omega_t:(X_R,\sigma)_{\sigma\in K}\to(f(X_R),f(\sigma))_{\sigma\in K}$ from $f$ to some PL homeomorphism through homeomorphisms. 
Here also $\omega_t$ is not necessarily fixing $X'_R$. 
We modify $\omega_t$ by a method similar to the above. 
Set $\xi=\omega_1^{-1}\circ f$ on $X'_R$, which is a PL homeomorphism of $X'_R$ preserving $K|_{X'_R}$. 
Extend $\xi$ to a PL homeomorphism $\tilde\xi$ of $X_R$ preserving $K$ by the Alexander trick. 
Set 
$$
\Xi(x,t)=\cases
(x,0)\quad&\text{for}\ (x,t)\in X_R\times\{0\}\\
(\omega_t^{-1}\circ f(x),t)\quad&\text{for}\ (x,t)\in X'_R\times[0,\,1]\\
(\tilde\xi(x),1)\quad&\text{for}\ (x,t)\in X_R\times\{1\}.
\endcases
$$
Then $\Xi$ is a definable homeomorphism of $X_R\times\{0,1\}\cup X'_R\times[0,\,1]$ preserving $\{\sigma\times\{0\},\sigma\times\{1\},\sigma'\times[0,\,1]:\sigma\in K,\sigma'\in K|_{X'_R}\}$. 
By the Alexander trick $\Xi$ is extended to a definable homeomorphism $\tilde\Xi$ of $X_R\times[0,\,1]$ preserving $K\times\{0,1,[0,\,1]\}$ and of the form $\tilde\Xi(x,t)=(\tilde\Xi'(x,t),t)$ for $(x,t)\in X_R\times[0,\,1]$. 
Consider $\omega_t\circ\tilde\Xi'(\cdot,t)$ in place of $\omega_t(\cdot)$ and denote it by $\eta_t(\cdot)$. 
Clearly $\eta_t,\ 0\le t\le 1$, is a definable isotopy $(X_R,X_{iR})\to(Y_R,Y_{iR})$ from $f$ fixing $X'_R$ and such that $\eta_1$ is PL. 
Indeed 
$$
\gather
\eta_t(x)=\omega_t\circ\omega_t^{-1}\circ f(x)=f(x)\quad\text{for}\ (x,t)\in X'_R\times[0,\,1],\\
\eta_0=\omega_0\circ\id=f,\quad\eta_1(\cdot)=\omega_1\circ\tilde\Xi'(\cdot,1)=\omega_1\circ\tilde\xi(\cdot).
\endgather
$$
Thus $\eta_t,\ 0\le t\le 1$, is the required isotopy, and surjectivity is proved. \par
  {\it Proof of injectivity.}
Let $f,g:(X,X_i)_i\to(Y,Y_i)_i$ be definable $C^0$ imbeddings such that $f_R$ and $g_R$ are definably isotopic through an 
isotopy fixing $X'_R$. 
As above we can reduce the problem to the case where $f$ and $g$ are PL. 
Then by lemma 1.2 it suffices to see $f_R$ and $g_R$ are PL isotopic through an isotopy fixing $X'_R$. 
Let $K$ be a simplicial decomposition of $X_R$ such that $X'_R$ and $X_{iR}$ are the unions of open simplexes in $K$ and $\phi_R$ is linear on each simplex in $K|_{X'_R}$, and $L$ one of $Y_R\times[0,\,1]$ such that $Y_{iR}\times[0,\,1]$ and $\phi_R(\sigma)\times[0,\,1]$ for each $\sigma\in K|_{X'_R}$ are the unions of some open simplexes in $L$. 
Let $F:(X_R,X_{iR})\times[0,\,1]\to(Y_R,Y_{iR})\times[0,\,1]$ be a definable $C^0$ imbedding of the form $F(\cdot,t)=(F'(\cdot,t),t)$ for $t\in[0,\,1]$ such that $F'(\cdot,0)=f_R(\cdot),\ F'(\cdot,1)=g_R(\cdot)$ and $F'(x,t)=\phi_R(x)$ for $(x,t)\in X'\times[0,\,1]$. 
Here we can assume $F'(\cdot,t)=f_R(\cdot)$ for $t\in[0,\,1/3]$ and $F'(\cdot,t)=g_R(\cdot)$ for $t\in[2/3,\,1]$. 
Then $F$ is PL on $X_R\times([0,\,1/3]\cup[2/3,\,1])$ and, hence, $F(\sigma\times([0,\,1/3]\cup[2/3,\,1]))$ are polyhedra for $\sigma\in K$. \par
  We first modify $F$ so that $F(\sigma\times[0,\,1])$ are polyhedra for all $\sigma\in K$. 
Let $p$ denote the restriction to $Y_R\times[0,\,1]$ of the projection $Y_R\times R\to R$. 
Note 
$$
p^{-1}(R-(1/3,\,2/3))\cap F(\sigma\times[0,\,1])=F(\sigma\times([0,\,1/3]\cup[2/3,\,1]))\quad\text{for}\ \sigma\in K. 
$$
Apply theorem 2.2,(3) to $Y_R\times[0,\,1],\ L,\ \{F(\sigma\times[0,\,1]):\sigma\in K\}$ and $p$. 
Then we have a definable homeomorphism $\pi$ of $Y_R\times[0,\,1]$ preserving $L$ such that $\pi=\id$ on $Y_R\times([0,\,1/4]\cup[3/4,\,1])$, $\pi(F(\sigma\times[0,\,1]))$ are polyhedra for $\sigma\in K$ and $p\circ\pi=p$, i.e., $\pi$ is of the form $\pi(\cdot,t)=(\pi'(\cdot,t),t)$ for $t\in[0,\,1]$ (we call $\pi$ level-preserving). \par
  Consider $\pi\circ F$ in place of $F$. 
Then the equality $\pi'\circ F(x,t)=\phi_R(x)$ for $x\in X'_R$ does not necessarily hold though the other properties---$\pi'\circ F(\cdot,0)=f_R(\cdot)$ and $\pi'\circ F(\cdot,1)=g_R(\cdot)$---continue to be true. 
We can modify $F'$ so that that equality also holds in the same way as in the above proof of surjectivity. 
We do not repeat. 
Hence we assume from the beginning that $F(\sigma\times[0,\,1])$ are polyhedra for $\sigma\in K$. \par
  Next we modify $F$ to a PL homeomorphism fixing on $X'_R\times[0,\,1]\cup X_R\times\{0,1\}$. 
Assume by induction that $F$ is PL on $|K^k|\times[0,\,1]$ for some $k\in\N$, and let $\sigma\in K^{k+1}-K^k$ with $\sigma\not\subset X'_R$. 
Then it suffices to modify $F$ so that $F|_{\sigma\times[0,\,1]}$ is PL. 
Consider two definable homeomorphisms $F|_{\sigma\times[0,\,1]}:\sigma\times[0,\,1]\to F(\sigma\times[0,\,1])$ and $\id:F(\sigma\times[0,\,1])\to F(\sigma\times[0,\,1])$ and compact polyhedra $F(\sigma'\times[0,\,1])$ for all proper faces $\sigma'$ of $\sigma$. 
Both are definable triangulations of $p|_{F(\sigma\times[0,\,1])}$ with a family of polyhedra. 
Hence by the complement of theorem 2.2 there is a PL homeomorphism $\omega:\sigma\times[0,\,1]\to F(\sigma\times[0,\,1])$ such that $\omega(\sigma'\times[0,\,1])=F(\sigma'\times[0,\,1])$ for proper faces $\sigma'$ of $\sigma$ and $p\circ\omega=p\circ F$ on $\sigma\times[0,\,1]$, i.e., $\omega$ is level-preserving. 
Here we can assume $\omega=F$ on $\partial(\sigma\times[0,\,1])$ for the following reason. \par
  Set $\xi=\omega^{-1}\circ F$ on $\sigma\times[0,\,1]$. 
Then $\xi$ is a definable homeomorphism of $\sigma\times[0,\,1]$, PL on $\partial(\sigma\times[0,\,1])$, level-preserving and preserving $\{\sigma'\times[0,\,1]:\sigma'\in K,\sigma'\subset\sigma\}$. 
By the Alexander trick there exists a level-preserving PL homeomorphism $\eta$ of $\sigma\times[0,\,1]$ such that $\eta=\xi$ on $\partial(\sigma\times[0,\,1])$. 
Replace $\omega$ with $\omega\circ\eta$. 
Then keeping the above properties of $\omega$ we can assume $\omega=F$ on $\partial(\sigma\times[0,\,1])$. \par
  Define $\xi$ as above. 
Now $\xi=\id$ on $\partial(\sigma\times[0,\,1])$. 
Apply the Alexander trick once more. 
Then we can extend $\xi$ to a definable homeomorphism $\Xi$ of $X\times[0,\,1]$ which is id outside of $\Int|\st(\sigma,K)|\times(0,\,1)$, preserving $\{\sigma'\times[0,\,1]:\sigma'\in K\}$ and level-preserving. 
Consider $F\circ\Xi^{-1}$ in place of $F$. 
It is a level-preserving definable $C^0$ imbedding of $(X_R,X_{iR})\times[0,\,1]$ to $(Y_R,Y_{iR})\times[0,\,1]$, $F'\circ\Xi^{-1}(\cdot,0)=f_R(\cdot)$, $F'\circ\Xi^{-1}(\cdot,1)=g_R(\cdot)$, $F'\circ\Xi^{-1}(x,t)=\phi_R(x)$ for $(x,t)\in X'_R\times[0,\,1]$, $F\circ\Xi^{-1}=F$ on $|K^k|\times[0,\,1]\cup(|K^{k+1}|-\sigma)\times[0,\,1])$ and $F\circ\Xi^{-1}$ is equal to $\omega$ on $\sigma\times[0,\,1]$ and, hence, PL. 
Thus $F\circ\Xi^{-1}$ is the required modification of $F$, which completes the proof of injectivity. \par
  {\it Injectivity of $\Aiso_\R(X,Y)\!\to\!\Aiso_R(X_R,Y_R),\,\Aiso_\R(X,Y)_\phi\!\to\!\Aiso_R(X_R,Y_R)_{\phi_R}$ and their 
relative maps.}
It suffices to prove the following statement. \par
  Set $Y'=\phi(X')$. 
Let $\pi_t:(Y_R,Y_{iR},Y'_R)\to(Y_R,Y_{iR},Y'_R),\ 0\le t\le 1$, be a definable isotopy of $Y_R$ such that $\pi_t=\id$ on $Y'_R$ for $t\in[0,\,1]$ and $\pi_1$ is the $R$-extension of a definable homeomorphism $\tau$ of $Y$. 
Then there exists a definable isotopy $\tau_t:(Y,Y_i,Y')\to(Y,Y_i,Y'),\ 0\le t\le 1$, of $Y$ such that $\tau_t=\id$ on $Y'$ for $t\in[0,\,1]$ and $\tau_1=\tau$. \par
  We can assume here $Y_i$ and $Y'$ are closed in $Y$. \par
  {\it Case of compact $Y$.} 
Assume $Y,Y_i$ and $Y'$ are polyhedra. 
Let $L$ be a simplicial decomposition of $Y$ compatible with $\{Y_i,Y'\}$. 
We have already seen in the proof of surjectivity of $\Iso(X,X_i;Y,Y_i)_\phi\to\Iso(X_R,X_{iR};Y_R,Y_{iR})_{\phi_R}$ that $\tau:(Y,Y_i,Y')\to(Y,Y_i,Y')$ is definably isotopic to some PL homeomorphism preserving $L$ and fixing $Y'$ through homeomorphisms. 
Hence we can assume $\tau$ is PL. 
Then what we prove is the following statement by lemma 1.2. \par
  Let $\pi_t:(Y_R,Y_{iR},Y'_R)\to(Y_R,Y_{iR},Y'_R),\ 0\le t\le 1$, be a definable isotopy of $Y_R$ fixing $Y'_R$ such that $\pi_1$ is PL. 
Then there exists a PL isotopy of $Y_R$ to $\pi_1$ fixing $Y'_R$. \par
  We proved this statement in the proof of injectivity of $\Iso(X,X_i;Y,Y_i)_\phi\to\Iso(X_R,X_{iR};Y_R,Y_{iR})_{\phi_R}$ without the condition that the PL isotopy is through homeomorphisms. 
This condition is clearly satisfied in that proof. 
Thus the compact case is proved. \par
  {\it Case of non-compact $Y$.} We can assume $\overline Y$ and $\overline Y-Y$ are polyhedra in $R^n$. 
Then there exist semialgebraic $C^0$ functions $\rho_1,...,\rho_k$ on $\overline Y$ such that $\rho^{-1}(0)=\overline Y-Y$ for any $i$ and the map $(\rho_1,...,\rho_k)|_Y:Y\to R^n$ is injective. 
Consider $(\rho_1,...,\rho_k)(Y)$ in place of $Y$. 
Then $\overline Y-Y=\{0\}$, and $\pi_t$ and $\tau$ are extended to a definable isotopy of $\overline Y_R$ fixing $Y'_R\cup\{0\}$ and a definable homeomorphism of $\overline Y$ fixing $Y'\cup\{0\}$ respectively. 
Thus we reduce the problem to the compact case. 
\qed\enddemo

\proclaim{Remark 3.4}
Consider always compact definable sets and definable maps between them. 
Then by lemma 3.3 singular (co)homology groups, homotopy groups, linking numbers, Whitehead groups, codordims groups, etc., do 
not depend on the choice of $R$. 
Moreover by lemmas 3.2 and 3.3 {\bf most} theorems and theories on differential topology hold over $R$, e.g., the Cairnes-Whitehead 
theorem on $C^r$ triangulations of $C^r$ manifolds, $r>0$, the Poincar\' e duality theorem, cobordism theory, the $h$-cobordism 
theorem, surgery theory. 
It also follows that a compact $\R$-PL manifold $X$ admits a $C^r$ manifold structure, $r>0$, if and only if $X_R$ admits a 
definable $C^r$ manifold structure. 
Some theorems on differentiable maps between differentiable manifolds cannot be, however, extended as follows. 
\endproclaim
  {\it Example.} Set $S^1=\{(x,y)\in\R^2:x^2+y^2=1\}$ and give to $\R$ an o-minimal structure so that $S^1$ and any notation of 
$S^1$ are definable, e.g., the o-minimal structure of restricted analytic functions. 
Let $R$ be a real closed field with an o-minimal structure which is a conservative extension of $\R$ and not equal to $\R$ (see 
[C$_2$] for existence). 
For $\epsilon\in\R$, let $\rho_\epsilon$ denote the rotation of $S^1$ with rotation number $\epsilon$. 
Then $\rho_\epsilon$ is not structurally stable, i.e., there exists a definable $C^\infty$ diffeomorphism $\rho$ of $S^1$ 
arbitrarily close to $\rho_\epsilon$ in the $C^\infty$ topology but not topologically conjugate to $\rho_\epsilon$. 
However $\rho_{\epsilon R}:S^1_R\to S^1_R$ is structurally stable for irrational $\epsilon$. 
(Here the homeomorphism of conjugation is not definable. 
If we admit only definable homeomorphisms, the theory of dynamical systems does not work because it requires infinite operations.) \par
  {\it Proof.} Regard $S^1$ as $\R/\Z$. 
Then $S^1_R=R/\Z$ and $\rho_{\epsilon R}(t)=t+\epsilon$ for $t\in R/\Z$. 
Let $\epsilon$ be irrational. 
Set $A=\{t\in R:nt<1\ \text{for}\ n\in\Z\}$. 
Let $\rho$ be a homeomorphism of $S^1_R$ so close to $\rho_{\epsilon R}$ that $\rho-\rho_{\epsilon R}$ has values in $A$. 
Then there exists a homeomorphism $\pi$ of $S^1_R$ such that $\pi\circ\rho_{\epsilon R}=\rho\circ\pi$ for the following reason. \par
  Let $p:R/\Z\to R/(\Z+\epsilon\Z)$ and $q:R/(\Z+\epsilon\Z)\to R/(\Z+\epsilon\Z+A)$ denote the projections, and let $\phi:R/(\Z+
\epsilon\Z)\to R/\Z$ and $\psi:R/(\Z+\epsilon\Z+A)\to R/(\Z+\epsilon\Z)$ be maps such that $p\circ\phi=\id$ and $q\circ\psi=\id$. 
We modify $\phi$ so that $\phi(s+a)=\phi(s)+a$ for $(s,a)\in(R/(\Z+\epsilon\Z))\times A$ as follows. 
Set 
$$
\tilde\phi(s)=\phi\circ\psi\circ q(s)+(s-\psi\circ q(s))\quad\text{for}\ s\in R/(\Z+\epsilon\Z). 
$$
Then $\tilde\phi$ is a map from $R/(\Z+\epsilon\Z)$ to $R/\Z$ with $p\circ\tilde\phi=\id$ because $s-\psi\circ q(s)\in A$, and 
$$
\gather
\tilde\phi(s+a)=\phi\circ\psi\circ q(s+a)+(s+a-\psi\circ q(s+a))=\qquad\qquad\qquad\qquad\quad\\
\phi\circ\psi\circ q(s)+(s+a-\psi\circ q(s))=\tilde\phi(s)+a\quad\text{for}\ (s,a)\in(R/(\Z+\epsilon\Z))\times A. 
\endgather
$$
Hence
$$
\tilde\phi\circ p(t+a)=\tilde\phi\circ p(t)+a\quad\text{for}\ (t,a)\in(R/(\Z+\epsilon\Z))\times A. \tag 1
$$
Using $\tilde\phi$ we will define $\pi$. 
For each $t\in R/\Z$ there exists uniquely $n(t)\in\Z$ such that $t-\tilde\phi\circ p(t)=n(t)\epsilon$. 
Note 
$$
\gather
t+a-\tilde\phi\circ p(t+a)\overset\text{by }(1)\to=t+a-\tilde\phi\circ p(t)-a=t-\tilde\phi\circ p(t)=n(t)\epsilon\tag 2\\
\qquad\qquad\qquad\qquad\qquad\qquad\qquad\qquad\text{for}\ (t,a)\in(R/\Z)\times A,\\
\rho(t+A)=t+\epsilon+A\quad\text{for}\ t\in R/\Z.\tag 3
\endgather
$$
Set $\pi(t)=\oversetbrace n(t)\to{\rho\circ\cdots\circ\rho}\circ\tilde\phi\circ p(t)$. 
Then 
$$
\gather
\pi\circ\rho_{\epsilon R}(t)=\pi(t+\epsilon)=\oversetbrace n(t)+1\to{\rho\circ\cdots\circ\rho}\circ\tilde\phi\circ p(t)=\rho\circ
\pi(t),\\
\pi(t+a)\overset\text{by }(2)\to=\oversetbrace n(t)\to{\rho\circ\cdots\circ\rho}\circ\tilde\phi\circ p(t+a)\overset\text{by }(1)
\to=\oversetbrace n(t)\to{\rho\circ\cdots\circ\rho}(\tilde\phi\circ p(t)+a),\\
\pi(t+A)=\oversetbrace n(t)\to{\rho\circ\cdots\circ\rho}(\tilde\phi\circ p(t)+A)\overset\text{by }(3)\to=\tilde\phi
\circ p(t)+n(t)\epsilon+A=t+A,
\endgather
$$
and, hence, $\pi|_{t+A}$ is a homeomorphism of $t+A$ for each $t\in R/\Z$. 
Thus $\pi$ is the required homeomorphism of $R/\Z$. 
\qed

\enddocument